\documentclass[14pt]{extarticle}
\usepackage[margin=1in]{geometry}
\usepackage{dsfont}
\usepackage{amsfonts}
\newtheorem{theorem}{Theorem}[section]
\newtheorem{pro}{Proposition}[section]
\newtheorem{lemma}{Lemma}[section]

\newtheorem{cor}{Corollary}[section]

\newcommand{\proof}[1]{\noindent{\it\bf Proof:#1\ }}
\newcommand{\QED}{\hfill$\Box$\medskip}

\begin{document}

\title{ Weakly Smooth  Structures  in Gromov-Witten Theory }
\author{Preliminary Version\\
\\
\\
Gang Liu }
\date{August 22,  2013}
\maketitle

\section{Introduction}
  In order to establish  Fredholm theory  on  stratified topological Banach manifolds in Gromov-Witten theory,  we have introduced   flat structures on such manifolds in [L4]. Such a structure is obtained from local  flat coordinate charts. The transformations between these charts  are only continuous in general. The purpose of this paper is two-fold: 
  firstly to show  that   on  the topological Banach manifolds and bundles appeared in GW theory, there are enough smooth functions and sections viewed in any admissible charts and trivializations; secondly to demonstrate 
 some finer aspects about the weakly smooth sections.
 As far as the Fredholm theory in [L4] is concerned,  the existence of sufficiently many  smooth functions and sections   makes    these topological Banach manifolds and bundles  behave as if they are the smooth ones. 
  The  motivation  of this work and [L4]
   is to overcome the analytic difficulty in Gromov-Witten and Floer type theories, the lack of differentiability of the orbit
   spaces of $L_k^p$ stable maps. These orbit spaces are quotient spaces of parametrized stable $L_k^p$-maps by the actions of the  reparametrization groups. For our purpose, they are the primary  examples of stratified topological  Banach manifolds. In GW and Floer type 
   theories,  each of them  comes along with a stratified topological  Banach bundle with a proper Fredholm section.
   As mentioned in [L4], the lack of differentiability shows up at   a few different levels: (A) Even there is only one stratum, the transition 
   functions between different coordinate charts are not smooth  due to the non-smoothness of the actions of the reparametrization group. But at least, the specified Fredholm section, the ${\bar {\partial}}_J$-section, is smooth viewed in any 
   chart due to the equivariancy of such a section. (B) When there are at least two strata,  the ${\bar {\partial}}_J$-section is  only stratified smooth. Moreover,  the pre-gluing joining different strata 
   introduces  coordinate charts near the "ends" of the 
   higher strata, which correspond   to the deformations and degenerations of the domains of the stable maps.  These charts together with the related bundles are fiberations
   over the  Deligne-Mumford moduli space of stable curves locally.  Even within a fixed stratum, locally   
   the ${\bar {\partial}}_J$-section should be considered as a family of smooth sections on the "central fiber" parametrized by the  Deligne-Mumford moduli space. 
   (C) 
   Unlike the usual stratified spaces, there is a pathological phenomenon in the stratum structure
 of the space of stable $L_k^p$-maps: each lower stratum viewed as end of a higher stratum appears to  have "more" dimensions than the higher stratum has. In other words, the family of the stable  $L_k^p$-maps parametrized by the gluing parameters is not "flat".
 Therefore,  it does not even make sense to talk about smoothness of the ${\bar {\partial}}_J$-section at a point in the lower stratum along the normal directions pointing to the  higher stratum simply because there is no local "product structure" near the end. 
   
   In this paper, we mainly  deal with difficulty (A) for  topological Banach manifolds and bundles in GW theory. This means that we only deal with  part of the ambient spaces in GW theory, which has only one stratum with trivial isotropy groups. Of course, a general ambient space is decomposed into its strata, and the theory of this paper is applicable to each stratum with a "fixed domain".
   However, in the general case,  each stratum  may have "ends"  defined by  the pre-gluing and moving double points. Some modifications are needed in order to apply the work  here   to these  cases.
   These modifications and   other cases and  difficulties for establishing a Fredholm theory on  stratified Banach manifolds     are treated in [L2] and [L3].  At the end of this section, we will outline the construction of the flat chart  so that readers can have a general idea  on  how the difficulties in (B) and (C) are resolved in [L2] and [L3].

   Large part of this paper is to use the ambient space of $L_k^p$-stable maps  as an example to illustrate on how to 
   use its weakly smooth structure  to establish the $C^{m_0}$-smoothness  of the
   moduli spaces  of perturbed  $J$-holomorphic maps. Here $m_0=k-\frac {2}{p}>1$ is the Sobolev  smoothness of an $L_k^p$-map.
   
As far as the smoothness of the moduli space is concerned, the main result  of this paper is the following  theorem.

\begin{theorem}   
Let $s:{\cal B}_{k, p}(A)\rightarrow {\cal L}_{k-1, p}$ be the  ${\bar{\partial }}_J $-section  of the bundle  ${\cal L}_{k-1, p}$ over the space of 
stable $L_k^p$-maps of class $A\in H_2(M, {\bf Z})$ from $\Sigma=S^2$ to a compact symplectic manifold $(M, \omega)$ 
with an ${\omega}$-compatible almost complex structure $J$, where
the fiber of ${\cal L}_{k-1, p}$ at $(f: \Sigma \rightarrow M)$ is $L_{k-1}^p(\Sigma , \Lambda^{0, 1}(f^*(TM))).$ Assume that
$s$ is proper in the sense that the  moduli space of $J$-holomorpic sphere of class $A$, ${\cal M}(J, A)=s^{-1}(0)$ is compact. Assume further  that
all isotropy groups are trivial and that the virtual dimension of  ${\cal M}(J, A)$ is less than $m_0$. Then there are generic small perturbations ${\nu}=\{\nu_i, i\in I\}$, which are compatible sections 
of the local bundles $ {\cal L}_i\rightarrow W_i$ of class
$C^{m_0}$  defined on the local uniformizers $W_i, i\in I$, such that the perturbed moduli space ${\cal M}^{\nu}(J, A)=\cup_{i\in I}(s+\nu_i)^{-1}(0)$ is a compact  manifold of class $C^{m_0}.$

\end{theorem}   

  Note that the main point of this theorem is   (i)  the degree of the smoothness of the moduli space, namely it has at least the same smoothness as the Sobolev differentiability of a generic element  in ambient space has; (ii) the smoothness is achieved by the coordinate
   transformations between the given "natural'' charts.

 In  view of (i) and (ii) above,   this is the best one can get presumably from the analytic set-up here. Of course,  this implies that if all geometric data are  of class $C^{\infty}$, we get  $C^r$-smoothness for the perturbed
  moduli space ${\cal M}^{\nu_r}(J, A)$ with $r$ arbitrarily   large. Using the comments below on regularity of perturbed 
  $J$-holomorphic maps,  it  is possible to show that in this case there is a common perturbation independent of $r$ so that the moduli space is in fact $C^{\infty}$-smooth. 
 On the other hand,   the theorem is still true even when  the geometric data is not  $C^{\infty}$-smooth, but only sufficient smooth comparing to $m_0$. The exact degree of the required smoothness for the geometric data can be determined from the proof of the theorem.
 Presumably $2k$ should be sufficient.
   
   The existence of such a  smooth or stratified smooth structure  on the moduli space
   have many  applications in $GW$ and Floer type theories. Note  that if we only want some weaker results such as the existence of a topological manifold structure on the perturbed  moduli space,  it is possible to avoid  part of the discussion on the finer aspect of weakly smoothness  in this section.  In fact, this    topological manifold structure on
   the  moduli space  is   already sufficient for   most of applications. For instance,  it  was used in 
   [LT] to prove the Arnold conjecture. 
    However, when [LT] was written, the issue of lacking smoothness was not yet addressed.  Nevertheless,  the matter was resolved  by  the author during the preparation  for his   seminar courses. 
   We now describe   what  modifications  to the global perturbation method of [LT] are needed to   get  a topological manifold structure on extended moduli spaces.  We will  mainly restrict ourself to the above case as it is well-known that this basic case already captures the main issue of the lack of differentiability. The main steps for generalization  will be only  outlined.
    Since the ${\bar{\partial}}$-section itself is smooth ( stratified smooth for the general case) viewed in any local uniformizer (which is called a local slice in this paper), we only need to make sure that the perturbation coming from the cokernel defined on one local slice
     is  at least $C^1$-smooth viewed in any other local slice. Those local perturbation are obtained by extending  the elements $\xi$ of the cokernel $K_f$ of linearization of the  ${\bar{\partial}}$-section at a $J$-holomorphic map $f$.

      The extension is obtained by extending  $\xi$ to a "constant" section ${\tilde \xi}$ over a local slice $W_f$ first by using the  trivialization of the local bundle ${\cal L}(f)\rightarrow W_f$ induced  by the   parallel transport  of $M$, then multiplying ${\tilde \xi}$ by a cut-off function supported on  $W_f.$

   Therefore, it is sufficient to show that  both cut-off function and  ${\tilde \xi}$  are  at least  of class $C^1$ viewed in any slice. 
 Any such functions or sections are called weakly smooth ones of class $C^1$ in this paper.
 
 \medskip
\vspace{2mm}
\noindent ${\bullet}$    Existence of the $C^{m_0}$-Smooth Cut-off Functions:

  \vspace{2mm}
\noindent
It is well-known   that when $p$ is a positive even integer, the $p$-th power of the $L_k^p$-norm of the Banach space
$L_k^p(\Sigma,  f^*TM)$ is smooth.
 This  smoothness of  $L_k^p$-norm of the Banach space and the related cut-off functions were already used in [LT].  In order to get the desired cut-off function on the 
 space of unparametrize stable $L_k^p$-maps, we need to show that the above $p$-th power of norm function is still smooth after composed  with the action map by  the Lie group $G$ of the reparametrizations of the domain. Denote the composed function by $ \Psi_p$.
 In the case that 
  $\Sigma={\bf R^1}\times S^1$ and $G$ is the group of translations, we have $ \Psi_p\circ g=\Psi_p$ for any $g\in G$.  Hence  $ \Psi_p$ is
$G$-invariant despite of the fact that $G$-action is only continuous generically. Because  the simplest case  above  is supposed to  capture the main difficulty of lack of differentiability already, this immediately suggests that in the general case 
 $\Psi_p$ should be  smooth  as well  and how this can be proved. Indeed the only difference of this special case with the general case is that $G$-action here preserves the volume form on $\Sigma,$   and the smoothness in $G$-direction for the general case can be  proved by a simple variable change formula in calculus. For completeness, we include the proofs of the above well-known  facts in the last section of this paper.

Since open sets in $L_k^p(\Sigma,  f^*TM)$ serve as 
coordinate charts of the space of parametrized  stable maps,  the smoothness of the  $L_k^p$-norm with respect to the $G$-actions should imply  that the push-forward of this function  to a local slice should be  $C^{m_0}$-smooth  viewed in  any other slice. 
Whether or not  this is true is not completely  obvious. 

We proceed slightly differently. First embed  the symplectic manifold $M$ with the induced metric isometrically into some ${\bf R}^d$. Then   consider the space of parametrized stable $L_k^p$-maps in $M$ as a closed Banach submanifold of the corresponding space of the maps
in ${\bf R}^d$. This later space is  a  Banach space so that the result above is applicable. 
Note that here we have used the fact that when $m_0$ is large enough, the  two Sobolev metrics on the space of $L_k^p$-maps from $\Sigma$ to $M$ are equivalent. 
 As a result of this method, in addition to the cut-off function that we are looking for, we also obtain a large collection of $G$-smooth functions on the space of parametrized  $L_k^p$-maps in $M$ by using pulling-backs.  They give rise  the corresponding weakly $C^{m_0}$-smooth functions defined on the space of unparametrized stable $L_k^p$-maps  of  $M.$

 Note that in general,   even above composed function  $\Psi_p$  is of class $C^{\infty}$ and  the cut-off function is of class $C^{\infty}$  in a given slice,  we still can not conclude that the cut-off function is $C^{\infty}$-smooth viewed in any slices (see the proof in  section 3 and the discussion below for the reason of this). In other words, it is weakly $C^{m_0}$-smooth but not  weakly $C^{\infty}$-smooth. Therefore we get the desired  cut-off functions.

\medskip
\vspace{2mm}
\noindent ${\bullet}$   $C^1$-Smooth  Perturbations: 
 
 As for the weakly smoothness  of ${\tilde \xi}$,   without using the discussion below on some finer aspect
 of weakly smooth sections, we can only prove that it is of class $C^1$.  
 Of course, there is no need to assume that $\xi $ is in the cokernel. But we do need to assume that it is an element of the fiber of $({\cal L}_{k-1, p}(f))_f$ of class at least $L_k^p$. Note that any element in $W_f$ is at least of class $L_k^p$.   Assume that $f$ and $\xi$ are of class $C^{\infty}$ first.

  Let ${\tilde {\cal B}}$ be the space of parametrized stable $L_k^p$-maps and  ${\tilde {\cal L}}\rightarrow {\tilde {\cal B}}$ be the 
  corresponding bundles.  Denote a small neighbourhood of $f$  in ${\tilde {\cal B}}$  by ${\tilde W}_f$ and the corresponding local bundle 
  by ${\tilde {\cal L}}(f).$   We have already denoted its restriction to  the slice  $W_f$ by ${\cal L}(f)$.   Assume that ${\tilde W}_f$ is 
  part of the $G$-orbit of $W_f$ near $f$.
 Then by abusing the notation,  we have a homeomorphism  ${\tilde W}_f \sim G\cdot W_f$ of two open  sets in a Banach space, and  two trivializations
 of  ${\tilde {\cal L}}$ on ${\tilde W}_f$: one is given by  using parallel transport of the central fiber; the other is obtained
 by using  parallel transport of the central fiber over the slice $ W_f$ first, then using the pull-backs of the $G$-actions to bring the 
 fibers over $ G\cdot W_f.$ The two trivializations are only topological equivalent. Use the first trivialization, the "standard" one, we 
 get a "constant"
   section over ${\tilde W}_f $ from $\xi$ in the central fiber, denoted by ${\bar \xi}_1$, which is smooth but not $G$-equivariant. To get a  $G$-equivariant extension,  we use the second trivialization to extend ${\tilde \xi}$,  the  restriction  of   ${\bar \xi}_1$  to  the 
  slice ${ W}_f $, over
 ${\tilde W}_f \sim G\cdot W_f$.   Clearly, this latter  extended  section, denote by  ${\bar \xi}_2$ is smooth with respect to the second 
 trivialization. The question is  about the degree of the smoothness of the section ${\tilde {\Phi}}_{1, 2}\circ {\bar \xi}_2$,  where
 ${\tilde \Phi}_{1, 2}$ is the topological identification  from the  second trivialization to the first one. In other words, we want to know the smoothness of ${\bar \xi}_2$ viewed in the standard trivialization.
 
 To this end,  we  first sketch a argument using $sc$-smoothness in the nested Sobolev spaces $L_1^p\supset L_2^p\cdots \supset L_k^p \cdots $  in the "ordinary sense" rather than the $sc$-smoothness in polyfold theory.  In other words,  instead of using the weaker topology to define derivatives, the usual operator norm is used for the definition of derivatives.
 Note that by letting $k$ varying, all identifications above  are of class $sc^{\infty}$. 
 Therefore,  ${\tilde {\Phi}}_{1, 2}\circ {\bar \xi}_2$ gives rise a $sc^{\infty}$-section since ${\bar \xi}_2$ is such a section. Under the trivializations, above sections become
 maps from the two domains to the central fiber, $L_{k-1}^p(\Sigma, \Lambda^{0, 1}(f^*(TM)))$, which give rise the corresponding
  $sc^{\infty}$-maps, denoted by the same notations. The  assumption that  $\xi$ is of class $C^{\infty}$ , hence in
   $L_{k}^p(\Sigma, \Lambda^{0, 1}(f^*(TM)))$ for any $k$ implies 
  that at level 
  $k$,  ${\bar \xi}_2:G\cdot W_f\rightarrow L_{k-1}^p(\Sigma, \Lambda^{0, 1}(f^*(TM)))$
   factors through
    ${\bar \xi}_2:G\cdot W_f\rightarrow L_{k}^p(\Sigma, \Lambda^{0, 1}(f^*(TM))). $ 
  So  does  ${\tilde {\Phi}}_{1, 2}\circ {\bar \xi}_2:{\tilde W}_f\rightarrow L_{k}^p(\Sigma, \Lambda^{0, 1}(f^*(TM)))\subset L_{k-1}^p(\Sigma, \Lambda^{0, 1}(f^*(TM))).$  Here we have used the fact that ${\tilde W}_f \sim G\cdot W_f$ is a $sc^{\infty}$-equivalent.  Note that while the  map $ {\bar \xi}_2$ can be lifted further into higher Sobolev space by our assumption,  ${\tilde {\Phi}}_{1, 2}\circ {\bar \xi}_2$ can not be.  Roughly speaking the reason for this  is  that the $sc$-regularity for the transition map between the two trivializations is of degree $k$=the degree of a generic element of the above two domains.

  It follows from the 
  definition of $sc$-smoothness, in this situation, ${\tilde \Phi}_{1, 2}\circ {\bar \xi}_2$ is of class $C^1$. 
  In other words, the $G$-equivariant extension ${\bar \xi}_2$ of ${\xi}$ is indeed of class  $C^1$ viewed in the standard chart and trivialization of the local bundle  ${\tilde {\cal L}}$ on ${\tilde W}_f$.
  Consequently, ${\tilde \xi}$ is of 
  class $C^1$ viewed in any slices (weakly smooth of class $C^1$).
   Having  faith in  that the smoothness obtained from  the considerations  in  $sc$-smoothness   should be  optimal, one may conclude that generically
    this $C^1$-smoothness is the best one can hope. Indeed, the following more direct and elementary argument obtained before the theory of the $sc$-smoothness  
     available  leads  to  the same conclusion.

To start the  direct proof,  denote  ${\tilde \Phi}_{1, 2}\circ {\bar \xi}_2$
  by $[{\bar \xi}_1]$ for short.
 We start with the two Banach coordinates on ${\tilde W}_f \sim G\cdot W_f$.   We have the standard coordinate given by $(w,t)$ with $w$ in the local slice $W_f$ and $t$ in $\hat {H}^*$, where $H^*$ is the complement  space of  $W_f$ in ${\tilde W}_f$ which  can be identified with the    direct sum of the orthogonal complements, denoted by ${\hat H}^*$ of the corresponding local hypersurfaces used to define the local slice  above (see Sec. 3 or [LT] for the  more  details). Note that dimension of $H^*$ is equal to dimension of $G$.      The second coordinate from  the above identification  is denoted by $(u, g)$ with $u$ in $W_f$ and $g$ in $G_e$. Here
 $G_e$ is a small neighbourhood of $e$ in $G$.  Note that according to Sec. 3, the coordinate $(w, t)$  for ${\tilde W}_f=W_f\times H^*$ with respect to the splitting 
  is $C^{\infty}$-compatible with the "standard" exponential coordinate for  ${\tilde W}_f $.  In particular,  different choices of the "base 
  point'' $h\in W_f$ (=origin), splitting factor $H^*$ give $C^{\infty}$-equivalent coordinates and related trivializations. Moreover since  for 
  a fixed $g\in G$,  its action  $\Psi_g$ is a smooth differomorphism  from 
  ${\tilde W}_f$ to its image and the action lifts to a smooth isomorphism between the two local bundles via pull-backs. It follows from 
  these considerations together with the fact that ${\bar \xi}_2$ is  $G$-equivariant  implies that  we only need to compute $D [{\bar 
  \xi}_1]_h$ with $h\in W_f$ using  the standard coordinate $(w, t)$ with $h$ as base point. Since 
   $[{\bar \xi}_1]|_{W_f}={\tilde \xi}$ is a "constant" section with respect to the trivialization and hence smooth, we have that along $w$-
   direction,  $D[{\bar \xi}_1]_h$ exists and equals to zero. 
  
To compute the partial derivatives along $H^*$ or $t$-direction,  consider a point $h^*\in H^*$ with coordinate $(w, t)=(0, {\hat h}^*).$
Then there is a unique $g=g(h^*)\in G$, such that $h^*\circ g(h^*)$ is in $W_f$. Hence the $(u, g)$-coordinate of $h^*$ is $(u,  g)=
(h^*\circ g(h^*),  (g(h^*))^{-1}).$ We now show that $h^*\rightarrow g(h^*)$ as a map form $H^*$ to $G$ is of class $C^{m_0}$. To this end, note that 
the identification $H^*\rightarrow {\hat H}^*$  is the restriction to $H^*$ of the $C^{\infty}$-smooth map $\pi_{{\tilde H}^*}\circ 
ev^l_{x}:{\tilde W}_f\rightarrow {\hat H}^*$. Here $ ev^l_{x}:{\tilde W}_f \rightarrow M^l$ is the $l$-ford evaluation map used to define the local 
slice $W_f$ with $2l=\dim (G)$, $\pi_{{\tilde H}^*}$ is the projection from a neighbourhood of $h(x)$ in $M^l$ to ${\hat H}^*$.
By definition, $(\pi_{{\tilde H}^*}\circ ev^l_{x})(h^*)$ is just the $t$-coordinate of $h^*$, and it has been denoted by ${\hat h}^*.$
 Now we have the local identification of class $C^{m_0}$ given by $h:D^{l}(x)\rightarrow {\hat H}^* \subset M^l.$ Here $D^{l}(x)$ is the $l$-fold 
 product of small discs centred  at the (minimal number of  ) marked points of the free components. Note that $h(x)$ is the ''origin'' of ${\hat H}^*$.
 Therefore ${\hat h}^*\rightarrow h^{-1}( {\hat h}^*)$ as a local map from ${\hat H}^*$ to $D^{l}(x)$ is of class $C^{m_0}.$  Clearly 
 $g(h^*)$ is determined by the location of  $h^{-1}( {\hat h}^*)$ in $D^{l}(x)$, and smoothly depends on it.  Put this together, we have 
 proved that $ g(h^*)$ depends on $h^*$ with $C^{m_0}$-smoothness. Now let $h^*=h_s$ be a smooth curve in $H^*$ such that $h_0=h\in W_f.$
 Then it has $(w(s), t(s))$-coordinate $(0, {\hat h_s}^*)$ as smooth functions of $s$.  The corresponding $(u, g)$-coordinate is  $(u(s), 
 g(s))=(h_s\circ  g(h_s), (g(h_s))^{-1}).$ It is of class $C^{m_0}$ in $s$.  Here we have abused notation by writing $h_s$ for its coordinate.   
 Let $v=\frac{ \partial h_s }{\partial s }|_{s=0}$. Then the partial derivative along $H^*$ is given by
 
 $$D[{\bar \xi}_1]_h(v)=\frac{ \partial[{\bar \xi}_1](u(s), g(s))}{\partial s }|_{s=0}$$
 $$=\lim_{s\mapsto 0}\frac{( [{\bar \xi}_1](u(s), g(s))-[{\bar \xi}_1](u(s), g(0)))}{s}$$
 $$+\lim_{s\mapsto 0}\frac{( [{\bar \xi}_1](u(s), g(0))-[{\bar \xi}_1](u(0), g(0)))}{s}.$$
 
 Since $[{\bar \xi}_1]$ is a constant along $W_f$, the second limit is equal to zero. Assume that $[{\bar \xi}_1](u, g)$ is of class $C^1$ 
 as the function on  $(u, g). $ Then the first limit is equal to $$\lim_{s\mapsto 0}\frac{ \partial[{\bar \xi}_1](u, g)}{\partial g }(u(s), g(\xi))\cdot 
 \frac{ \partial g}{\partial s}(\xi)=\frac{ \partial[{\bar \xi}_1](u, g)}{\partial g }(u(0), g(0))\cdot 
 \frac{ \partial g}{\partial s}(0). $$ Here $0<\xi<s.$ 
 
Therefore,  we only need to  show existence of the partial derivatives of $[{\bar \xi}_1]$ as a continuous function  in $(u, g)$.  
We will only prove the existence of the partial derivatives and  leave it to the readers  to prove its continuity.

For the same reason as above,  we only need to consider the case that $h$ is in $W_f$. In this case the derivatives along $u$-directions  are  still equal to zero since $[{\bar \xi}]
 $ is a constant on $W_f$.

Therefore, we only need to show that $(\frac{ \partial [{\bar \xi}_1]}{\partial g})_h$  exists.

To this end, assume that $\frac{ \partial }{\partial g}|_h=\frac{ \partial h_s }{\partial s }|_{s=0}$. Here $h_s=h\circ g_s:\Sigma\rightarrow M$  with $s\in (-\epsilon, \epsilon)$   is  a  curve in ${\tilde W}_f$ with $h_0=h$ representing $\frac{ \partial }{\partial g}|_h$ and defined by the corresponding curve $g_s$ in $G$.

Then we have $(\frac{ \partial [{\bar \xi}_1]}{\partial g})_h=(\frac{ \partial [{\bar \xi}_1](h\circ g_s)}{\partial s})_{s=0}$.
Now being considered  as s-dependent sections on the fixed bundle $\Lambda^{0, 1}(f^*TM)\rightarrow \Sigma, $  $$[{\bar \xi}_1](h_s)=[{\bar \xi}_1](h\circ g_s)=
\Pi^{-1}_{h_s}(((g_s)^{-1})^*(\Pi_{h} (\xi))) .$$   Here these $\Pi_{h_s}$ are bundle morphisms from $\Lambda^{0, 1}(\Sigma,f^*TM)\rightarrow \Sigma  $ to $$\Lambda^{0, 1}(\Sigma, h_s^*TM)\rightarrow \Sigma. $$ Since $h_s$ is of class $L_k^p$,  so is the bundle $\Lambda^{0, 1}(\Sigma, h_s^*TM)\rightarrow \Sigma   $. Therefore these differomorphisms are of class  $L_k^p$ too. Note  that the section   $\Pi_{h} (\xi)$ of the bundle $\Lambda^{0, 1}(h^*TM)\rightarrow \Sigma   $ can only be of class $L_k^p$ even $\xi$ is of class $C^{\infty}$. From this explicit formula for $[{\bar \xi}_1](h_s)$, one concludes that its derivative with respect to $s$ is  of class $L_{k-1}^p$. This proves the existence  of  $D[{\bar \xi}_1]_h$ at any point $h$. Again we leave it to the readers  to prove its continuity  with respect to $h$ so that $[{\bar \xi}_1]$ is of class $C^1.$ 

This concludes the necessary modifications to [LT] for the case there is only one stratum. The general case is not really much harder. In 
this case, we want to show that the same is true in each stratum. As mentioned above,  main  point  here is to note that near the end of a 
higher stratum,   a local slice $W_f$ is a fiberation  over a local chart (uniformizer) of ${\bar {\cal M}}_{0, k}$ with $\Sigma_f$  in  a 
lowest stratum of ${\bar {\cal M}}_{0, k}$. Hence within a stratum, it splits locally as a product of $W^{\alpha_0}_f$ and a neighbourhood  
$\Lambda_{\delta}(\alpha_0)$ of $[\Sigma_{\alpha_0}]$ in the corresponding stratum of ${\bar {\cal M}}_{0, k}.$ Here $\alpha$ is the local 
parameter of ${\bar {\cal M}}_{0, k}$ near $\alpha_0$ and $\Sigma_{\alpha_0}=\Sigma_{f_{\alpha_0}}, $  where $f_{\alpha_0}$ is one of the "base" points in the given stratum obtained from $f$ by moving its double points first, then make the corresponding pre-gluing.
There is also a corresponding (only $C^{\infty}$) product structure for 
the universal curve ${\cal U}|_{\Lambda_{\delta}(\alpha_0)}\simeq \Sigma_{\alpha_0}\times  \Lambda_{\delta}(\alpha_0)$ so that we may view   $\Sigma_{\alpha} $ with $\alpha\in \Lambda_{\delta}(\alpha_0)$    as a family of complex structures defined on  the  fixed 
$ \Sigma_{\alpha_0}$ parametrized by $\alpha.$ Using this product 
structure, the ${\bar {\partial }}_J$-section on $W_f$ is translated into a family of sections ${\bar {\partial }}_J^{\alpha}$ defined on  
the  fixed fiber  $W^{\alpha_0}_f$. Similar interpretations  are applicable for those perturbations. The metric on the domains used to 
define the space of $L_k^p$ maps and related bundles become a family of metrics on a fixed domain. Moreover, the action of the  automorphism 
group acting on the free components of $\Sigma_f$ extends naturally to the product $\Sigma_{\alpha_0}\times  \Lambda_{\delta}(\alpha_0)$, 
which, in turn, induces an action on $W^{\alpha_0}_f\times \Lambda_{\delta}(\alpha_0)$. This essentially put us in similar situation as 
above, and what one needs to do is the show the  corresponding  statements  accordingly. At this point, the proofs of these statements are the straight forward generalizations of what we have  done  above.  The details will be given in a subsequent paper.

 Therefore,  with the  supplement here, for the case that all isotropy groups  being  trivial,  what proved by  the argument of [LT]  is the following theorem.

\begin{theorem}

Let $s$ be the ${\bar {\partial }}_J$-section.
Assume that all isotropy groups are trivial.  Let $B(r)=\oplus_{i\in I}B(r_i)$  be the ball of "radius" $r$ in ${\cal R}=\oplus_{i\in I}{\cal R}_i$=collection of all compatible perturbations ${\nu}=\{\nu_i, i\in I\}$.  Here  $\nu_i\in r_i$   is a stratified weakly $C^1$-section supported on  the local uniformizers $W_i, i\in I$ obtained from the cokernel as above with $\|\nu_i\|<B(r_i)$, and $\{W_i, i\in I\}$ covers ${\cal M}(A, J)$ in ${\cal B}.$
Then for $r$  small enough, each  extended  local moduli space $E{\cal M}^{i}(J, A)=(s+ev_{(i,B(r))} )^{-1}(0)$  is a stratified submanifold of class $C^{1}$ in ${\cal R}\times W_i.$ 
In fact elliptic regularity implies that above local moduli space $E{\cal M}^{i}(J, A)=(s+ev_{(i,B(r))} )^{-1}(0)$  is of  class $C^{1}$ in ${\cal R}\times W^{\infty}_i$  if all the geometric data are  of class $C^{\infty}$.  Here $W^{\infty}_i\subset W_i$ consists of  the  corresponding $C^{\infty}$ elements.
Moreover, these local extended  moduli spaces patch together to form a   stratified topological manifold  $ E{\cal M}(J, A)=\cup_{i\in I}E{\cal M}^{i}(J, A)$, the total space of the obstruction sheaf.

\end{theorem}

See the discuss later in this section on how to  obtain  the perturbed moduli space from  the extended moduli space.
This concludes the  proof of the following theorem  and   fills in the gap  in [LT].

 \begin{theorem}
  With the supplement above, the Floer homology  in [LT] is well-defined. Similarly using the method in [L0], the work in [LT] 
  establishes the existence of $GW$-invariants and quantum cohomology for a general symplectic manifold.

  \end{theorem}

 \medskip
\vspace{2mm}
\noindent

 \medskip
\vspace{2mm}
\noindent ${\bullet}$ "Geometric"  $C^{m_0}$-Perturbation:

As mentioned before, the best we can get for extending an element  $\xi$ in the central fiber ${\cal L}(f)$ to a constant section ${\tilde \xi}$ over $W_f$ is to get a weakly smooth section of class $C^1$ even we start with a  $C^{\infty} $ element $\xi$  like ones in the cokernel.

In order to get the desired $C^{m_0}$ smoothness for perturbations used in the main theorem of this paper, we proceed differently.

  Note that  by linearity  we only need to extend each element in a basis of the cokernel (or any prescribed finite dimensional space of the central fiber).
 
The  key observation then is that instead of extending each element in the basis using the standard process above, we decompose each element $\eta$
 into a finite sum of elements, each of them is localised near a point of the domain $\Sigma$. Of course, these new elements are not in the cokernel anymore. In the case that $f$ and $\eta$ are smooth, each of these new elements localized in a disc $D_{\delta}(x_0)\subset\Sigma$ can be written as a finite sum the terms of the form $\phi\otimes_{(i,J)}\xi$. Here   $\phi$ is a smooth $(0,1)$-form on $\Sigma$ supported in  $D_{\delta}(x_0)$ and $\xi$ is a restriction of a smooth global vector field ${\tilde \xi}$ on $M$.
Therefore, we only need to extend the elements of the above form. This can be done for $\phi$ and $\xi$ separately in the corresponding bundles. More specifically,  let ${\tilde {\cal L}}\rightarrow {\tilde {\cal B}}$ be the bundle used to define ${\bar{\partial}}$-section for parametrized stable $L_k^p$-maps.  Then the  fibre of ${\tilde {\cal L}}={\tilde{\cal L}}_{k-1,p}$ at $f$,  $L^p_{k-1}(\Sigma, \wedge^{0,1}(\Sigma) \otimes f^*(TM)),$ is linearly
homeomorphic to  $L^p_{k-1}(\Sigma,  f^*(TM)) \otimes_{L^p_{k-1}(\Sigma)}L_{k-1}^p(\wedge^{0,1}(\Sigma)) $ when $k$ is large enough.
 This gives rise a bundle isomorphism ${\tilde{\cal L}}\simeq {\tilde{\cal T}}\otimes {\tilde \Omega}^{1}.$
 Here   the fibre of ${\tilde{\cal T}}={\tilde{\cal T}}_{k-1,p}$ at $f$ is $ L^p_{k-1}(\Sigma,  f^*(TM)), $  and ${\tilde \Omega}^{1}$ is the  trivial bundle whose fiber is $L_{k-1}^p(\wedge^{0,1}(\Sigma)) .$ Note that the  $G$-actions on ${\tilde{\cal T}}$ and ${\tilde \Omega}^{1}$ are compatible with the one on ${\tilde{\cal L}}.$ So are the local trivializations and $G$-equivariant local trivializations induced by the ones on  a local slice. This implies that we really can deal with the two components of $\eta$ separately.

Clearly, the vector field ${\tilde \xi}$ induces a global section on  ${\tilde{\cal T}} $, denoted by ${\bf \xi}_0$ and defined by  ${\bf \xi}_0
(g)={\tilde \xi}(g)$ for any $g \in {\tilde {\cal B}}$. Since  $G$ only acts on the domain,   the section ${\bf \xi}_0$ is clearly $G$-equivariant. It is easy to see that it is also smooth (see Sec. 4 for  the proof). This completes the desired extension for $\xi$.

To extend $\phi$, use the trivialization ${\tilde \Omega}^{1}\simeq {\tilde {\cal B}}\times L_{k-1}^p(\wedge^{0,1}(\Sigma)).$  We get a smooth 
global section ${\tilde {\phi}}$ from $\phi$ and   its restriction to a local slice $W_1(f)$, denoted by  ${\tilde {\phi}}_1$. 
As above, the question is if it is still smooth viewed in the other slices.   We follow the same idea,  using  the  $G$-action to  obtain a $G$-equivariant section, denoted by  ${\bf {\phi}}_1$
over the open set $G\times W_1(f),$ then   deciding  if  ${\bf {\phi}}_1$ is smooth on  $G\times W_1(f).$
Since the situation here is much better than the general extension problem we discussed before, in stead of getting $C^{1}$-smoothness, the similar computation there shows that ${\bf {\phi}}_1$ is of class $C^{m_0}. $ 
In Sec.3, we give a direct proof that the extension  ${\tilde {\phi}}_1$
 on  $W_1(f)$ is $C^{m_0}$-smooth viewed in any other local slices.

 Put this together, we get  the desired extensions
of the elements in the coknernel over a local slice, which are  $C^{m_0}$-smooth viewed in any local slices. We will show in Sec. 3 that
these local extensions  can be used to achieve local transversality. Since these extensions  are  of   geometrical nature,    we will refer the sections and perturbations so obtained as geometric sections and  perturbations though they are not the usual
 geometric  perturbations  in the sense  used in GW theory.

 We have already known that the "constant" extension induced by the standard parallel transport can only be of class $C^1$. If we allow the parallel transport used here depending on $f$,  then     the "constant extension"  of a smooth   element in the central fiber over a local slice is  still $C^{m_0}$-smooth viewed in any slices.

  \medskip
\vspace{2mm}
\noindent ${\bullet}$ ${\bullet}$  $C^{m_0}$-smoothness of   constant extensions defined by $f$-dependent connections:

\medskip
\vspace{2mm}
\noindent
Assume that a  point-section $\eta\in C^{\infty} (\Sigma, \Lambda^{0, 1}(f^*(TM)))\subset  L_{k-1}^p(\Sigma, \Lambda^{0, 1}(f^*(TM)))$ over $f:\Sigma \rightarrow M$ satisfies the condition that 
 (i) $\eta$  is supported in an open disc $D_{\delta}(x_0)$ in $\Sigma$; (ii) $f$ is an embedding on $D_{\delta}(x_0)$; (iii) over  $D_{\delta}(x_0)$, $\eta$ can be written as $\eta =\phi \otimes_{(i, J)} \xi$, where $\phi $ is a smooth $(0,1)$-form on $\Sigma$ supported in  $D_{\delta}(x_0)$ and $\xi$ is  the restriction of a smooth vector field ${\tilde \xi}$ on an open set $U$ of $M$ containing $f(D_{\delta}(x_0))$ with the property that ${\tilde \xi}$ is 
a covariant constant vector field  with respect to a $J$-invariant connection on $U$.
 
  The similar discussion as above (see  Sec. 3 ) 
  shows that in this case the extension of $\eta$ to a local slice containing $f$ by parallel transport is  $C^{m_0}$-smooth viewed in any other slices. 
  
  On the other  hand, one can show that by using a local parallel frame of the complex bundle $(TM, J)$ over $U$,  the condition (i) and (ii) implies that the condition (iii) can be arranged. 
  Of course, for a symplectic manifold $(M, \omega)$ with $\omega$-compatible
  almost complex structure $J$, the standard  $(J,\omega)$-invariant connection 
  with Neijinhaus tensor as its torsion  does not have  such a local $J$-flat frame    unless  $J$  is integrable over $U$.   However, if we regard $(TM, J)$ as a abstract complex vector bundle rather than as the tangent bundle of $M$, any local complex trivialization of $TM$ over $U$  gives rise    such a local $J$-flat frame  over $U$. This means that we have to use a family of $(f, x_0)$-dependent connections on $M$ to give the local trivializations of  the bundle $({\cal L}\rightarrow {\cal B})$ at least for $ f$ in the moduli space of $J$-holomorphic maps.

  We conclude that if $f$ satisfies the condition (i) and (ii), by using a $ (f, x_0)$-dependent connection on   $M$,     the constant extension of $\eta$  with respect to the connection is   $C^{m_0}$-smooth. 
  
 Note that the point-sections to be extended in [LT] are obtained from the cokernel by multiplying all  elements in the cokernel by a fixed cut-off function supported in an open set of $\Sigma$ away from its double points (see page 29 of [LT], under the assumption that the evaluation map at double points is transversal to the multi-diagonal there at  the  $J$-holomorphic map $f$ ).
  By the unique continuation principle  used in [LT],   the same argument there implies  that the cut-off  function  can be chosen  to be  supported in a disc  $D_{\delta}(x_0)$ for each component of the domain $\Sigma$.
  Since $f$ in this case is in the moduli space of $J$-holomorphic curves, it is automatically smooth and has some point $x_0$ on each component of $\Sigma$ such that $f$ is a local embedding near $x_0$. In fact such points are open and sense in $\Sigma$. Therefore as long as we choose the center of each  $D_{\delta}(x_0)$ to be one of the above "good" points,
 both (i) and (ii) are automatically true.   Consequently, with above modification, the extension used in [LT] is $C^{m_0}$-smooth viewed in any slices if we are willing to use a family of $f$-dependent connections to obtain "constant" extensions. Note that since  the moduli space is compact, for the purpose of the constructions in [LT], there are only finitely many such connections are involved.


\medskip
\vspace{2mm}
\noindent ${\bullet}$ Regularity Assumption:

 We have mentioned that  it is possible to establish the $C^{m_0}$-smoothness of the perturbed moduli space without assuming that all geometric data are $C^{\infty}$-smooth  but only assuming that they are only, for instance, $C^{2k}$-smooth.
 
 This seemly very technical   point  concerns the  general  philosophy on how to deal with the main  difficulty of lack of differentiability in   the current research. The reader might   have been  aware of that  in the polyfold theory,  the Frechet manifold of smooth stable maps lying inside  the sc Banach manifold of stable maps plays a fundamental role in the formulations of various notions and constructions of the theory.     The basic requirement that the resulting moduli space should lie inside the  Frechet manifold  is behind many considerations in the theory. 
 
 This motivates our efforts in this and subsequent  papers to explore the   possibilities: (i) to get the main construction of virtual moduli space without using the $C^{\infty}$ regularity results  in the case of geometric data are smooth; or (ii) to get the same conclusion even
  without assuming the geometric data are smooth.    
  As  the proofs in this paper show, under the assumption that all geometric data are $C^{2k}$-smooth, one can eliminate the role played by the  Frechet manifold and still obtain the extended moduli space, the total space of the obstruction bundle (sheaf),  with  $C^{m_0}$-smoothness,  the same regularity as its elements. 
  
  As this stage,  to  get  the perturbed moduli space, one has to use Smale-Sard theorem.  Therefore one may simply  require that the virtual dimension of the moduli space is less than $m_0$. The other possible way  to deal with this last problem is to use Whitney's theorem to give the extended moduli space a compatible $C^{\infty}$-smooth  structure first, then  to   deform the "projection" map from the extended moduli space to the linear space  of the perturbations into a smooth one. This last method works without the assumption on the virtual dimension, but it requires the true $C^{m_0}$-smoothness of the extended moduli spaces (in [L3])rather than just the stratified smoothness in the general case.

 There is  a related
 question about the  regularity  of individual perturbed $J$-homomorphic maps. This was used  from the main theorem about $C^{m_0}$-smoothness of the moduli space  to infer its $C^{\infty}$-smoothness under the assumption that geometric data are smooth.
 It has been suspected that the solution of a perturbed ${\bar {\partial }}_J$-equation may  not be a $C^{\infty}$ smooth map any more but only $C^{m_0}$-smooth if we use $L_k^p$-maps as the  ambient space to work with. In other words, as far as the regularity of each individual map is concerned, it may happen that nothing is special  about  a map being a solution of the perturbed  ${\bar {\partial }}_J$-equation.  Here of course our assumption is that all geometric data  are  of class $C^{\infty}.$
 
 To see the smoothness of the solution of the perturbed  ${\bar {\partial }}_J$-equation,  we may  use local coordinate charts centred at  smooth curves since  such curves are dense in the space of $L_k^p$-curves. In the  case of using "constant" perturbations,   each element of the perturbations is obtained by the parallel transporting a smooth section at the center over a local slice first, then multiplying the cut-off function to make the constant section localized. In other words in the local chart and local trivialization, the the perturbed  ${\bar {\partial }}_J$-equation takes the following form.
  $[{\bar {\partial }}]_J \xi+\gamma (\xi)\nu(\xi)=0$. Here $[{\bar {\partial }}]_J:L_k^p(\Sigma, f^*(TM))\rightarrow L_{k-1}^p(\Sigma, \Lambda ^{0,1}(f^*(TM)))$ is the ${\bar {\partial }}_J$-operator written in the local chart and trivialization centered
   at $f$, $\gamma$ is the cut-off function supporter in a small neighbourhood of $f$ and $\nu(\xi)=\nu(0)$ is a constant section. Therefore, if $\xi$ is a solution of this perturbed  ${\bar {\partial }}_J$-equation, it is the solution of the equation for $\eta$: $
  [{\bar {\partial }}]_J \eta+\gamma (\xi)\nu(0)=0$. This is  just a inhomogeneous  equation  for the usual  quasi-linear elliptic  operator,
    ${\bar {\partial }}_J$, with smooth data. Note that in this equation, $\xi$ is fixed, $\gamma (\xi)$ is a fixed real number. Therefore, the solution is smooth by the elliptic estimate for ${\bar {\partial }}_J$-operator.
  In   our  case  of using localized geometric perturbations, we get  similar equation with the term $\gamma (\xi)\nu(\xi)$ being replaced by 
   linear combinations of the terms $\gamma (\xi) \omega \cdot [X(\xi)]$, where $\omega$ is a smooth  $(0,1)$-form on $\Sigma$, $X$ is a smooth vector field on $M$,   $ X(\xi)$ is   the  pull-back of $X$ by the map $Exp_f\xi$ and $ [X(\xi)]$ is $ X(\xi)$ written in the local trivialization at $f$.   Since    a generic
  $\xi$ is of class $L_k^p, $ $ X(\xi)$, and hence $\omega \cdot [X(\xi)]$ is  of class  $L_{k}^p$.  The same argument  above implies that   the solution is in $L_{k+1}^p(\Sigma, f^*(TM))$. Therefore, we get the same conclusion by bootstrapping.


As mentioned above,   the method  in [LT]  gives a stratified topological manifold structure on the extended  moduli space.  On the other hand   in this paper and its sequels [L2] and [L3], we go further to establish the 
$C^{m_0}$-smoothness of the perturbed moduli spaces. This will be done by   using   the  induced  weakly smooth  structures from the ambient space.



The details  on how this can be done is  given  in Sec. 3. Here we just mention briefly how to get a $C^{m_0}$-smooth moduli space from a compatible collection of perturbed
moduli space defined on local slice by using the weakly smooth structure defined in Sec. 2.

\medskip
\vspace{2mm}
\noindent ${\bullet}$  $C^{m_0}$-Smooth Moduli Space and Smoothness of Evaluation Maps  on the Ambient space:
 
 In Sec. 2 and 3, we have defined a  function or section on the ambient space of unparametrized stable $L_k^p$-maps to be weakly smooth
  if it  is  smooth  with respect to some covering slices. The germs of such functions give rise the weakly smooth structure on the ambient space.  At first sight, this obvious notion is just a convenient way to talk about   smoothness in the present situation  if we do not want to put any further structure like  sc-smoothness   on the ambient space but only regard  it  as a topological Banach manifold.
  
  What is not so  obvious is that   
   in the case of GW and Floer type theories, the weakly  smooth functions defined on the ambient spaces induce  honest 
  smooth structures on any finite dimensional topological manifolds embedded  in the ambient spaces as weakly smooth submanifolds. In particular,  the perturbed moduli space is one of such topological manifolds by the construction in  sec. 3. In fact,  the smooth structure on the moduli space is determined by the obvious weakly smooth maps, the evaluation maps, on the local slice of the ambient space.

More specifically, using the "global" geometric perturbations in Sec. 3, one obtains a collection of compatible local moduli spaces for a generic perturbation. Each of such local moduli spaces is lying inside a local slice and is a $C^{m_0}$-smooth submanifold
inside the local slice. They are compatible in the sense that on the overlaps of the  images of theses local slices regarded as charts in the ambient space of the unparametrized stable  $L_k^p$-maps, these local moduli spaces  are mapping to each other under the ambient transition maps of the local slices. Since these transition functions are continuous, it follows  that the "global" moduli space
is at least a topological manifold. This together with some of related   properties is sufficient for most of all applications. At first sight, it seems
 that this  the best one can get since these transition functions  are only  continuous. However, the global moduli space above is an example of a finite dimensional  weakly smooth submanifold inside the ambient space. We proved in Sec. 3 that in this case, the ambient space is a weakly smooth Banach manifold which is effective with respect to the equivalence class of the "covering data". This means that there are sufficiently many weakly $C^{m_0}$-smooth functions defined on the
 ambient space so that they give any  finite dimensional  weakly smooth submanifold as above an honest $C^{m_0}$-smooth structure. 
 In fact, it is easy to see that the $l$-fold evaluation  map $ev_{x}:{\tilde  {\cal B} }\rightarrow M^l$ from the  space ${\tilde {\cal B}}$ of parametrized $L_k^p$ maps to $M^l$ given by $ev_{x}(f)=(f(x_1), \cdots, f(x_l))$ is $C^{\infty}$ smooth. Here $x=(x_1, \cdots, x_l)$ with $x_i\not =x_j$ for $i\not = j.$  Of course,   the restrictions of $ev_{x}$ to    any  local slice is still $C^{\infty}$ smooth.   It is proved in Sec.3 that  it is  $C^{m_0}$-smooth viewed in any other slices.  Therefore we get  a collection of weakly smooth functions of class $C^{m_0}$ on open sets of ${ {\cal B}}$. By choosing  $l$ and $x$ properly, the collection of such functions serves as "coordinates" for these local perturbed moduli spaces.
 In other words,
 in the case of the ambient topological Banach manifold appeared in GW and Floer type theories, these evaluation maps on the slices of the ambient  space are already  sufficient to detect the honest smoothness of an
  embedded finite dimensional object. 
 
\
\medskip
\vspace{2mm}
\noindent ${\bullet}$ Flat Charts

A local flat  chart  inside a local slice $W(f)=W_{\epsilon}({f}, H)$  of the space of  stable $L_k^p$-maps is the collect all $L_k^p$-maps in $W_{\epsilon}({f})$ which are $J$-holomorphic on a collection of prescribed  disc neighbourhoods of the double points of the domain $\Sigma$ of $f$ and the corresponding  prescribed annulus   of  the domain $\Sigma_{\alpha}$ of the pre-gluing $f_{\alpha}.$

To give a more precise definition, 
 let $W^{D_1}_{\epsilon}({f})$ be one of the  strata of $W(f)$. As above, by considering $W (f)$ as a local fiberation over the Deligne-Mumford moduli space, we   decompose $W (f)$   as the union of slices (fibers) with fixed  gluing parameters, $W(f)=\cup _{\alpha  \in \Lambda(D_1)}W^{\alpha}({f_{\alpha}})$, where ${f_{\alpha}}$ is a pre-gluing of $f$  and $\Lambda (D_1) $ is the set of the "gluing" parameters within the stratum $D_1$.  Note that here gluing parameter $\alpha=(\alpha_t, \alpha_n) $ such that $f_{\alpha_n}$ is the  pre-gluing along "normal" direction while  the parameter $\alpha_t$ describes the local  deformations of the domain $\Sigma_f$ by the motions of its double points. Following [LT],   collection the all pre-gluing maps $f_\alpha, \alpha \in \Lambda_{\delta}$ obtained from $f$  will be  called "base" maps.
 
  To define the flat charts and the related bundle, we start with one of the  lowest strata, denoted by $D.$ 
 
 Then $W^D_{\epsilon}({f})$ corresponds to 
the collection of above   $W^{\alpha}({f_{\alpha}})$ with ${\alpha_n =0.}$
  
  On this  lowest stratum ,   the local "bundle"  ${ {\cal L}}(f)$ restricted to $W^D_{\epsilon}({f})$, denoted by ${{\cal L}}_D(f)$ has a trivialization so that each fiber can be identified with $({\widetilde {\cal L}}(f))_f=L_{k-1}^p(\Sigma_0, \Lambda^{0, 1}( f^*(TM))).$
  More specifically, one each  $W^{\alpha}({f_{\alpha}})$ with a fixed $\alpha \in \Lambda(D_1)$, the domain $\Sigma_{\alpha}$  is fixed, 
  the standard  trivialization induced by parallel transport in this paper  with some modification is still applicable (see [LT] for 
  instance).   On the other hand, the "standard" differomorphisms in [LT] to identify the domain $\Sigma_f=\Sigma_0$ with $\Sigma_{\alpha}$ 
  induce the identifications of the corresponding fibers. These differomorphisms are used to define the local  $C^{\infty}$ product structure of the 
  universal curve over Deligne-Mumford moduli space mentioned  before. They can be realized as time $1$-maps of the corresponding flows 
  generated by  vector fields supported on  small annuli around double points. Note that the vector fields are not holomorphic on these  
  annuli. Consequently,  this last identification is only over ${\bf R}.$ Though this does not affect our definition below, it  seems better 
  to use a fiberwise version define the objects for each fixed ${\alpha}$.

 In any case,  combing these together we get the trivialization of the local bundle over the lowest stratum $D$.

 Consider  the following sub-bundles:
 ${{ {\cal L}}}_{D, Loc, \delta}(f)$  of
${{{\cal L}}}_{D}(f)$ defined by
$$({{ {\cal L}}}_{D, Loc, \delta}(f))_g=\{\xi \,|\, \xi\in
({{ {\cal L}}}_{D}(f))_g, \, \xi=0\, on\, \, each\,\, \delta-disc
\,\, around\,\,a\,\, double\, \,\,point\}.$$
Let
${{W}}^D_{\epsilon, Loc,\delta}(f)= {\bar{\partial}}_{J}^{-1}
({{ {\cal L}}}_{D, Loc, \delta}(f)).$ This is one of the local flat chart in  $W(f )$  restricted to the its lowest stratum.
Hence any element in ${{W}}^D_{\epsilon, Loc,\delta}(f)$ has the property that it is $J$-holomorphic on the $ \delta-discs$ above.

To define the corresponding objects on a higher stratum, consider the restriction of the bundle ${{ {\cal L}}}_{D, Loc, \delta}(f)$ to the base maps $f_{\alpha_t}$ with ${\alpha_t}$ in $D$.
According to what we have done  for a fixed stratum, we only need to know how to extend the fiber $({{ {\cal L}}}_{D, Loc, \delta}(f))_{f_{\alpha_t}}$ 
of the  bundle above into the fiber ( of the bundle to be defined) over  $f_{\alpha}.$  
Here $\alpha=(\alpha_t, \alpha_n)$ and $f_{\alpha}$ is obtained from $f_{\alpha_t}$ by the normal pre-gluing with the gluing parameter 
$\alpha_n.$ The key point is that when $\alpha_t$ and $\delta$  are fixed, for $\|\alpha_n\|<< \delta,$ the pre-gluing from the domain 
$\Sigma_{\alpha_t}$ to $\Sigma_{\alpha} =\Sigma_{(\alpha_t, \alpha_n)}$ identifies $\Sigma_{\alpha_t}\setminus\{\delta$ -discs at double 
points $\}$ with $\Sigma_{\alpha}\setminus\{"\delta$ -annuli" around double points $\}.$  Here for simplicity,  we have assumed that  non of 
the component of $\alpha_n$ is equal to zero so that $\alpha$ is in the top stratum.

This completes the construction of the bundle over the base maps $f_{\alpha}, \alpha\in \Lambda_{\delta}. $ As above by parallel transporting  these  "central fiber" over $f_{\alpha}$, we get the sub-bundle  over $W^{\alpha}_{\epsilon}(f)$, denoted by  ${{ {\cal L}}}_{\alpha, Loc, \delta}(f_{\alpha})$

 We   define ${{W}}^{\alpha}_{\epsilon, Loc,\delta}(f_{\alpha})= {\bar{\partial}}_{J}^{-1}
({ {\cal L}}_{{\alpha}, Loc, \delta}(f_{\alpha})).$ Then in the general case any element in ${{W}}^{\alpha}_{\epsilon, Loc,\delta}(f_{\alpha})$  is  $J$-holomorphic on the  $ \delta-$discs or $ "\delta-$annuli"  of $\Sigma_{\alpha}$.

The union of these ${{W}}^{\alpha}_{\epsilon, Loc,\delta}(f_{\alpha})$, denoted by ${{W}}_{\epsilon, Loc,\delta}(f)$ is one of the flat charts  inside the local slice $W(f).$

Of course, when $\delta$ is getting smaller,  so is the size of $\Lambda$.  Since a  transformation of a flat chart in one slice into another slice 
is induced by transformations on the domains, the size and "shapes" of these charts given by  parameters $"\delta"$ may change under the transformations.  Since for our purpose here, there are only finitely many fixed local uniformizers are used to covering ${\cal M}(A, J)$, by allowing different shapes of ${\delta}$-discs in the obvious sense in this situation,  the collection of all such flat charts are closed under
coordinate transformations between uniformizers. They define the flat structure that we are looking for.

 Let  ${{W}}^0_{\epsilon}(f)=\cup_{\delta}{{W}}_{\epsilon, Loc,\delta}(f)=\cup_{\alpha,\delta}{{W}}^{\alpha}_{\epsilon, Loc,\delta}(f)$ be the union  of these local flat charts and ${ {\cal L}}^0=\cup_{\delta}{ {\cal L}}_{ Loc, \delta}(f) =\cup_{\alpha,\delta}{ {\cal L}}_{{\alpha}, Loc, \delta}(f_{\alpha}).$

Here ${ {\cal L}}_{ Loc, \delta}(f) =\cup_{\alpha}{ {\cal L}}_{{\alpha}, Loc, \delta}(f_{\alpha}).$

Note:  Part of the above construction, the definition of ${{W}}^D_{\epsilon, Loc,\delta}(f)$  was  already introduced in [LT], page 25. However,  the key  part of the construction on flat charts above,  the definition of  ${{W}}^{\alpha}_{\epsilon, Loc,\delta}(f_{\alpha})$   dealing with the norm  gluing parameters  was  missing in [LT].

  \medskip
\vspace{2mm}
\noindent 
$\bullet$ Properties of the flat charts:

\medskip
\vspace{2mm}

(I)  The indices of  the ${\bar {\partial}}_J$-section on any strata of  a local slice are the same as the indices of the  section restricted to  the corresponding strata of a flat chart considered as a section of the corresponding sub-bundle.  

(II) The part of moduli space of stable $J$-holomorphic maps inside a given slice near an "end" is contained in any of  these flat charts. 

(III)  For a fixed ${\alpha}$ with $\|\alpha\|<<\delta$, each  flat chart ${{W}}^{\alpha}_{\epsilon, Loc,\delta}(f)$ of "size" $\delta$ is 
a Banach manifold.

(IV)  For a fixed ${\delta}$ small enough, all ${{W}}^{\alpha}_{\epsilon, Loc,\delta}(f)$ fit together to form a Banach manifold.

The property (III) follows from the following facts:  (i) the  cokernel $K_{ f_{\alpha}}$ of $D{{\bar \partial }}_{J, f_{\alpha}}$ is of finite dimensional; (ii) 
elements in $K_{ f_{\alpha}}$  satisfy the unique continuation principle. In deed, (i) and (ii) imply that 
${\bar {\partial }}_{J}:{W}^{\alpha}_{\epsilon}(f)\rightarrow {\cal L}(f)$ is transversal to the sub-bundle ${ {\cal L}}_{{\alpha}, Loc, \delta}(f_{\alpha})$ 
 so that its inverse image a Banach manifold. Note that here we have used the fact that under the local trivializations of the two bundles, 
 the central fiber of the subbundle is a closed subspace of the corresponding central fiber. This is true  under  our assumption that
$m_0>1.$  The property (I) can be  derived , for instance,   by using (III).

\medskip
\vspace{2mm}
\noindent 
$\bullet$  $\bullet$  Gluing and Shape of  flat charts:
It follows from (III) that in the standard local coordinate chart 
$Exp_{f_{\alpha}}:V^{\alpha}\subset\L_k^P(\Sigma_{\alpha}, f_{\alpha}^*(TM), h)\rightarrow {{W}}^{\alpha}_{\epsilon}(f)$ of the local slice, 
the flat chart $Exp_{f_{\alpha}}^{-1}({{W}}^{\alpha}_{\epsilon, Loc,\delta}(f)), $ denoted by ${{V}}^{\alpha}_{\epsilon, Loc,\delta}(f)$ in  
$V^{\alpha}$ is  realized as a graph 
of a map from the open  ball $B(r_\alpha)$ in the tangent space of ${{V}}^{\alpha}_{\epsilon, Loc,\delta}(f)$ at $f_{\alpha}$ to its orthogonal complement (assuming that $p=2$ so that the $L_k^p$-space 
is a Hilbert space). Here $r_\alpha$ is the radius of the ball depending on $\alpha.$  Therefore the  key point to prove (IV) is  to prove that there is a positive lower bound for $r_{\alpha}$ independent of $\alpha$.

Note that it follows from a parametrized version of implicit function theorem that when ${\alpha}$ is moving within  same stratum slightly there is such a bound. Therefore, the real question is about the motion form $\alpha_0=(\alpha_t, 0)$ to $\alpha=(\alpha_t, \alpha_n)$  with  normal gluing parameter $\alpha_n$. 

In this case, we need to show that the radius $r_{\alpha}$ does not goes to zero as $\|\alpha_n\|$ goes to zero so that  the size and shape of each ${{V}}^{\alpha}_{\epsilon, Loc,\delta}(f)$ is comparable  with the fixed ${{V}}^{\alpha_0}_{\epsilon, Loc,\delta}(f).$

This is exactly the  situation  that the so called  gluing  technique in [L0] and [LT] can be used to  deal with. The main estimate in [LT], page 32-41, with some modifications
implies that (IV) above holds.  The detail of  the proof of  (IV)  is given  in [L4]. We will give a outline of its proof below.  As a by-product of the proof, we have also obtained two  different  new proofs for the gluing of $J$-holomorphic maps.

 The following theorem  summarizes  these properties of the flat charts and  their consequences.

 \begin{theorem}\label{2.7}

For any fixed ${\delta}, $ the local flat chart ${{W}}_{\epsilon, Loc,\delta}(f)=\cup_{\|\alpha\|<<\delta}{{W}}^{\alpha}_{\epsilon, Loc,\delta}(f)$  is a Banach manifold ("ball") and  ${ {\cal L}}_{ Loc, \delta}(f) =\cup_{\alpha}{ {\cal L}}_{{\alpha}, Loc, \delta}(f_{\alpha})\rightarrow {{W}}_{\epsilon, Loc,\delta}(f) $ is    a Banach bundle. The restriction of the ${\bar{\partial}}_{J}$-section $s$ to the flat chart,   still denoted by $s$ is a  smooth Fredholm  section,   whose linearizations have the same indices as the ones for the  "old" ${\bar{\partial}}_{J}$-section  along
 all strata.

The zero sets
${\bar\partial}^{-1}_{J}(0)$ in ${{W}}_{\epsilon, Loc,\delta}(f)$,
when projected to
${\cal B}$  is  just ${ {{\cal M}}}(J,A)\cap
{ [ W]}_{\epsilon}(f;\bf H)$.

Furthermore, the $L_k^p$-topology and "weak"-topology on  any flat chart above are equivalent.

\end{theorem}

\medskip
\vspace{2mm}
\noindent \\
$\bullet$  $\bullet$  Exponentially weighted $L_k^p$-norm and further fibration:

Note that  in the last statement of the above theorem, the $L_k^p$-topology here is measured in term of cylindrical coordinates rather 
than the "spherical" ones. The same is applied to   the discussions on property (IV).  The necessity of using cylindrical coordinates on the $\delta$-
discs $D^{\pm}_{\delta}(\pm d)$ centred at each double point $d={\pm} d$ of the domain $\Sigma=\Sigma_{\alpha_0}=\Sigma_f$ is clear since as 
$f_{\alpha}:\Sigma_{\alpha}\rightarrow M$ degenerates into $f:\Sigma_{\alpha_0}\rightarrow M$, in the "spherical" coordinate the injective 
radius of $\Sigma_{\alpha}$ goes to zero, which is not adequate for using
general $L_k^p$-norms for analytic set-up here unless $k=1$ that implies that $m_0=0$. 

In any lowest stratum ${  W}^D_{\epsilon}(f;\bf H)$, consider one of its fibers,  for instance, for simplicity,  the   central fiber ${  W}^{\alpha_0}_{\epsilon}(f;\bf H).$ Then all the elements have the same  domain $\Sigma=\Sigma_{\alpha_0}=\Sigma_f$.

Now  we identify each pair of   $\delta$-discs at a double point $d=\pm d$ of $\Sigma$ with  a pair of half cylinders $C_{\pm}\simeq {\bf R}^{+}\times S^1$ of infinite length  with $\{\pm\infty\}\times S^1$ corresponding to the double point. Let $w=\pm w$ and $(s, t)=\pm (s, t)$ be the corresponding spherical and
cylindrical coordinates. Then $w=exp (-(s+it)).$ Clearly,  any $h=Exp_f\xi$ in ${  W}^{\alpha_0}_{\epsilon}(f;\bf H)$ satisfying the condition that $h(d)=f(d)$ satisfies the exponential decay estimate  that $\|D^i\xi\|\sim O(\exp (-s))$ for  $ i\geq 0.$  Moreover, if   $f$ and 
$h$ are smooth,  then for
 any $0<\kappa<1$, $\|\xi\|_{k, p;\kappa}=\|e^{\kappa  s}\xi\|_{k, p}<\infty.$   This means that we should introduce a new fibration ,  $ev^l:{  W}^{\alpha_0}_{\epsilon}(f;{\bf H})\rightarrow  M^{l}$ given by the evaluation map at all double points of $\Sigma$. Here $l$ is the number of these fixed double point. One can  show that both above $ev^l$ and its restriction to any of flat charts inside the stratum,  are  in deed  fibrations over a neighbourhood of $f(d)$ in $M^l$. Therefore,  for any  $m\in M^l $ near $f(d)$ we get the corresponding fibers, ${  W}^{(\alpha_0, m)}_{\epsilon}(f;\bf H)$ and 
${{W}}^{(\alpha_0, m)}_{\epsilon, Loc,\delta}(f)$. As  mentioned above, all elements in ${  W}^{(\alpha_0, m)}_{\epsilon}(f;\bf H)$ are exponential decay and have  a finite $L_k^p$-norm with exponential weight ${\kappa}$ for any $0<\kappa<1$. One can show that at least for 
elements $h=exp_f\xi$ in ${{W}}^{(\alpha_0, m)}_{\epsilon, Loc,\delta}(f), $ we have that the usual $L_k^p$-norm  dominates the exponential  
$L_{k, \kappa}^p$-norm. Moreover, the indices  of the linearlization of ${\bar{\partial }}_J$-operator, $Ds_h$ are the same for any $h$ in 
${{W}}^{\alpha_0, m}_{\epsilon}(f) $ with the property that if it is surjective with respect to the usual $L_k^p$-norm, it is still  
surjective with respect to the exponential norm.   Therefore,  switching  to the exponential norm  all the desired properties are preserved. 
In the rest of the discussion  of this section on flat chart, we will assume that the $L_k^p$-maps are measured  by exponentially  weighted $L_k^p$-norms. The corresponding (enlarged) spaces will be denoted by $^{\kappa}{{W}}^{(\alpha_0, m)}_{\epsilon}(f) ,$ etc.

  We now  extend above discussion to higher stratum.  Note that  by using normal coordinate at $f(d)$ to identify  its  neighbourhood in $M$ 
  with an open ball in $T_{f(d)}M$ and replacing the double point $d$ by the corresponding $S^1$ at infinity in the cylindrical 
  coordinates, the evaluation map $ev^l(h)$ can be calculated by the integral $\int_{S^1}\xi.$  The point  is that this later formula for 
  the fibration given by $ex^l$ is sill applicable  to  the higher strata ${  W}^{\alpha}_{\epsilon}(f;\bf H)$. For simplicity, we assume 
  that $\alpha=(0, \alpha_n)$ with ${\alpha_0}=(0, 0).$ Then we get a corresponding fibration, denoted by  $ev^l_{\alpha}:{  
  W}^{\alpha}_{\epsilon}(f;{\bf H})\rightarrow  M^{l}$ given by the "evaluation map" at all  middle circle $S^1_{\alpha}$ of $\Sigma_{\alpha}$.
  Here $S^1_{\alpha}$ is the middle circle of the $\delta$-neck in $\Sigma_{\alpha}$ obtained from the pre-gluing, and $ev^l_{\alpha}$ is defined by 
  the similar  integral formula above by  replacing $S^1$ by $S^1_{\alpha}$. Again, one can show that $ev^l_{\alpha}$  are   fibrations for the
   both cases.
 
 Now using the corresponding cylindrical coordinates of finite length $-log \|\alpha\|$ defined on $\Sigma_{\alpha}$ to define the $L^p_{k, ;\kappa}$-norm , we obtain the corresponding fibers  $^{\kappa}{{W}}^{(\alpha, m)}_{\epsilon}(f) $ and $^{\kappa}{{W}}^{(\alpha_0, m)}_{\epsilon, Loc,\delta}(f)$.
 
Let  $Exp_{f_{\alpha}}:$ $  ^{\kappa}V^{(\alpha, m)}\subset L_{(k;\kappa)}^p(\Sigma_{\alpha}, f_{\alpha}^*(TM); h, S^1_{\alpha})\rightarrow$ $ 
^{\kappa}{{W}}^{(\alpha, m)}_{\epsilon}(f)$ be the local coordinate of the local slice with fixed parameter $(\alpha, m).$ Here the new 
constrain  for an element $\xi$ in $L_{(k;\kappa)}^p(\Sigma_{\alpha}, f_{\alpha}^*(TM); h, S^1_{\alpha})$ is that the integral  
$\int_{S^1_\alpha}\xi=0.$
 
 Let $[s]^{(\alpha, m)}:$ $ ^{\kappa}V^{(\alpha, m)}\subset L_{(k;\kappa)}^p(\Sigma_{\alpha}, f_{\alpha}^*(TM); h, S^1_{\alpha})\rightarrow $
 $
L_{(k;\kappa)}^p(\Sigma_{\alpha}, \Lambda^{0, 1}(f_{\alpha}^*(TM))$ be the ${\bar {\partial}}_J$-operator written  in the standard local coordinate chart and trivialization. Note that here we have assumed that $\alpha=(0, \alpha_n)$ is the normal gluing parameter.
 
 We now make two assumptions on $$D[s]^{(\alpha_0, m)}_f:L_{(k;\kappa)}^p(\Sigma_{\alpha_0}, f_{\alpha_0}^*(TM); h, S^1_{\alpha_0})\rightarrow 
L_{(k;\kappa)}^p(\Sigma_{\alpha_0}, \Lambda^{0, 1}(f_{\alpha_0}^*(TM))$$.
\medskip
\vspace{2mm}
\noindent \\
$\bullet$  $\bullet$ (A1): It is surjective.
 
 \medskip
\vspace{2mm}
\noindent \\
$\bullet$  $\bullet$ (A2): It is injective.
 
 The second assumption can be achieved locally by introducing more marked points $y$ located on the 'fixed part' of $\Sigma_{\alpha}$  and requiring that $h(y)$ lies on  some prescribed local hypersurfaces of $M$ for $h$ in  $^{\kappa}{{W}}^{(\alpha_0, m)}_{\epsilon}(f)$. Here as before, the 'fixed part'  $\Sigma^0_{\alpha}$=$\Sigma_{\alpha}\setminus \{$ all $\delta$-annuli around double points$\}$, 
 which can be identified with the 'fixed part'  $\Sigma^0=\Sigma\setminus \{$ all $\delta$-discs around double points$\}$  by the pregluing.   
 
Then the key to prove (IV) is the following proposition

  \begin{pro}\label{2.7}

There exists a  fixed ${r_0}>0 $ independent $\alpha$ such that  for any $0<r\leq r_0,$  when $\|\alpha\|$  is small enough,
 $[s]^{(\alpha, m)}:B_r^{(\alpha, m)}\subset L_{(k;\kappa)}^p(\Sigma_{\alpha}, f_{\alpha}^*(TM); h, S^1_{\alpha})\rightarrow 
L_{(k-1;\kappa)}^p(\Sigma_{\alpha}, \Lambda^{0, 1}(f_{\alpha}^*(TM))$ is a diffeomorphism form a ball of radius $r$ to its image.

\end{pro}

Note that for $\alpha=\alpha_n$ with all non-zero entry, all such  balls $B_r^{(\alpha, m)}$ can be identified each other. However, their dimension  are "infinitely smaller" than the one for $\alpha=(0, 0).$ 

By applying Picard method, the proof of this proposition follows from the following "main estimate" in the gluing technique mentioned above. 
 \begin{pro}\label{2.7}

Consider  a family of linear operators $$D[s]^{(\alpha, m)}_{f_{\alpha}}: L_{(k;\kappa)}^p(\Sigma_{\alpha}, f_{\alpha}^*(TM); h, S^1_{\alpha})\rightarrow $$ $$
L_{(k-1;\kappa)}^p(\Sigma_{\alpha}, \Lambda^{0, 1}(f_{\alpha}^*(TM)).$$ 
There exists a  fixed $c_0>0 $ independent $\alpha$ such that for $\|\alpha\|$ small enough,$\|D[s]^{(\alpha, m)}_{f_{\alpha}}\|>c_0.$

\end{pro}
\proof 

The easy half  proof of the main estimate in [LT], page 38-41 or the half of the proof in [L0], page 276-277 implies the Proposition.

\QED

This is  the first "new" proof of the  gluing we mentioned above.  This proof  still goes along the same line as the old proof. But because the  fibration above, we have eliminated the extra terms  in the exponential weighted $L_k^p$-norms used in [LT] and [L0] coming from the integral
 along $S^1_{\alpha}$ above.  This simplifies the argument there and make the proof here almost as simple as Floer's gluing for non-degenerate broken connection orbits  in [F], page 599.

 \medskip
\vspace{2mm}
\noindent 
$\bullet$  $\bullet$ Heuristic reasoning for the second proof of the gluing:

For our purpose, it is sufficient  to establish  a weaker form of the first proposition above for the flat charts.  

 Let $Exp_{f_{\alpha}}^{-1}(^{\kappa}{{W}}^{\alpha, m}_{\epsilon, Loc,\delta}(f)), $ denoted by $^{\kappa}{{V}}^{(\alpha, m)}_{ Loc,\delta} $  be the inverse image of flat chart in $^{\kappa}V^{(\alpha, m)}.$
 
  Denote the restriction of $[s]^{(\alpha, m)}$  to $^{\kappa}{{V}}^{(\alpha, m)}_{ Loc,\delta} $  by the same notation,
 $[s]^{(\alpha, m)}:$ $ ^{\kappa}V^{(\alpha, m)}_{ Loc,\delta}\rightarrow         
L_{(k;\kappa; Loc,\delta)}^p(\Sigma_{\alpha}, \Lambda^{0, 1}(f_{\alpha}^*(TM)). $  Here the right hand side is the fiber at $
f_{\alpha}$
of  the local bundle used to define the flat chart.

Then   the assumptions (A1) and (A2) imply the two corresponding statements for this new $D[s]^{(\alpha_0, m)}_f .$

The weaker form of the first proposition above is the following proposition.

\begin{pro}\label{2.7}

There exists a  fixed ${r_0}>0 $ independent $\alpha$ such that  for any $0<r\leq r_0,$  when $\|\alpha\|$  is small enough,
 $[s]^{(\alpha, m)}:$ $^{\kappa}V^{(\alpha, m)}_{ Loc,\delta}\cap B_r^{(\alpha, m)}\rightarrow 
L_{(k-1;\kappa; Loc,\delta)}^p(\Sigma_{\alpha}, \Lambda^{0, 1}(f_{\alpha}^*(TM)) $ is a diffeomorphism form a ball of radius $r$ to its image.
Here $ B_r^{(\alpha, m)}$ is the ball of radius $r$ in $ L_{(k;\kappa)}^p(\Sigma_{\alpha}, f_{\alpha}^*(TM); h, S^1_{\alpha})$

  Moreover, $^{\kappa}V^{(\alpha, m)}_{ Loc,\delta}\cap B_r^{(\alpha, m)}$ can be realized as   the graph of a smooth function from the ball 
  of radius $r'$, denoted by ${\hat B}_{r'}^{\alpha}$,   in the tangent space of ${{V}}^{(\alpha, m)}_{ Loc,\delta}$ at $f_{\alpha}$, 
  denoted by $T_{\alpha}$,  to its orthogonal complement (assuming that $p=2$) in 
$L_{(k;\kappa)}^p(\Sigma_{\alpha}, f_{\alpha}^*(TM); h, S^1_{\alpha}) , $ denoted by $O_{\alpha}.$ Consequently, its image in
 $$L_{(k-1;\kappa; Loc,\delta)}^p(\Sigma_{\alpha}, \Lambda^{0, 1}(f_{\alpha}^*(TM)) $$ contains a  ball ${\bar  { B}}_{r''}^{\alpha}$
 of fixed radius $r''$ independent of $\alpha$. Note that by our construction, these images can be thought inside  the same space 
 $$L_{(k-1;\kappa; Loc,\delta)}^p(\Sigma_{\alpha_0}, \Lambda^{0, 1}(f_{\alpha_0}^*(TM)),$$  the central fiber denoted by ${\cal L}_{\delta}$ for short.

\end{pro}

\proof

(Sketched):

Let $G_{\alpha}: 
 {\hat B}_{r'_{\alpha}}^{\alpha}\subset T_{\alpha}\rightarrow  O_{\alpha} $ be the map representing 
 $^{\kappa}V^{(\alpha, m)}_{ Loc,\delta}\cap B_r^{(\alpha, m)}$. Denote its graph $$id\oplus G_{\alpha}: 
 {\hat B}_{r'_{\alpha}}^{\alpha}\rightarrow  {\hat B}_{r'_{\alpha}}^{\alpha}\oplus  O_{\alpha} $$ $$\subset L_{(k;\kappa)}^p(\Sigma_{\alpha}, f_{\alpha}^*(TM); h, S^1_{\alpha}) $$  by    $F_{\alpha}.$ 
  Let  
  $S^{\alpha}=[s]^{\alpha, m}\circ F_{\alpha}:{\hat B}_{r'_{\alpha}}^{\alpha} \rightarrow {\cal L}_{\delta}.$ We need to show that there  is a positive $r'\leq r'_{\alpha}$ independent of $\alpha$ such that  $S^{\alpha}$ restricted to  
 ${\hat B}_{r'}^{\alpha} $  is a diffeomorphism to its image for all $\alpha$. Clearly we may assume that $S^{\alpha_0}$ does so.

It is sufficient to  prove the following.

(I)$DS^{\alpha}_{f_{\alpha}}:T_{\alpha}\rightarrow {\cal L}_{\delta}$    has an uniform inverse.

This  is the  most difficult step for the usual gluing estimate in [l0] and [LT]. But it is almost trivial here simply because the sections 
$\xi$ in $T_{\alpha}$ restricted to the $\delta$-annulus $A_{\delta,\alpha} \subset \Sigma_{\alpha}$  satisfy the linearised 
${\bar {\partial}}_J$-equation at $f_{\alpha}$. In cylindrical coordinates, they  are exponential decay in the sense that their norms satisfy
$$\|\xi|_{A_{\delta,\alpha}}\|_{k, p;\kappa}\leq C exp(-(-ln\delta))\|\xi|_{C_{\delta,\alpha}}\|_{k, p;\kappa}= C\delta\|\xi|_{C_{\delta,\alpha}}\|_{k, p;\kappa}.$$ 
Here $C$ is a constant independent of  $\alpha$, $C_{\delta,\alpha}$ is the two cylinders of length one ( in cylindrical 
coordinates) at the ends of $A_{\delta,\alpha}$. Note that all $C_{\delta,\alpha}$ can be identified with $C_{\delta,\alpha_0}$  in $
\Sigma_{\alpha_0}.$ Therefore, $\|\xi|_{A_{\delta,\alpha}}\|_{k, p;\kappa}$  are ignorable for $\delta$ small enough. For such small $\delta$ we may assume that for all $\alpha,$ 
$f_{\alpha}|_{\Sigma^0_{\delta,\alpha}}$  are  equal to the restriction of $f_{\alpha_0} $   to $\Sigma^0_{\delta,\alpha_0}$ by the construction the pre-gluing. 
Here  
$\Sigma^0_{\delta,\alpha}=\Sigma_{\alpha}\setminus A_{\delta,\alpha}$ which can be identified with $\Sigma^0_{\delta,\alpha_0}$.
Therefore,   the  operator norm $$\|(DS^{\alpha}_{f_{\alpha}})^{-1}\|\sim \|(DS^{\alpha}_{f_{\alpha}}|_{\Sigma^0_{\delta,\alpha}})^{-1}\|$$ 
$$\sim \|(DS^{\alpha_0}_{f_{\alpha_0}}|_{\Sigma^0_{\delta,\alpha_0}})^{-1}\|\sim \|(DS^{\alpha_0}_{f_{\alpha_0}})^{-1}\|.$$ 

(II)  Assume that the $r'>0$  exists first. It  will be proved later.

(III) Second order estimate:   Since the desired second order estimate for $[s]^{\alpha, m}$ for Picard method is well-known by Floer's work and is independent of 
${\alpha}$, it easy the see that we only nee to show that there is a constant $C$ independent of ${\alpha}$ such that $\|
F_{\alpha}(\xi)\|\leq C\|\xi\|$ for any $\xi\in 
{\hat B}_{r'}^{\alpha} \subset {\hat B}_{r'_{\alpha}}^{\alpha}.$  For any $\xi_{\alpha}\in {\hat B}_{r'}^{\alpha} $, by using cut-off 
functions one can construct a corresponding  $\xi_{\alpha_0}\in {\hat B}_{r'}^{\alpha_0} $ such that upto an "uniform ignorable" 
exponential   decay factor $\|F_{\alpha}(\xi_{\alpha})\|\sim \|F_{\alpha_0}(\xi_{\alpha_0})\|$. The proof of this last statement is straight 
forward but takes quite a few steps. It says that the shape of the graph of $F_{\alpha}$ is similar to the fixed  graph of  $F_{\alpha_0}.$  we refer reader to [L4] for the details of this step, but just mention that  here one needs to use the surjectivity of $D[s]_f$ stated in (A1) not just the corresponding one in (I) above. In other words, we can not give a complete "intrinsic" proof for this proposition.

(II) follows form (III):  Note that in the proof of (III) above, we may replace $r'$ by $r'_{\alpha}$. We only need to show that 
$r'_{\alpha}$ is bounded below by a positive number. In other words, 
 the graph of $F_{\alpha}$ does not get off $^{\kappa}V^{(\alpha, m)}$ too fast. This is indeed the case since  by (III) it has the "same" increasing 
 rate as the fixed $F_{\alpha_0}$.

\QED

This second proof is not necessarily shorter in details. But it is conceptually more elementary and almost trivialize the gluing analysis. It  says 
 that in the set-up above one essentially can freely manipulate the contributions from the neck areas by multiplying them with cut-off function without  affecting the any related estimates, a conclusion that one  tried  to get for the usual gluing estimate  with quite effort.

 \medskip
\vspace{2mm}
\noindent 
$\bullet$  $\bullet$ Fredholm  theory and tautological coordinate on flat charts:

In this setting to establish the Fredholm for the local flat charts and their related  bundles essentially amounts to  know how to take derivatives along "normal" direction $\alpha=(\alpha_0, \alpha_n)$ at a point $h$ in the $\alpha_0$-stratum 
 $^{\kappa}V^{(\alpha_0, m)}_{ Loc,\delta}$ for 
$[s]^{(\alpha, m)}:$ $^{\kappa}V^{(\alpha, m)}_{ Loc,\delta}\cap B_r^{(\alpha, m)}\rightarrow 
L_{(k-1;\kappa; Loc,\delta)}^p(\Sigma_{\alpha}, \Lambda^{0, 1}(f_{\alpha}^*(TM)) $ and related sections used to obtain local perturbed moduli space.
Denote $^{\kappa}V^{(\alpha, m)}_{ Loc,\delta}\cap B_r^{(\alpha, m)}$ by $^{\kappa}V^{(\alpha, m)}_{ Loc,\delta}$ for short.

Still work under the assumptions (A1) and (A2). The key point is to introduce a product structure near the "end" 
$^{\kappa}V^{(\alpha_0, m)}_{ Loc,\delta}$. That is that we need to identify $^{\kappa}V^{(\alpha, m)}_{ Loc,\delta}$ with all entries  $\alpha_n\not= 0$ with $^{\kappa}V^{(\alpha_0, m)}_{ Loc,\delta}.$  By the above proposition, for these flat charts,  we have much better chance to have such  identifications. However, the identification is still not immediately since we are not in the finite dimensional situation. On the other hand
by definition when $\alpha_t=\alpha_0$ is fixed,  for  all $\alpha$ we have
$$L_{(k-1;\kappa; Loc,\delta)}^p(\Sigma_{\alpha}, \Lambda^{0, 1}(f_{\alpha}^*(TM))=L_{(k-1;\kappa; Loc,\delta)}^p(\Sigma_{\alpha_0}, \Lambda^{0, 1}(f_{\alpha_0}^*(TM)).$$  Therefore, by the above proposition, we may assume all  $^{\kappa}V^{(\alpha, m)}_{ Loc,\delta}$ can be identified  with $^{\kappa}V^{(\alpha_0, m)}_{ Loc,\delta}$ since we have the identifications
$^{\kappa}{\bar V}^{(\alpha, m)}_{ Loc,\delta}$ with $^{\kappa}{\bar V}^{(\alpha_0, m)}_{ Loc,\delta}$. Here
 $$^{\kappa}{\bar V}^{(\alpha, m)}_{ Loc,\delta}=[s]^{(\alpha, m)}(^{\kappa}{ V}^{(\alpha, m)}_{ Loc,\delta})$$ and 
 $$^{\kappa}{\bar V}^{(\alpha_0, m)}_{ Loc,\delta}=[s]^{(\alpha, m)}(^{\kappa}{ V}^{(\alpha_0, m)}_{ Loc,\delta})$$ 
 are the "balls"  of almost the same sizes in the fixed target space. Under these identification, we get the 
  product structure for
   $$\cup_{\alpha}^{\kappa}V^{(\alpha, m)}_{ Loc,\delta}\sim $$ $$^{\kappa}{\bar V}^{(\alpha_0, m)}_{ Loc,\delta}\times \Lambda (\alpha_0)$$ where  $\Lambda (\alpha_0)= \{\alpha\}$ is the collection of the local normal gluing parameters.

Any element in $$^{\kappa}{\bar V}^{(\alpha_0, m)}_{ Loc,\delta}\times \Lambda (\alpha_0)$$   is  called a
"tautological coordinate". Therefore we get a new coordinate system for a slice  of the flat chart. Clearly, it is tautologically true that
with respect to this coordinate chart, $[s]$ is automatically smooth along $\alpha$-direction. In fact, it is constant along $\alpha$-direction. 

This resolves the  difficulty  (C) for $[s]$ in a "strange" way and under the assumptions (A1) and (A2). We  already know  how to remove  (A2). We  will not give a complete argument  here on how to deal with (A1) (in [L4]). For our purpose to establish the smoothness of extended moduli spaces by Fredholm  theory here, it is sufficient to know that in the general case, the tautological coordinate chart together with the evaluation map at $y$  introduced above to deal with (A2) gives us the desired coordinate  on the product of the flat chart  with  the  cokernel of $D[s]$ at $f_{\alpha_0}.$ It follows from this that at least in this
fixed flat chart,  the perturbations obtained from the cokernel is also smooth
 at point $f$ along normal direction. Here we may assume that any element 
 in the "cokernel" has the property that  it is equal to zero on  the  $\rho$-discs on ${\Sigma}_{\alpha_0}$ around
 double points  by analytic continuation principle. 
 In order to show that the same is
 true for these perturbations  viewed in any other flat charts, we assume that $\rho$ and $\delta$ are in the same order with  $\delta<< \rho$. Of course, to  obtain  a global section from the local  perturbations here, we need to have  the cut-off function which is obtained from the corresponding $L_{k,\kappa}^p$-norm.  In this setting above, the $L_{k,\kappa}^p$-norm is the measured in the target space of
 $[s]$.

\medskip
\vspace{2mm}
\noindent
$\bullet$  $\bullet$ Local  Fredholm  theory  in the  natural  coordinate chart  of a flat end:

\medskip
\vspace{2mm}
It still remains to directly  establish the local  Fredholm theory  in the natural coordinate chart and trivialization for the bundle on the flat end.
As mentioned above, despite of the fact that each element in the flat chart is rigid along the neck, it is still not immediate to give the required product structure for the local Fredholm theory.  In stead of trying to define a product structure  directly, we construct a new
chart  for  the end which  automatically has a product structure and is "close" to the given flat chart. 
Then we prove that with respect to this new chart the section $s^0$, the restriction of ${\bar \partial}_J$-section is smooth of class $C^{1}$ with respect to the product structure.

To  this end,  recall that the flat chart near the end at $f$ is denoted  by ${{W}}^0_{\epsilon}(f)=\cup_{\delta}{{W}}_{\epsilon, Loc,\delta}(f)=\cup_{\alpha,\delta}{{W}}^{\alpha}_{\epsilon, Loc,\delta}(f)$. It comes  with a bundle ${ {\cal L}}^0=\cup_{\delta}{ {\cal L}}_{ Loc, \delta}(f) =\cup_{\alpha,\delta}{ {\cal L}}_{{\alpha}, Loc, \delta}(f_{\alpha})$ and a section $s^0 :{{W}}^0_{\epsilon}(f)\rightarrow { {\cal L}}^0 $ that is the restriction of ${\bar \partial}_J$-section $s$. We defined a new "product" chart, denoted by  $({W')}^0_{\epsilon}(f),$  by the same formula  but  denoting  all objects involved by the same letters with a "prime". Start with the definition of a lowest stratum, ${{W'}}^{\alpha_0}_{\epsilon, Loc,\delta}(f).$ It is defined to be the collection of all
elements $g'.$  Each $g'$ is obtained from an element $g$ in ${{W}}^{\alpha_0}_{\epsilon, Loc,\delta}(f)$ defined as following: (i) on  each  small disc of radius $\rho<<\delta$ centred at one of the double points, $g'$  is a constant function equal to the value of $g$ at the double point; (ii) away from these discs of radius $2\rho, $  $g'=g;$ (iii) on the neck areas, $g'$ is obtained by joining  the constant function with $g$  using a fixed cut-off function supported in the neck area. Note that $\rho=\rho(f)$ needs to be chosen sufficiently small so that $g'$ is still in the neighbourhood of $f$ required in the main gluing estimate above (or see [LT]).  For a normal gluing parameter $\alpha$ with $|\alpha|<< \rho, $  $({W}')^{\alpha}_{\epsilon, Loc,\delta}(f)$ is defined to be the collection of maps $g'_{\alpha}:\Sigma\rightarrow M$,  each  coming   from a map $g':\Sigma_{\alpha_0}\rightarrow M$ in ${{W'}}^{\alpha_0}_{\epsilon, Loc,\delta}(f).$ More specifically, since $g'$ is constant on those $\rho$-discs centred at double points, for $|\alpha|< \rho/2, $  $g'_{\alpha}$ is  defined  to be $g'$ on
the part of $\Sigma _{\alpha}$ that corresponds to $\Sigma _{\alpha_0}\setminus \{$   $\rho/2$-discs  $\}$ and constants on the rest.

Now recall the assumptions  (A1) and (A2) that $D(s_0)_f=D({{\bar \partial } }_J)_f$  has no trivial kernel and cokernel. As far as the smooth structure and associated local Fredholm theory  around a  flat end are concerned, this is  the essential case although from the point view of regularizing the moduli spaces of $J$-holomorphic stable maps the situation here is trivial.

Under this assumption,  let  $\xi$ be an element  in  the central fiber ${ {\cal L}}_{{\alpha_0}, Loc, \delta}(f_{\alpha_0})$ 
and $\alpha=(0, \alpha_{n})$ be a normal gluing parameter.  Then there is an unique 
$g_{\alpha}\in {{W}}^{\alpha}_{\epsilon, Loc,\delta}(f)$ such that $s_0(g_{\alpha})=(\xi, \alpha)$ (or $[s]_0(g_\alpha)=\xi$).  Similarly, there is an unique 
$g_{\alpha_0}\in {{W}}^{\alpha_0}_{\epsilon, Loc,\delta}(f)$ such that $s_0(g_{\alpha_0})=(\xi, \alpha_0)$ (or $[s]_0(g)=\xi$). Recall that $(\xi, \alpha_0)$
and $(\xi, \alpha)$ are the tautological coordinate of $g_{\alpha_0}$ and  $g_{\alpha}$ respectively.
Then the above process given rise the corresponding $g'_{\alpha_0}$ and $g'_{\alpha}$.

Now the key point is that there is a gluing process similar to the one in [L0]. But in the gluing here,   the linear model ${{V'}}^{\alpha}_{ Loc}(f_{\alpha})$ for  ${{W'}}^{\alpha}_{\epsilon, Loc,\delta}(f)$, an open ball in  an infinite dimensional Hilbert space plays the role
of the open ball in the kernel $K_f$ of $D(s_0)_f$ before.  
The equation to be solved is $\pi_{\alpha}^{\perp}\circ [s_{\alpha}] (\xi+\xi^{\perp})=0.$
Here $\pi^{\perp}_{\alpha}$ is  the projection to $({ {\cal L}}_{{\alpha}, Loc, \delta}(f_{\alpha}))^{\perp}$ for the orthogonal  decomposition of $${ {\cal L}}_{\alpha}(f_{\alpha})={ {\cal L}}_{{\alpha}, Loc, \delta}(f_{\alpha})\oplus ({ {\cal L}}_{{\alpha}, Loc, \delta}(f_{\alpha}))^{\perp};$$ and $\xi+\xi^{\perp}$ are the two components with respect to the decomposition of $L_k^p(f_{\alpha}^*TM)=
{{V'}}^{\alpha}_{ Loc}(f_{\alpha})\oplus ({{V'}}^{\alpha}_{ Loc}(f_{\alpha}))^{\perp}.$

The gluing process  here  gives rise a family of  maps parametrized by ${\alpha}$, $\Phi_{\alpha}:{{W'}}^{\alpha}_{\epsilon, Loc,\delta}(f)\rightarrow {{W}}^{\alpha}_{\epsilon, Loc,\delta}(f)$. This gives an other coordinate chart for 
${{W}}^0_{\epsilon}(f)=\cup_{\alpha,\delta}{{W}}^{\alpha}_{\epsilon, Loc,\delta}(f)$.

Now come back to  $(g'_{\alpha}, \alpha)$  constructed before from $(g_{\alpha}, \alpha)$ that corresponds to $(\xi, \alpha) $  with $\xi \in  { {\cal L}}_{{\alpha_0}, Loc, \delta}(f_{\alpha_0})$ and 
$\alpha=(0, \alpha_{n})$.  By the construction, $\Phi_{\alpha}(g'_{\alpha}, \alpha)=
(g_{\alpha}, \alpha)$. It is easy to see that under new natural coordinate with product structure, along a $\alpha$-curve $\alpha \rightarrow (g'_{\alpha}, \alpha)$ constricted above from a fixed $\xi\in { {\cal L}}_{{\alpha_0}, Loc, \delta}(f_{\alpha_0})$ the value of the section $s_{0}$ is just $\xi_0$ after  the natural identifications of the fibers so that the partial derivatives  of ${s}_0$ is identically equal to zero as expected.

This proves that the in the new natural  coordinate with product structure, $s_0$ is of class at least $C^1$ for the case  that the assumptions (A1) and (A2) are satisfied. The more details on this  as well as treatment for the general  case with be given in [L4].

\medskip
\vspace{2mm}
\noindent 
$\bullet$  $\bullet$ $\bullet$  Higher level natural coordinate charts and push-forward of bundles:

Still work locally on a fixed flat chart.
We have  the "natural" chart $^{\kappa}V^{(\alpha, m)}_{ Loc,\delta}$ for a slice of the flat chart which is realized as a graph of a map 
from a ball of $T_{\alpha}$ to its orthogonal complete $O_{\alpha}$ in 
$ L_{(k;\kappa)}^p(\Sigma_{\alpha}, f_{\alpha}^*(TM); h, S^1_{\alpha}) $. 

In particular,  any elements in $^{\kappa}V^{(\alpha_0, m)}_{ Loc,\delta} $ or $^{\kappa}V^{(\alpha, m)}_{ Loc,\delta}$
 are  $J$-holomorphic on the the ${\delta}$-discs or ${\delta}$-annuli around double points.  Therefore, any such elements are determined by their restrictions on the fixed part $\Sigma^0(\alpha)=\Sigma^0(\alpha_0)$. Recall that
$\Sigma^0(\alpha)$ and $\Sigma^0(\alpha_0)$  are defined by removing
  the ${\delta}$-annuli or ${\delta}$-discs from $\Sigma_{\alpha}$ and $\Sigma_{\alpha_0}$ respectively, and are identified by the 
 pre-gluing. Of course, the surface  $\Sigma^0(\alpha)=\Sigma^0(\alpha_0)$ depends on $\delta$. We denote it by $\Sigma^0(\delta).$

 Now fix a $\delta$ and consider $\Sigma^0(\gamma)$
for all   $0<\gamma<<\delta $  as a trivial family  of open curves inside the universal curve   ${\cal U}$ locally over $\Lambda (\alpha_0).$ In fact, we should consider all the open surfaces "generated" by these $\Sigma^0(\gamma)$. But we will suppress this point in the discussion here.
 
As before by using exponential coordinate charts at $f_{\alpha_0}$ and $f_{\alpha}$ and abusing notations,  each elements in $$ L_{(k;\kappa)}^p(\Sigma^0(\gamma), f_{\alpha}|_{\Sigma^0(\gamma)}^*(TM); h)$$ 
  or  $$ L_{(k;\kappa)}^p(\Sigma^0(\gamma), f_{\alpha_0}|_{\Sigma^0(\gamma)}^*(TM); h) $$ can be thought as a $L_k^p$-map from the fixed  open Riemann surface $\Sigma^0(\gamma)$ to $M$. As before for any fixed $\gamma$, 
 with this interpretation here we get a corresponding space of $L_k^p$-maps, denoted by ${{\hat W}}_{\gamma}(f_{\alpha})={{\hat W}}_{\gamma}(f_{\alpha_0})$ near
 the restriction of $f_{\alpha}$ to $\Sigma^0(\gamma)$ which is the same as the restriction of $f_{\alpha_0}.$  As usual,  there  is a   corresponding bundles over ${{\hat W}}_{\gamma}(f_{\alpha})$  whose fibers are the corresponding $L_{k-1}^p$-sections and a smooth ${\bar{\partial}}_J$-section ${\hat s}$. By varying $\gamma$, we get a functorial system of bundles together with  specified sections.
  In general, this system here is not very useful since ${\hat s}$ is not Fredholm any more. 
 
 However,   since   any element in flat chart $^{\kappa}V^{(\alpha, m)}_{ Loc,\delta}$ is determined by its values on any ${\Sigma^0(\gamma)}$,  we get an obvious embedding of $^{\kappa}V^{(\alpha, m)}_{ Loc,\delta}$ into ${{\hat W}}_{\gamma}(f_{\alpha})$ for any $\gamma<<\delta.$ Moreover,  the bundles used to define the flat charts are push-forwarded  into  the corresponding bundles such that the pull-back of $ {\hat s}$ is $s$.  Similarly, the section used to define perturbations can be lifted into sections of these higher bundles. In any case, for the flat charts and related bundles and sections, above seemly useless system  serves as a systems of ambient coordinate charts so that
  one can talk about smoothness in the similar way as the usual definition of  the smoothness for a function defined on a "bad" set inside ${\bf R}^n$ described, for instance at the beginning of Milnor's book [M]. Note that in this setting, we
   have no difficulty to take derivatives along the normal directions.
   
   This brings us a new set of problems to verify the compatibility of  the smoothness of  an  object like moduli space or a section obtained from various different "coordinate"  systems.
   
  Before we close this discuss here, just mention one technical point. Since we  use open Riemann surfaces here, the $L_k^p$-maps usually have some pathology
   behaviour near the boundary, it is better to introduce a modifier, a fixed 
   cut-off function $\beta_{ \delta, \gamma}$on $\Sigma^0(\gamma)$ which is 
   compact supported on  $\Sigma^0(\gamma)$ and equals to 
   $1$ on $\Sigma^0(\delta).$ Then  for any $L_k^p$  map
   $\xi:\Sigma^0(\gamma)\rightarrow M$,  the pseudo-norm 
   $\|\xi\|_{k, p; \kappa; \beta_{\delta, \gamma}}$ is defined to be
    $\|\beta_{\delta, \gamma}\xi\|_{k, p; \kappa}$ and becomes a real norm on the flat chart. This gives a  family of  induced norms  on the flat chart. This essentially completes the construction of the local  flat charts and the local Fredholm theory near the ends, and  hence resolves the difficulty (C).

 We will  generalize the work  on weakly smoothness  in   Sec.2 and Sec.3 to the case of flat charts above in [L2] and [L3]. Combining  these together we get the following results on the smoothness of the moduli spaces in [L2] and [L3].

 \begin{theorem}   
Let $s$ be the ${\bar {\partial }}_J$-section.
Assume that all isotropy groups are trivial.  Let $B(r)=\oplus_{i\in I}B(r_i)$  be the ball of "radius" $r$ in ${\cal R}=\oplus_{i\in I}{\cal R}_i$=collection of all compatible perturbations ${\nu}=\{\nu_i, i\in I\}$.  Here  $\nu_i\in B(r_i)$   is a section of the bundle over a local flat chart $W_i={{W}}_{\epsilon_i, Loc,\delta_i}(f_i)$, supported on   $W_i, i\in I$ and  $\{W_i, i\in I\}$ covers ${\cal M}(A, J)$. These sections  ${\nu}=\{\nu_i, i\in I\}$ are  obtained from the cokernel  either  (i) as "constant" extensions given by [LT],
or (ii) as  "geometric"  extensions given by  this paper. In the first case, they are weakly smooth of class $C^1$, while  in the second
 case, they are of class $C^{m_0}$.

Then for $r$  small enough, each  extended  local moduli space $E{\cal M}^{i}(J, A)=(s+ev_{(i, B(r))})^{-1}(0)$  is  either (i) a submanifold of class $C^{1}$ in ${\cal R}\times W_i$, or (ii) a submanifold of class $C^{m_0}.$

Moreover, these local extended  moduli spaces patch together to form a  manifold, denoted by $E{\cal M}(J, A)$ of class $C^1$ or $C^{m_0}$ respectively.

\end{theorem} 

 As mentioned before to get the required smoothness for the perturbed moduli space, we have to deal with the above two cases separately.
 
 \medskip
\vspace{2mm}
\noindent 
$\bullet$  $\bullet$ Case (i): 

In this case, Sard theorem is not applicable, we give 
$E{\cal M}(J, A)$ a compatible $C^{\infty}$ structure first,   then deform the $C^1$ projection map $\pi':E{\cal M}(J, A)\rightarrow {\cal R}$ slightly into a $C^{\infty}$-map $\pi$.

Then we have the following theorem for case (i).

\begin{theorem}   
For generic choice of $\nu$ in $B(r)$ with $\|\nu\|$  small enough, the perturbed moduli space ${\cal M}^{\nu, \pi}(J, A)=(\pi)^{-1}(\nu)$
is a compact $C^{\infty}$-manifold whose cobordism class is well defined independent of  all the choices  made.

 \end{theorem} 

 The proof of this  is in [L2].

\medskip
\vspace{2mm}
\noindent 
$\bullet$  $\bullet$   Case (ii): 
 
We may assume that $m_0>$ the index of $Ds_f.$ Then Sard's theorem is applicable to $\pi'$. We have the following theorem in [L3].
\begin{theorem}   
For generic choice of $\nu$ in $B(r)$ with $\|\nu\|$  small enough, the perturbed moduli space ${\cal M}^{\nu}(J, A)=(s+\nu)^{-1}(0)$
is a compact $C^{ m_0}$-manifold.
\end{theorem}

 In this paper we will also discuss in detail a few  topics that  may be  considered  as  "side" issues since they are not directly related to the main theme, the lack of differentiability  discussed here. 
These include   the  proof of the compactness of the 
perturbed moduli space; the existence of a positive lower bound  for   $\|{\bar{\partial}}_J\|_{k-1, p}$ on the boundary of a $L_k^p$-neighbourhood of the moduli space inside the ambient space; the sufficiency for transversality of the extended and perturbed moduli space inside the $L_k^p$-neighbourhood when the space of perturbations   is  only of finite dimensional coming from  the localized geometric sections derived from the cokernels. In particular, for the last  issue, unlike polyfold theory where there are  much more  abstract perturbations  that can be used to achieve transversality,  there is a possibility  that in our case, we may get into Zeno type paradox.    Perhaps all these issues are  obvious to the experts. However, as  these   "side" issues are quite of general nature,        we    give a complete treatment  on them in this paper.  Of course, there are overlaps of  the discussions here with [LT], and our treatment here is  not necessarily simpler mainly  because of  the last side issue.  
Note that while  most of the results  in Sec. 3  are only applicable for the moduli space in GW and Floer type theories with only one stratum, our treatment for these side issues in Sec. 3  works  for the general cases.

   This paper is organized as following.  
   
    Sec. 2 introduces the notion of  weak smoothness as a general frame work to  overcome the difficulty (A).
   
   Sec.3     establishes the  $C^{m_0}$-smoothness of the perturbed moduli space  stated in the main theorem. 

  The argument in this section with some modifications is applicable to the  other cases in GW and Floer type theories. In fact, after a quick review on  some of  notations in GW theory,   
  we  immediately  start to use "generic"  terminologies and notations in GW and Floer theories.  Therefore, almost all  the arguments and results are applicable to  the part of moduli spaces in $GW$ and Floer theories  that have a fixed stratum  with trivial  isotropy  groups.

   
    Sec. 4 collects some well-known analytic facts. In particular, we include an elementary proof of the  smoothness of the $p$-th power of $L_k^p$-norms of the Sobolev spaces when $p$ is an even positive integer.


 This paper and subsequent ones are written based on author's notes. The electronic files  of part of these notes  are listed as [C] and [L1] in the reference. 
 All ideas and almost all results in this paper are already written  in those notes.

 \section{Weakly Smooth  Banach  Manifolds and Bundles}  
  \medskip
\vspace{2mm}
\noindent
 In this section, we  introduce the notion of weakly smoothness structure associated with  a covering  on a topological Banach manifold.  
  We start with a topological Banach manifold ${\cal B}$  and Banach bundle ${\bf L}\rightarrow {\cal B}.$

 
 Fix   a coving ${\cal U}=\{ U_i, i\in I\}$ of ${\cal B}$ and the collection of  locally trivialized  bundles $\{ {\bf  L}_i\rightarrow U_i, i\in I\}. $  They will be called a collection of admissible charts and trivializations.    We 
define a function on an open set of ${\cal B}$ to be weekly smooth with respect to ${\cal U}$ if  it is smooth viewed in any admissible coordinate charts of ${\cal U}$.  The germs of such  weakly smooth functions gives rise the weakly smooth structure ${\cal O} $ on ${\cal B}$.    The weakly smooth sections  of ${\bf L} $ can be defined similarly  by using the local bundles ${\bf L}_i$.

The real question is if this  obvious  notion of weak smoothness is useful.

 We will show in this and subsequent papers that  in  the case of  GW and Floer type theories, the stratified topological manifolds appeared there  do have enough smooth functions and section so that they behave as if they are honest  smooth manifolds. Moreover,     any finite dimensional weakly submanifold in such a manifold with induced weakly smooth structure is  in fact a honest smooth  one.  This gives  a way  to  have perturbed  $C^{m_0}$-smooth moduli spaces in GW and Floer theories that do not lie inside the space of smooth stable maps, as we mentioned in the introduction.

   
There are  several  immediate  questions on the above definition of weakly smooth functions and sections associated with the covering ${\cal U}:$  (I) To what extend does this notion depend on covering ${\cal U}$?  Since we only have a topological Banach manifold, we expect that  there are not too many  weakly smooth functions with  respect to the fixed covering. It is desirable to have more weakly smooth functions. On the other hand, replacing  ${\cal U}$ by a compatible refinement does   give rise more smooth functions with respect to new covering.
This suggests that at least one should incorporate  the effect of compatible refinement in order to formulate the notion of weakly smoothness. (II) What is the functorial behaviour of this notion? Is it possible to form a category of weakly smooth Banach manifolds so that the usual functorial constructions in smooth category can be performed?

The answer to  these and related question motivates the definition below.

 ${\bullet}$ Weakly   smooth Banach  manifolds  and   bundles:
 
 We mimic  the usual way to use ringed space to specify analytic or algebraic-geometric structures on  topological spaces by gluing
 the corresponding local structures. We  start with  the formal definition for the weak smooth structure on 
 ${\cal B}$.
  Then  we define  the notion of weakly smooth sections  of  subbundles of ${\bf  L}$. Those subbundles will be used as obstruction bundles to perturb the Fredholm section $s$ to achieve the transversality.  
 
 \noindent 
Let ${\cal B}$ be a paracompact  topological Banach manifold locally modelled on a separable Banach space $E$. Assume that it is   covered by open sets ${\cal B}=\cup_i U_i$ with  coordinate chart $\phi_i:U_i\rightarrow W_i, $   where $W_i$ is a open set of $E$.

The collection $({\cal U}, \Phi)=(\{ U_i, i\in I\},\{ \phi_i:U_i\rightarrow W_i, i\in I\})$  forms  the fixed covering data for ${\cal B}.$

On each $U_i, i\in I$, we consider the "structure sheaf"  ${\cal O}_i$ which is "admissible" by the "induced"covering  data of $({\cal U}, \Phi)$ on $U_i$. More specifically, each ${\cal O}_i, i\in I$  is    a sheaf of   subrings   consisting of germs of  smooth functions on $U_i$ with respect to the "induced" covering data. Here  the "induced" covering data on $U_i$, denoted by $({\cal U}, \Phi)(U_i)$ is just the restriction of $({\cal U}, \Phi)$
to $U_i.$
 As usual,  the choices of the subrings depend on (or give) the "structures" that we want to define associated to the covering.  For our purpose here, the obvious canonical choice for  ${\cal O}_i$ is  the one consisting of germs of all consistent
  smooth functions on   $U_i$. 
  More precisely, any element in ${\cal O}_i$
defined on  an open subset $W_i$ of $U_i$   is  weakly smooth if and only if it is smooth viewed  on  any open subset $ W_i\cap U_j, j\in I$ after composing with the corresponding transition function. Strictly speaking, we should write the  structure sheaf ${\cal O}_i$ on $U_i$ as ${\cal O}_i(({\cal U}, \Phi)(U_i))$ to indicate its dependence on the induced covering data.  Once each structure sheaf ${\cal O}_i$ is obtained,  the structure sheaf on  ${\cal B}$ associated to the covering data $({\cal  U},\Phi)$, denoted by  ${\cal O}_{\cal  U}$, can be obtained by the  "gluing"  these ${\cal O}_i, i\in I$  together in the usual manner.


Now assume that ${\cal U'}$ is a "compatible"  refinement of the covering ${\cal U}$ given by $i:{\cal U'}\rightarrow  {\cal U}$ in the sense that each coordinate chart
$\phi'_j:U'_j\rightarrow {\cal B}$ is the restriction of $\phi_{i(j)}:U_{i(j)}\rightarrow {\cal B}.$
Clearly   the obvious pull-back map $i^*:{\cal O}_{\cal  U}\rightarrow {\cal O}_{\cal  U'}$ makes ${\cal O}_{\cal U}$   a  subsheaf of ${\cal O}_{\cal  U'},$ and it is a strict subsheaf in general. In other other
words, there are more weakly smooth functions  with respect to ${\cal U'}$ in general.

 The collection  of all possible coverings of ${\cal B}$ that are compatible refinements of a fixed ${\cal U}$ with respect to the partial order given by inclusion  form a "directed" set and the collection of  locally finite ones are cofinal. We denote the collection by $[{\cal U}].$
We define the structure sheaf ${\cal O}_{[{\cal U}]}$ of the weakly smooth structure  on ${\cal B}$ to be the inverse limit of ${\cal O}_{\cal U}$ with respect to all ${\cal U}$-compatible locally finite coverings on ${\cal B}.$ In summary,  a element in  ${\cal O}_{[{\cal U}]}$ is represented by a function $f$  defined on an open sect of ${\cal B}$ such that there is a locally finite covering which is a compatible refinement of 
 ${\cal U}$ so that $f$ is smooth after pulling back to each admissible  coordinate chart of the covering.

The  weakly smooth structure ${\cal O}$ on ${\cal B}$ so defined are obtained by selecting all possible compatible smooth  sections of the trivial bundle ${\cal B}\times {\bf R}^1\rightarrow {\cal B}. $  We define the notion of weakly smooth sections for a topological Banach bundle
${\bf  V}\rightarrow {\cal B}$ in  a similar way. Here we assume that locally over  each $U_i, i\in I, $ ${\bf V}_i={\bf V}|_{U_i}$ has a fixed  trivialization so that it makes sense to talk about local smooth sections of ${\bf V}$ over $U_i, i\in I$ with respect to the  covering ${\cal U}$ and local trivializations.  Let ${\cal V}_{\cal U}$ be the  maximum of  compatible 
 ${\cal O}_i$-modules ${\cal V}_i, i\in I $, where each ${\cal V}_i$ is   a sheaf of submodules over ${\cal O}_i$ consisting of the germs of smooth sections of ${\bf V}_i.$ Any  global section of ${\cal  V}$ over an open set $W$   is  called weakly smooth if  it is a smooth section   viewed in any "admissible" local bundles ${\bf V}_i, i\in I$ for some  compatible covering in $[{\cal U}].$

As before, strictly speaking in the above definition of the  ${\cal O}_{\cal B}$-module,  ${\cal V},$ in addition to the covering data on ${\cal B}$,  we should   also specify more carefully on the  data for bundles.

Here is  another reason that   we need to introduce
the compatible  refinements to a covering.

 A continuous map $ F$ between $({\cal B}', {\cal U}', { \Phi}')$ and $({\cal B}, {\cal U}, { \Phi})$ is said to be weakly smooth
 if $F^*$ maps  ${\cal O}_{[{\cal U},{ \Phi} ]}  $ into the corresponding one with respect to the compatible refinement of 
 $({\cal U}',{ \Phi}') $ by intersecting it  with $F^{-1}{\cal U}.$  Clearly, $F$ induces a map on the inverse limits. Therefore, the collection 
 $({\cal B}, {\cal O}_{[{\cal U},{ \Phi} ]})$ with the morphism just defined form a category of weakly smooth Banach manifold with a 
 specified compatible  weakly smooth charts $[{\cal U},{ \Phi} ]$. The equivalence between two such objects are defined in the obvious 
 way. In particular, each ${\cal O}_{[{\cal U},{ \Phi} ]}  $ determine a equivalent class, the weakly smooth structure  associated  to 
 the covering data  $[{\cal U},{ \Phi} ]$. In other words, when the covering data $[{\cal U},{ \Phi} ]$ of ${\cal B}$ is fixed, the 
 effect of  equivalent class  is to collecting all possible "admissible" coordinate charts  not just the ones that are  compatible refinements of the given covering data.

 This can be done similarly for the $({\bf V}\rightarrow {\cal B})$
  to get the corresponding category of ${\cal O}$-modules. The functorial nature of these constructions makes it possible to define various
  familiar  notions. We mention a few relevant ones.

  There is a well defined cotangent functor and hence, a corresponding tangent functor, in the category  weakly  smooth Banach manifolds.  On any  $U_i$ of the admissible cover of ${\cal B},$ the  the cotangent sheaf $T^*U_i$ is simply defined to be the ${\cal O}|_{U_i}$-modules generated by the derivatives $d\phi,$ for a   weakly smooth function $\phi \in {\cal O}(W), $ where $W$ is an open sub set of $U_i$. Unlike the  category of smooth Banach manifolds, in which the cotangent functor gives rise Banach bundles, in the case of   weakly  smooth Banach manifolds, 
 the cotangent functor only gives a sheaf of modules over ${\cal O}_{\cal B}$ in the sense defined above in general.

 To justify our definition,  we will call  a weakly smooth structure on ${\cal B}$ is effective if any finite dimensional topological
 submanifold of ${\cal B}$ with the induced weakly smooth structure from ${\cal B}$  is a honest smooth manifold. In other words, 
 if a weakly smooth structure is effective, then there are enough smooth  functions on ${\cal B}$ to detect a finite dimensional object in ${\cal B}$. In next section, we will show that the weakly smooth structures appeared in GW  are effective in  the context of this  paper.  
 
 Therefore, if the weakly smooth structure on ${\cal B}$ is effective, then  for any finite dimensional topological, hence smooth  submanifold $M$  of ${\cal B}$,
the restriction the cotangent bundle to $M$ is just the usual cotangent bundle $T^*M$. Of course, as usual, once $T^*{\cal B}$ is defined, the  sheaf of differential
 forms $\Omega^*({\cal B})$ is well defined. It is again a functor on the category of weakly smooth Banach manifolds. Consequently,   integration of differential forms over a finite dimensional submanifold in ${\cal B}$ is well defined.

 
%


 \section{Weakly Smooth Structure on the Ambient Space  of the Moduli Space of  $J$-holomorphic Maps.}  
  \medskip
\vspace{2mm}
\noindent

 Let $(M, \omega)$ be a  compact symplectic manifold of dimension $2n$ with $\omega)$ compatible  almost  complex structure $J$. Denote the associated metric
 by $g_J=\omega(-, J-).$
Assume that all geometric data above are of class $C^{\infty}.$ Fix an effective second homology class $A\in H_2(M, {\bf Z})$. Let $(\Sigma, i)={\bf P}^1$  with its standard complex structure and metric.



 
 Recall that a map $f:(\Sigma, i)\rightarrow (M, J)$ is said to be $J$-holomorphic of class $A$ if (1) the homology class of $f$, $[f]$ is equal to $A$; (2) the map $f$ satisfies $ {{\bar \partial }}_{J}f=0$. Here  the operator $ {{\bar \partial }}_{J} $ is defined by 
 $ {{\bar \partial }}_{J}u=d u+J(u) \circ du \circ i .$

Let ${\tilde {\cal M}} (A)$ be the collection all such  $J$-holomorphic maps. The virtual dimension of ${\tilde {\cal M}} (A)$
 can be calculated by the index of the linearization 
$D{\bar {\partial }}_{J}|_u$ at any $u\in {\tilde {\cal M}} (A)$ and is equal to $2c_1(A)+2n$.   The group  $G={\bf PSL}(2, {\bf C})$  acts on ${\tilde {\cal M}}(A)$  as the group of reparametrizations. Under the assumption that $A$ is effective, each element in ${\tilde {\cal M}}(A)$ is stable in the sense that it has no infinitesimal automorphism. Consequently, in this case,  the action $G$  is proper  with finite isotropy group. For the purpose of this paper, we assume that the isotropy groups are trivial for all elements in ${\tilde {\cal M}}(A)$.   Denote  the quotient space of the unparametrized  $J$-holomorphic maps of class $A$ 
 by ${\cal M} (A)$.  Assume that  it is compact. This  is the moduli space that we want to regularize by perturbing the defining  section $ {{\bar \partial }}_{J}$. 
 We will outline a proof  in this section that in this situation, the resulting perturbed moduli space is a $C^{m_0}$-smooth manifold, where $m_0=k-\frac{2}{p}$    is the  Sobolev  differentiability  of the ambient space ${\cal B} (A)$ of stable  $L_k^p$-maps of class $A$.

Our next goal is to describe the weakly smooth structure on ${\cal B} (A)$ and to show that there are enough weakly smooth functions to make it effective.   To this end, we need  to give more details on ${\cal B} (A)$ first.  To simplify our notation, from now on, we will  drop the $"(A)"$ in all the notations and write ${\cal B} (A)$ as ${\cal B} ,$ etc.

Note that  each $u$ in ${\tilde {\cal M}} $ has  some   point $x_0\in \Sigma$  such that $u$ is a local embedding near $x_0.$ 
This  actually implies that $u$ has  no infinitesimal  automorphisms and hence stable.  It also implies that    we can find a local  slice for  the $G$-action near $u$  by using a local hypersurface $H_{u}$ of codimension 2 transversal to the local  image of $u$ near $x_0$. 


Now embed ${\tilde {\cal M}}$ into ${\tilde {\cal B}}, $ the space of parametrized stable $L_k^p$-maps of class $A$. 
Here ${\tilde {\cal B}}={\tilde {\cal B}}_{k,p}(A)$ is defined to be the collection of all maps $u:\Sigma\rightarrow M$
of class $L_k^p$  measured by the metrics $g_J$ on $M$ and the standard one on $\Sigma$ such that (i) $[u]=A;$ (ii) $u$ is stable.
 In [L1], the stability for $L_k^p$-maps  is formulated in more general setting.
Here for simplicity, we  assume that  each element $u \in {\tilde {\cal B}},$  like $J$-holomorphic ones,  has a local embedding   
point $x_0$.  Again, we will assume that the stabilizer of $u$ is trivial. These assumptions  imply that the $G$-action on  ${\tilde {\cal B}}$ has local slices.



 

 

 Here we assume that the $m_0=k-\frac {2}{p}>1 $ so that each element in 
${\tilde {\cal B}}$ is at least of class $C^1.$ Another assumption that we use in rest of the paper is that $p$  is a positive even integer. This assumption implies that  
 the $p$-th power of the $L_k^p$-norm is a smooth function on $L_k^P(\Sigma, f^*(TM)).$ 
 

It is well-know that
when the Sobolev index $m_0= k-\frac{p}{n}$ is greater than zero, the notion of an element  $f$ being  in $L^p_{k,loc}({\bf R }^n, {\bf R}^m)$ is invariant with respect the local differomorphisms of both the domain and target. This implies that  $ {\tilde {\cal B}}$ is well-defined.
 To specify the topological and smooth structures  on  $ {\tilde {\cal B}}$, we introduce local  coordinate charts for $ {\tilde {\cal B}}$.

 
  \medskip
\vspace{2mm}
\noindent

$ {\bullet}$ The local coordinate of ${\tilde {\cal B}}$:

 For any $f$ in ${\tilde {\cal B}}$, when $\epsilon$ is small enough,
 we have a local coordinate
chart
$Exp_f: U_{\epsilon}(f)\rightarrow {\tilde {\cal B}}$ for an $\epsilon $-ball $ U_{\epsilon}(f)$ in $L_k^P(\Sigma, f^*(TM)).$ Here $Exp_f$ is defined to be $Exp_f(\xi)(x)=exp_{f(x)}\xi(x)$ for $\xi\in  U_{\epsilon}(f)$ and $x\in \Sigma$, and $exp_y :T_y M\rightarrow  M$ is the exponential map of $M$.
It is well-know in Gromov-Witten and Floer theories that smoothness of exponential map in $M$ implies  that the transition  functions  between  these coordinate charts  are $C^{\infty}$ even although the elements in ${\cal M}_{k,p}(\Sigma, M)$ are only in  $L_k^p$.
 This makes ${\tilde {\cal B}}={\tilde {\cal B}}_{k}^{p}(\Sigma, M)$  a $C^{\infty}$  Banach manifold.

Next Lemma is  well-known   in Gromov-Witten theory. We include  the simple  proof here
for completeness.


\begin{lemma}\label{first}
For any fixed $x=(x_1, \cdots, x_l)$ in $\Sigma^l$ with $x_i\not= x_j$ for $i\not =j, $ define the evaluation map $e_{x}:{\tilde {\cal B}}\rightarrow
M^l$ by $e_x(f)=(f(x_1), \cdots f(x_l))$ for any $f\in {\cal M}.$  Then $e_x$ is a $C^{\infty}$ submersion.

\end{lemma}

 \proof

Recall that for any $f\in {\tilde {\cal B}}$,   we have the coordinate chart $Exp_f:  U_{\epsilon}(f)\rightarrow {\tilde{\cal B}}$ for an $\epsilon $-ball $ U_{\epsilon}(f)$ in $L_k^P(\Sigma, f^*(TM))$,   and $Exp_f(\xi)(x)=exp_{f(x)}\xi(x)$ for $\xi\in  U_{\epsilon}(f)$ and $x\in \Sigma.$ Similarly we introduce the coordinate chart $exp_y:  U_{\epsilon}^l(y)\rightarrow { M}^l$ for an $\epsilon $-neighbourhood  $ U_{\epsilon}^l(y)$ of the zero in $T_yM^l. $  It is defined by $exp_y(\eta)=(exp_{y_1}(\eta_1),\cdots, exp_{y_l}(\eta_l))$ for any $\eta=(\eta_1, \cdots, \eta_l)$ in $T_yM^l $ with $y=(y_1, \cdots, y_l)$ in $M^l.$ With respect to these coordinate charts of ${\tilde{\cal B}}$ and $M^l$, with $y=f(x)$, the evaluation map $e_x$ has the form $exp_{f(x)}^{-1}\circ e_x\circ Exp_f(\xi)=exp_{f(x)}^{-1}(exp_{f(x_1)} \xi (x_1), \cdots, exp_{f(x_l)} \xi(x_l))=(\xi(x_1), \cdots,  \xi(x_l))$. It is induced from a linear  map  from  $L_k^P(\Sigma, f^*(TM))$ to $T_{f(x)}M^l$.
Under the assumption that $m_0>0$, it is continuous linear map, hence smooth. The surjectivity of $De_x$ at $f$ follows from the fact that
for any $\gamma=(\gamma_1, \cdots, \gamma_l)$ with $\gamma_i\in T_{f(x_i)}$  $i=1, \cdot, l,$ there exits a $\xi$ in  $L_k^P(\Sigma, f^*(TM))$ such that $\xi(x_i)=\gamma_i$ for $i=1, \cdots, l.$

\QED

\begin{lemma}\label{second}
For any fixed $g\in G,$  $\Psi_g:{\tilde { \cal B}}\rightarrow {\tilde {\cal B}}$ is smooth. Here $\Psi_g$ is  the action map for a fixed $g\in G$ and defined 
by $\Psi_g(f)= h\circ g.$

\end{lemma}

\proof

Let $h=Exp_f(\xi)$ for $\xi $ in  $U_{\epsilon }(f) \subset L_k^P(\Sigma, f^*(TM)).$ Then $\Psi_g(h)=h\circ g=Exp_{f\circ g}(\xi\circ g)=Exp_{f\circ g} g^*(\xi)$ in $U_{\epsilon' }(f\circ g)\subset L_k^P(\Sigma, (f\circ g)^*(TM)).$ Therefore,  in the local coordinate charts
$U_{\epsilon }(f)$ and $U_{\epsilon' }(f\circ g),$  $\Psi_g$ can be identified with ${\tilde \Psi_g}: L_k^P(\Sigma, f^*(TM))\rightarrow  
L_k^P(\Sigma, (f\circ g)^*(TM))$ defined by $ {\tilde \Psi_g}(\xi)=g^*(\xi), $  which is just the pull-back on sections. Clearly  ${\tilde \Psi_g} $ is linear and continuous, and hence smooth.

\QED

\begin{pro}\label{first}
The action map $\Psi:G\times {\tilde {\cal B}} \rightarrow {\tilde {\cal B}}$ composed with the evaluation map $e_x$, denoted by
$\Phi_x: G\times {\tilde {\cal B}} \rightarrow M^l$, is of class $C^{m_0}$.

\end{pro}

\proof

This  follows from the $sc^k$-smoothness of the action map. A direct proof was given in [L1] and will be given in [L3]. 

 Note that it follows from above two Lemmas that $\Phi_x$  is $C^{\infty}$ along ${\tilde {\cal B}}$-direction and is at least of class $C^{m_0}$ along $G$ direction by Sobolev
embedding theorem. The only question is about the  smoothness of the mixed derivatives.  This can be reduced  as follows. Note that $\Phi_x: G\times {\tilde {\cal B}} \rightarrow M^l$ is  the composition of the follow tow maps: $\Psi_x: G\times {\tilde {\cal B}} \rightarrow \Sigma^l \times {\tilde {\cal B}}$  given by $\Psi_x:(g, f)= (g(x), f)$ and $ev^l: \Sigma^l \times {\tilde {\cal B}} \rightarrow M^l$ given  by $ev^l(y, f)=f(y).$
The first map is obvious smooth. The total evaluation map  $ev^l$ is $C^{\infty}$-smooth along $  {\tilde {\cal B}}$-direction,  and  it is $C^{m_0}$-smooth along  $\Sigma^l$-direction by Sobolev embedding theorem. Again the question is about the mixed derivatives.

\QED

  \medskip
\vspace{2mm}
\noindent

$ {\bullet}$ The local slice of the  $G$ action  on ${\tilde {\cal B}}$: 

The argument here is  only outlined. The more details for more general cases were in [L] and will be given in [L3]. For simplicity, we only consider the case that ${ {\cal B}}$ is a small neighbourhood of the moduli space ${\cal M}$ of stable $J$-holomorphic maps.
Given a stable  $J$-holomorphic  map $f$ in  ${\tilde {\cal B}}$ with trivial isotropy group,  consider its $G$-orbit $O(f)$ as a injective map $O_f: G\rightarrow  {\tilde {\cal B}}$. The following Lemma was proved in [L]  and  the detail   will be given in [L3].

\begin{lemma}\label{first}

The condition above   implies that there are open neighbourhoods $V$  of $e$ in $G$ and  $W(f)$ of $f$ in  ${\tilde {\cal B}}$ such that for any point $g$ in  $G\setminus V$, $O_f(g)$ is not in $W(f)$. Here we require that $V$ is pre-compact.

Moreover,   when  $W(f)$ is small enough, for any $h$ in $W(f)$, the same is true for the orbit map $O_h$ with a fixed $V$ independent of $h$.

\end{lemma}

The lemma implies that we  can replace the $G$-space  ${\tilde {\cal B}}$ by  $W(f)$, and   we only need to get local slice for the action from the pre-compact set $V$ on $W(f)$.

To this end, fix a  point $p=(p_1, p_2, p_3)$ on  ${\Sigma}^3$ such that $f$ is a local $C^1$ embedding  near $p_i, i=1, 2, 3$. We assume further that each $p_i$ is a injective point.
 Let $H=H_f=(H_1, H_2, H_3)$  and $H_i$ be a local hypersurface of codimension  2 at $f(p_i)$ transversal to $f$ locally.  The assumption
 implies that when $H$ and $W(f)$  are  small enough, there are open discs $D_i$ centred  at $p_i$ such that for any $h$ in $W(f)$,  (i) $h(\Sigma\setminus D_i)$ does not intersection $H_i$; (ii) the  restriction of $f$  to each $D_i$ intersects with  $H_i$  transversally only at $f(p_i)$.

To define the local slice, note that  since $ev_p:{\tilde {\cal B}}\rightarrow M$ is a $C^{\infty}$ surjective map, it follows from the implicit function theorem that when ${\epsilon}$-small
enough, the inverse image $W_{\epsilon }(f, H)=W_{\epsilon }(f)\cap (ev_p)^{-1}(H)$ is a $C^{\infty}$-submanifold of the open set
$W_{\epsilon }(f)=Exp_f (U_{\epsilon }(f)).$  Since the action of each element  $ g$ in $V\subset G$ on $h$ is determined by $h\circ g (p_i),$ the above conditions imply that $W_{\epsilon }(f, H)$ is indeed a  slice of the $G$-action on $W_{\epsilon }(f)$. Another way to get a local slice is the following.
Note that although the $G$-actions are only continuous, the map $\Phi_p: G\times W_{\epsilon }(f)\rightarrow  M$ is at least of class $C^1$ and is transversal to $H$ at $(e, f)$ even restricting to $G\times \{f\}$ by the local  injectivity assumption on $p$. The same is true for
the restriction of $\Phi_p$ to  $ G\times W_{\epsilon }(f, H)$. By implicit function theorem,  this at least  implies that $W_{\epsilon }(f, H)$ is a slice  for  the action of a small neighbourhood of $G$ on $W_{\epsilon }(f)$ when $\epsilon $ is small enough. The fact that the stabilizer of $f$ is trivial implies that it is a slice for the $G$-action on  $W_{\epsilon }(f)$.

Therefore, the quotient
 space of ${\tilde {\cal B}}$ by the $G$-actions, denoted by ${ {\cal B}}$, locally is  modelled on  an open  ball  $U_{\epsilon }(f, h)$ in the Banach $L_k^p(\Sigma, f^*(TM), h)$.  Here $L_k^p(\Sigma, f^*(TM), h)$ is  the closed subspace of $L_k^p(\Sigma, f^*(TM))$, whose element
 $\xi$  is subject to the condition that $\xi(p)$ is in $h$,  which is the tangent space of $H$ at $f(p)$. 
 More precisely, it follows from the implicit function theorem that in the local coordinate  $U_{\epsilon }(f)$, 
 the inverse image of the local slice $Exp_f^{-1}(W_{\epsilon }(f, H))$ is realized as a graph of a function from the open  ball  $U_{\epsilon }(f, h)$ in $
  L_k^p(\Sigma, f^*(TM), h)$ to  the  $L^2$ orthogonal complement of $
  L_k^p(\Sigma, f^*(TM), h)$ in $
  L_k^p(\Sigma, f^*(TM))$. Note that $L_k^p(\Sigma, f^*(TM), h)$ is the tangent space of $W_{\epsilon }(f, H)$ at $f$.
 We may assume that $H$ are geodesic submanifolds so that $Exp_{f(p)}h=H.$ Then $Exp_f^{-1}(W_{\epsilon }(f, H))$  becomes  the  open ball $U_{\epsilon }(f, h)$ in $L_k^p(\Sigma, f^*(TM), h)$ and $ Exp _f:U_{\epsilon }(f, h)\rightarrow  W_{\epsilon }(f, H)) $ is  a local coordinate chart (slice )  for  ${\cal B}$.


 To see  that ${\cal B}$ is a topological Banach manifold, we need to  find the transition function between two local slices $W_{\epsilon }(f_1, H_1)$ and $W_{\epsilon }(f_2, H_2).$ We will  abuse notation by  writing the transition function as $t_{21}:W_{\epsilon }(f_1, H_1)\rightarrow W_{\epsilon }(f_2, H_2).$

\begin{lemma}\label{first}
There is a $C^{m_0}$-smooth function $T_{21}:W_{\epsilon }(f_1, H_1)\rightarrow G$ such that  $t_{21} (f)=f\circ T_{21}(f).$

\end{lemma}

\proof

For simplicity and without lose the generality, we may assume that $f_1$ and $f_2$ are in the same $G$-orbit with $f_2=f_1\circ g_0. $ In fact, one can reduce further by assuming that $g_0=e$ so that we are looking at the same neighbourhood $W_{\epsilon }(f)$ with two different
slices  centered at $f=f_1=f_2.$
Let ${\bf O}(W_{\epsilon }(f))$ be the open set of the orbit of $W_{\epsilon }(f)$. Then
 $\Psi:G\times W_{\epsilon }(f_1, H_1)\rightarrow {\bf O}(W_{\epsilon }(f))$ is a  homeomorphism, which induces   another (incompatible )    
 smooth structure on $W_{\epsilon }(f). $ By the Lemma above, $\Phi_{p_2}=ev _{p_2}\circ \Psi:G\times W_{\epsilon }(f_1, H_1)\rightarrow  M$ is of class $C^{m_0}.$ With respect this new  smooth structure, $\Phi_{p_2}$ is transversal to $H_2$ along $G$-direction at the point $(f_1, e).$ It follows from the implicit function theorem that the inverse image $\Phi_{p_2}^{-1}(H_2)$ is a  $C^{m_0}$-smooth submanifold of $G\times W_{\epsilon }(f_1, H_1)$ and  it  can be realized as  the graph $\{ (T_{21}(\xi), \xi)|\,\, \xi\in W_{\epsilon }(f_1, H_1)\}  $ for a function  $T_{21}: W_{\epsilon }(f_1, H_1)\rightarrow G.$
  Clearly,  $T_{21}$ is the function that we are looking for, which is of class $C^{m_0}$.


\QED

\medskip
\vspace{2mm}
\noindent
Having given ${\cal B}$ a topological Banach manifold structure, we now specify  a weakly smooth structure on this space.

\medskip
\vspace{2mm}
\noindent
 ${\bullet}$ Weakly Smooth Structure on ${\cal B}$:  
 
 We will take the covering ${\cal U}$ of ${\cal B}$ to be the collection of  all possible local  slices $W_{\epsilon }(f, H)$,
  and define ${\cal O}={\cal O}({\cal U})$  to be generated by the germs all functions which are smooth  viewed in any $W_{\epsilon }(f, H).$ 
  Clearly, for any $x=(x_1, \cdots,  x_l ) $ in $\Sigma^l$, the evaluation map $e_x: W_1=W_{\epsilon_1 }(f_1, H_1)\rightarrow M^l$ is smooth.
  Next lemma shows that it is $C^{m_0}$-smooth
  viewed in any slice $ W_2=W_{\epsilon_2 }(f_2, H_2)$. Therefore, its composition with any smooth function on $M^l$ is a weakly smooth function of class $C^{m_0}$. In particular by using  local coordinate functions on $M^l$, we obtains a collection of such   smooths function associated to each $e_x$. For our purpose, it is more convenient  to regard $e_x$ as a $M^l$-valued weakly smooth function.

  \begin{pro}\label{first}
The evaluation map $e_x: W_1=W_{\epsilon_1 }(f_1, H_1)\rightarrow M^l$ 
  is  $C^{m_0}$-smooth
  viewed in any  other slices.

\end{pro}

\proof

Let $ W_2=W_{\epsilon_2 }(f_2, H_2)$ be another slice and $t_{12}: W_2\rightarrow W_1$ be the coordinate transformation. Then
$e_x\circ t_{12}: W_2\rightarrow {M}^l$  can be written as $e_x\circ t_{12}=e_x\circ \Psi \circ (T_{12}, Id), $ where $\Psi:G\times {\tilde {\cal B}} \rightarrow {\tilde {\cal B}}$ is the action map. We have already proved that $e_x\circ \Psi$ is of class $C^{m_0}.$ Since $
T_{12}: W_2\rightarrow G$  is of class $C^{m_0}, $  so is $e_x\circ t_{12}$.
  \QED
  
  Recall  that  in local coordinates 
  $Exp_f: U_{\epsilon }(f)\rightarrow W_{\epsilon }(f)$ and $exp_{x}:T_xM^l\rightarrow M^l$, the evaluation map $e_x$ has the form $exp_{f(x)}^{-1}\circ e_x\circ Exp_f(\xi)=(\xi(x_1), \cdots,  \xi(x_l))$. That  is a just a collection of linear functionals.
 In particular  the derivative  $D(e_x)|_f $  at  $f$ is given by   $D(e_x)|_f (\xi)=(\xi(x_1), \cdots, \xi(x_l))=0.$

 We now prove one of the main results of this section.

 \begin{pro}\label{first}

 The weakly smooth structure so defined on ${\cal B}$ is effective. That is any  finite dimensional weakly smooth  submanfold N in ${\cal B}$
 with respect to any subcovering ${\cal U'}$ is in fact a $C^{m_0}$-smooth manifold with respect to the induced weakly  smooth structure.

\end{pro}

\proof

 We only give the  proof for the full covering ${\cal U}.$ Our assumption implies that for any point $f$ in $N$, there is neighbourhood 
 $V(f)=N\cap W_{\epsilon }(f, H)$ such that $V(f)$ is a closed smooth submanifold of $W_{\epsilon }(f, H).$ Then $N$ is covered such neighbourhoods. Assume that dimension $N=n.$ 
 Since any weakly smooth function on ${\cal B}$ is smooth when  it is restricted to any of such $V(f)$,  the restrictions of these weakly smooth functions to $N$ give it  a smooth structure if we can show that on each $V(f)$ there are  enough such functions to give $V(f)$ a local coordinate system for sufficient small $\epsilon.$ Using those $M^l$-valued evaluation map, we only need to show that  when  $l$ is large enough, $e_x:V(f)\rightarrow M^l$ is a local embedding for suitable choice of  $x$.  To this end, we only need to show that $D(e_x)|_f: T_f V(f)\rightarrow T_{f(x)} M^l$ is  injective.
 Now $T_f V(f)$ is a $n$-dimensional linear subspace of $T_fW_{\epsilon }(f, H)=L_k^p(\Sigma, f^*TM, h)$  and the derivative $D(e_x)|_f$  is the restriction of the corresponding one of for the evaluation map on the ambient space  $W_{\epsilon }(f, H)$ , which is  given by the formula
 $D(e_x)|_f (\xi)=(\xi(x_1), \cdots, \xi(x_l)).$    When $l=1$,  $x=x_1$ is just a point on ${\Sigma}.$  For proper choices of $x_1$, the condition that $\xi(x_1)=0$ will cut the dimension of $T_fV(f)$ at least  by one. Therefore for the proper choices of $x$ with $l\geq n$, 
 the equation $D(e_x)|_f (\xi)=(\xi(x_1), \cdots, \xi(x_l))=0$ can only have the trivial solution. This implies that $D(e_v)|_f$   is  injective.
  
 \QED

 \medskip
\vspace{2mm}
\noindent
   ${\bullet}$  Note:

     \medskip
\vspace{2mm}
\noindent
(i)  In applications, the submanifold  $N$ is only of class $C^{m_0}$ viewed in each admissible charts of ${\cal B}$, the conclusion and the proof 
  above   remain the same.
 
  \medskip
\vspace{2mm}
\noindent
(ii)  Consider the case that in a given stratum of ${\cal B}$, there are coordinate charts  of the form  $\cup_{\alpha\in \Lambda} W_{\epsilon }(f_{\alpha}, H)$ where $\Lambda$ is the collection of the gluing parameters near a lower stratum. Then  there is a similar statement, which   includes the case of 
 the regularized moduli space. This implies that the  moduli space is a stratified cornered manifold. A proof of this is given in [L2]. In [L2] and [L3], we will 
 show that there is a Fredholm theory for the stratified Banach manifolds appeared in GW and Floer theories so that the induced weakly smooth structure on the perturbed moduli space is in fact $C^{m_0}$-smooth rather than just stratified smooth.

  \medskip
\vspace{2mm}
\noindent

\medskip
\vspace{2mm}
\noindent
 ${\bullet}$ More  weakly smooth  functions on ${\cal B}$:

Let $M$ be a compact Riemanian  manifold and $i: M\mapsto {\bf R}^d$ be an isometric embedding.  Consider   the induced embedding 
$i_*: {\tilde {\cal B}}( M)={\tilde {\cal B}}_k^p(M)\rightarrow {\tilde {\cal B}}({\bf R}^d)= {\tilde {\cal B}}_k^p({\bf R}^d)$, where 
${\tilde {\cal B}}( M)$ and ${\tilde {\cal B}}({\bf R}^d)$ are the spaces of the  parametrized $L_k^p$-stable maps from the domain $\Sigma$ to $M$ and ${\bf R}^d$ respectively.  It is well known  that the  induced    Sobolev  metric  on ${\tilde {\cal B}}( M)$  from the embedding  is  equivalent to the intrinsic one when $m_0=k-\frac {2}{p}$ is large. Note that  ${\tilde {\cal B}}({\bf R}^d)$ is an open  set of  the Banach space $L_k^p(\Sigma, {\bf R}^d).$
For any $x$ in $\Sigma$, the evaluation map $e _x:{\tilde {\cal B}}({\bf R}^d)\rightarrow {\bf R}^d $ is linear and $M$ is closed in ${\bf R}^d$. This implies that $e_x^{-1}(M)$ is closed for each $x$, and that ${\tilde {\cal B}}( M)$ as the intersection of all such inverse images is also closed in ${\tilde {\cal B}}({\bf R}^d)$. Clearly any smooth function
$\phi$ on $L_k^p(\Sigma, {\bf R}^d) $ pulls back to a smooth function function $i^*\phi$. We are particular interested in those smooth $\phi$  with property  that the composition $\phi\circ \Phi:\Sigma\times L_k^p(\Sigma, {\bf R}^d)\rightarrow {\bf R}^1$ is still smooth. Here 
$\Phi:\Sigma\times L_k^p(\Sigma, {\bf R}^d)\rightarrow L_k^p(\Sigma, {\bf R}^d)$ is the action map. Any such $\phi$  will be called $G$-smooth.
Since  the $G$-actions  are  compatible with $i_*, $  if $\phi$ is $G$-smooth, so is $i^*\phi$.

\medskip
\vspace{2mm}
\noindent
 ${\bullet}{\bullet}$  Weakly smooth  cut-off functions on ${\cal B}$:  
 In the last section,  we give an elementary proof of  the well-known  fact that  for $p$ is a positive even integer, the $p$-th power of the $L_k^p$-norm on 
 $L_k^p(\Sigma, {\bf R}^d)$, denoted by ${\tilde \phi}_{k, p}$ is a $G$-smooth function.  Composed with a bump-off function $\beta$, we get a cut of function  $\beta\circ\phi_{k, p} $ on $L_k^p(\Sigma, {\bf R}^d)$ denoted by  ${\tilde \psi}_{k, p}.$ Let ${ \psi}_{k, p}$ be its pull-back to ${\tilde {\cal B}}$. Then for any fixed local slice $W_{\epsilon }(f, H)$,  the restriction of ${ \psi}_{k, p}$ to the slice is smooth.   Using the  explicit  form of  the transition function $T_{f,g}:W_{\epsilon }(g, H_g)\rightarrow W_{\epsilon }(f, H_f)$ defined before, it is easy to see that ${ \psi}_{k, p}\circ T_{f,g}:W_{\epsilon }(g, H_g)\rightarrow {\bf R}^1 $ is of class $C^{m_0}$.
This implies the existence of $C^{m_0}$-smooth cut-off function on ${\cal B}$. Then the standard argument in Lang's book implies that 
 ${\cal B}$ has partition of unit subordinated to any locally finite covering.

\medskip
\vspace{2mm}
\noindent
 ${\bullet}{\bullet}$  Weakly smooth   functions on ${\cal B}$ induced from $G$-linear functionals on $L_{-k'}^q(\Sigma, {\bf R}^d)$:

As usual, any ${\bf R}^d$-valued smooth function $\phi:\Sigma\rightarrow {\bf R}^1 $ induces a  continuous linear functional $T_{\phi}$ on 
$L_k^p(\Sigma, {\bf R}^d)$ by $L^2$-paring. The same elementary  argument  in last section show that $T_{\phi}$  is $G$-smooth on 
$L_k^p(\Sigma, {\bf R}^d)$. Therefore its pull-back to ${\cal B}$  defines a weakly  $C^{m_0}$-smooth function. By completing   all such $T_{\phi}$  in $L_{-k'}^q(\Sigma, {\bf R}^d)$ and pulling them back to ${\cal B}$, we get the corresponding weakly  smooth functions of certain regularity induced from $L_{-k}^q(\Sigma, {\bf R}^d).$ For instance, when $q'=-k+1$,  the regularity is of class $C^1$. This gives a large supplies of basic $C^{r}$-smooth function.
 One  can   use  these functions to prove the  $C^{m_0}$ effective smoothness of ${\cal B}$ instead of using those evaluation  maps.

\medskip
\vspace{2mm}
\noindent
 ${\bullet}$ The local trivializations  of the bundle $({\tilde {\bf L}}\rightarrow {\tilde {\cal B}})$:

Recall that in GW theory discussed here, the fiber of the bundle ${\tilde {\bf L}}_f$ is defined to be $L_{k-1}^p(\Sigma, \Lambda^{0, 1} (f^*(TM))).$

To give  ${\tilde {\bf  L}}$ a Banach bundle structure rather than just a family of Banach  spaces, we recall one of the standard ways to give it  a local trivialization by using the $(J, G_J)$-invariant connection $\nabla=\nabla^M$ of $M$ defined before. 
Given two points $x$ and $y$ in $M$ sufficient close to each other so that they can be jointed by a unique short geodesic, 
we will denote   the parallel transport of $TM$ from $x$ to $y$ along the connecting  geodesic  by $\pi_{x y}$.
Then a local $C^{\infty}$-trivialization of ${\tilde {\bf  L}}$ near a point $f$ can be given by
$\Phi: U_{\epsilon}(f)\times L^p_{k-1}(\Sigma, \wedge^1 (f^*(TM))\rightarrow  {\cal L} ,$  $\Phi (\xi, \eta)(x)=\pi_{f(x)exp_{f(x)}(\xi(x))}(\eta(x))$ for any $(\xi, \eta)\in  U_{\epsilon}(f)\times L^p_{k-1}(\Sigma, \wedge^1 (f^*(TM)).$  Note that  the parallel  transport
only acts on the part of $\eta(x)$ involving its "value" in $T_{f(x)}M$ and has no effect on the $1$-form part.

By restricting  $\Phi$  to the local slices $W(f, H)$, we get local trivializations for the local bundle   ${\bf  L}|_{W(f, H)}$.
 Note that if two of  these local slices $W(f, H), $ being considered as coordinate charts for ${\cal B}, $ have non-empty intersections, they may  not intersect
in ${\tilde {\cal B}}$ in general. They  are only related by  some  (non-constant) actions of $G$. To define a locally trivial bundle structure for the quotient bundle
 ${\cal L}\rightarrow {\cal B}, $ it is necessary  to compare the above  local trivializations
over these slices that  are only related by $G$-actions. To this end, we introduce  different  coordinate charts and local trivializations
for ${\tilde {\bf L}}\rightarrow{\tilde {\cal B}},$ which are $G$-equivariant but only continuous.


\medskip
\vspace{2mm}

\medskip
\vspace{2mm}
\noindent
${\bullet}$ $G$-equivariant local trivialization for the bundle ${\tilde {\bf  L}}\rightarrow {\tilde {\cal B}}.$

\medskip
\vspace{2mm}
Note that the local trivialization for ${\tilde {\cal L}}$ above  is given by  parallel transport that has no effect on $1$-forms.

This is not
adequate for some part of our discussion. We now introduce another system of  local chart and local trivialization which behave better with respect
to the $G$-action.


We  cover ${\tilde {\cal B}}$ by the charts $ W(f, H)\times G$, and consider it as a topological Banach manifold.
With respect to  each of these  coordinate charts $ W(f, H)\times G$, we define the trivialization  of ${\tilde {\bf L}} $  in the obvious
way  as follows.
Along $ W(f, H)$ direction, we trivialize ${\tilde {\bf L}} $ as before, but along $G$-direction, we simply use the pull-backs  induced by the $G$-actions as reparametrizations of the domain to identifying fibres along a $G$-orbit.   This gives rise  $G$-invariant coordinate
 charts  and $G$-equivariant local  trivializations for the bundle ${\tilde {\bf L}} \rightarrow {\tilde {\cal B}}. $ 
 Of course, the  $G$-equivariant transition functions between
  charts and trivializations so defined  are only continuous. 
 

 By  restricting   ${\tilde {\bf  L}}$ to  local slices $W(f, H)$, or to its $G$-orbits $G\times W(f, H), $ we get a  collection of trivial bundles.
 As we mentioned before, these local bundles together define a bundle  ${\bf  L}$ that can be thought as a topological Banach  bundle on ${\cal B}.$

Recall that the   transition functions $t_{21}$ between two local slices $W(f_1, H_1)$ and 
$W(f_1, H_1)$ have been given by using (non-constant) $G$-actions obtained from the $C^{m_0}$-map $T_{21}:W(f_1, H_1)\rightarrow G$. Then transition functions between any two  local bundles over these  local slices are induced by pull-backs from the actions. This gives  the local trivializations for the topological Banach bundle ${\bf  L}$ over ${\cal B}$.


 
 
\vspace{2mm}

\noindent
${\bullet}$ Weakly smooth sections  of the topological Banach bundle  $( { {\bf  L}} , { {\cal B} }). $ 

\vspace{2mm}
\noindent
 
 Recall that we use the "full" covering ${\cal U}$ of the collection of all local slices $W(f, H)$ to define the weakly smooth structure on ${\cal B}.$ We now use the same coving to define the ${\cal O}_{\cal B}$-module ${\cal L}$ associated to the bundle $( { {\bf  L}} , { {\cal B} }). $ 
 
  We first describe a way to construct weakly smooth section of ${\cal L}.$ 
  Consider  a smooth and $G$-equivariant section $s: {\tilde {\cal B}}\rightarrow
 {\tilde {\cal L}} .$   Then $s$ restricted to any slice $S(f, H)$ is still smooth,  and
 any two such restrictions are related by a  transition function for the bundle ${\bf   L}$. In other words, $s$ can be thought as a section of the quotient bundle  ${\bf   L}$, which    is smooth viewed in any admissible local trivialization despite of the fact that $G$-action is only
 continuous.  Recall that such a section  $s$  was called an  weakly smooth section of ${\cal  L}.$

   Here  are  some  details on above  construction in term of local charts and trivializations.
 
 Let $W_1$ and $W_2$ are two local slices of ${\tilde {\cal B}}.$ The $C^{m_0}$-coordinate transformation 
 $t_{21}:W_1\rightarrow W_2$ is given by $t_{21}(\xi)=\xi\circ T_{21}(\xi).$
 Assume that $s$ is smooth and $G$-equivariant.
 Then for any $g\in G$, $g^*(s)=(g^{-1})^*(s\circ {g})=s.$ 
  Let $s_i=s|_{W_i}, i=1,2 $.  After composed  with the transition functions  between  the two coordinate charts and trivializations,   the section $s_1 $  is transformed into     a section in $W_2$, denoted  by $s'_1$. Then for any $\eta$ in $W_2,$ we have
  $$s'_1(\eta)= (T_{21})^*(\eta)  s_1(T_{12}(\eta)\cdot(\eta))$$$$=(g^{-1})^*(s|_{W_1}\circ {g})=(g^{-1})^*(s\circ {g})|_{g^{-1}(W_1)}=
  s|_{W_1}.$$  Here $g= T_{12}$.      Therefore, $s'_1=s_2$ and hence is of class
  $C^{\infty}$. 
  
   \medskip
\vspace{2mm}
\noindent
Note that $T_{12}$ is only of class $C^{m_0}$.

We already mentioned before, in GW and Floer theory , we need to  decide if a  smooth section defined on a local slice $W_f$  is still 




 One important example of  a $G$-equivariant $C^{\infty}$-section is the 
 section ${\bar {\partial }}_{J, H}:{\tilde {\cal B} } \rightarrow {\tilde {\bf  L}} $  used to define the Floer trajectories,  or  ${\bar {\partial }}_{J}$-section used to define $J$-holomorphic maps.
 
 \noindent
${\bullet}$ The section ${\bar {\partial }}_{J}:{\tilde {\cal B} } \rightarrow {\tilde {\bf  L}}.$ 

\vspace{2mm}
\noindent
It is well-know that   ${\bar {\partial }}_{J} $ is a  $G$-equivariant $C^{\infty}$-section. 
Therefore, we get a weakly  smooth section of the bundle ${\bf L}\rightarrow{\cal B},$ denoted by $s$, which is a smooth Fredholm section viewed in any local slice. 




The question is if the assumed Fredholm section $s:{\cal B}\rightarrow {\cal L}$  has enough weakly  smooth perturbations to achieve
 transversality for the moduli space. The lack of the transversality of $s$  at a point $f$  with $f$ in  the moduli space ${\cal  M}= s^{-1}(0),$ is measured
 by the co-kernel of the  derivative  of $s$ at $f$. Here the  derivative is a bounded linear map  $Ds_{ f}: L^p_{k}(\Sigma,  f^*(TM), h)\rightarrow L^p_{k-1}(\Sigma, \wedge^{0,1}(\Sigma) \otimes f^*(TM))$. 
 
 The elliptic regularity for the  non-linear ${\bar{\partial }}_{J, H}$-operator
 or  ${\bar{\partial }}_{J}$-operator implies that $f$  is smooth and the co-kernel $K_f$ is a finite dimensional subspace of  $L^p_{k-1}(\Sigma, \wedge^{0, 1}(\Sigma) \otimes f^*(TM))$
consisting of smooth elements.  
However, as  mentioned in the introduction, our construct below does not use the smoothness of $f$.

Each element $\eta$ in $K_f$  can be thought as a section of ${\bf L}$ over  the single point $f$.  We need to extend $\eta$ into a weakly $C^{m_0}$-smooth section over
${\cal B}$. These extension  satisfy the conditions described below.

  \medskip
\vspace{2mm}
\noindent
(I) Let $GE: K_f\rightarrow GE(K_f)$ be the desired extension. We require that it is a linear map (in fact it is a linear isomorphism). Therefore, we only need to describe  the extension for finitely many elements that form a basis of $K_f$.

  \medskip
\vspace{2mm}
\noindent
(II) Let $W_f$ be a local slice and $[s_f]:W_f\rightarrow  L^p_{k-1}(\Sigma, \wedge^{0,1}(\Sigma) \otimes f^*(TM))$ be the section $s$ written in the local slice and trivialization. Assume that $GE(K_f)$ is the extension over $W_f$. Then it gives  the  evaluation map   written in the local trivialization,   $[ev]:W_f\times GE(K_f)\rightarrow L^p_{k-1}(\Sigma, \wedge^{0,1}(\Sigma) \otimes f^*(TM))$. We require that $D[s_f]_g\oplus [ev]_g:L^p_{k}(\Sigma,  f^*(TM), h)\oplus GE(K_f)\rightarrow L^p_{k-1}(\Sigma, \wedge^{0,1}(\Sigma) \otimes f^*(TM))$ is  "quantitatively" surjective over $W_f$. The exact meaning of the 
"quantitative"  surjectivity will be explained later in this section.

  \medskip
\vspace{2mm}
\noindent
(III)  We require that $GE(K_f) $ is localized near $f$. In other words, the support of each element of  $GE(K_f) $ is contained in $W_f$. This can be done by multiplying a weakly $C^{m-0}$-smooth cut-off function supported in $W_f$ and equal to one on $W'_f\subset \subset W_f$.  Still  denote the resulting extension 
by $EG(K_f). $ Then the tansversality condition in (II)  is still satisfied over
$W_f.$

Note that in (III) we have used the existence of the weakly  $C^{m_0}$-smooth cut-off function.
We will show in this section that the extension $GE(K_f) $ above is sufficient
to achieve the local transversality for the perturbed $s$-section.

As mentioned in the introduction, the desired extension can not be obtained  by using the "standard" method by  regarding the elements in $K_f$ as "constant" sections over $W_f$, then moving them over the orbit of $W_f$ by pull-backs.



To get the desired extension,   fix  a $C^{\infty}$ partition of unit $1=\beta_i, i\in I,$ subordinated to a finite by sufficiently fine (depending on f) covering of ${\Sigma}=\cup_{i\in I}D_{\delta_i}(x_i).$ Then the  cokernel $K_f$  is contained  the  finite sum  $\Sigma_{i\in I} K_i$ as   finite dimensional spaces
 of  $L^p_{k-1}(\Sigma, \wedge^{0,1}(\Sigma) \otimes f^*(TM)).$  Here  $K_i=\beta_i\cdot K_f$  and $D_{\delta_i}(x_i)$ is  the disc of radius $\delta_i$ on ${\Sigma}$ centered at $x_i.$ 
 We will assume that the covering is fine enough and $W_f$ is small enough such that for any $i\in I $ ,  $g\in W_f $  and $x\in D_{\delta_i}$, $|g(x)-f(x_i)|<\epsilon_i$. Since $\|f-g\|_{C^0}$ is bounded by the corresponding $L_k^p$-norm. The assumption can be satisfied. Here  $\epsilon_i$ is chosen such that  there is a
 fixed trivialization of $TM$ over the ball of radius $2{\epsilon_i}$ centered
 at $f(x_i)$. Note that the image  $g(D_{\delta_i}(x_i))$ is in this ball for all
 $g$ in $W_f.$  Clearly, we only need to  construct the extension for the elements in  a basis  of  $K_i, i\in I.$

 Now we change the notation. Simply  use $K_f$  to denote  one of these $K_i.$
 
  
 Clearly,  it is sufficient to show that  for each ${\eta}$ in the basis of $K_f$,  there exists an approximated extension $GE(\eta)$ over $W_f$ such that the $L^2$-norm of $(GE(\eta)(f)-\eta)$ is less than a prescribed positive number. 
  This is what we are going to do next. After this we will show that it is possible to give a true extension rather than just a approximate one under the assumption that all the initial data are $C^{\infty}$-smooth.
  
  Let $\eta$ be  a element in the basis of $K_f$. Then it is a finite linear combinations  of the elements of the form  $\phi\otimes \xi$, where $\phi$ is a smooth  $(0,  1)$-form supported  in  $D_{\delta_i}(x_i)$ and $\xi$ is a  local $L_k^p$-section of the bundle $(f^*(TM)\rightarrow \Sigma).$



 We now show that $GE(\eta)$ can be obtained   by combining two kinds of $G$-equivariant  sections. Both are  "geometric" nature.

\medskip
\vspace{2mm}
\noindent
 ${\bullet}$ The "geometric" perturbations   of the  ${\bar {\partial }}_{J}$-section:

Note that the fibre of ${\tilde{\cal L}}_{k,p}$ at $f$,  $L^p_{k}(\Sigma, \wedge^{0,1}(\Sigma) \otimes f^*(TM)),$ is linearly
homeomorphic to  $L^p_{k}(\Sigma,  f^*(TM)) \otimes_{L^p_{k}(\Sigma)}L_k^p(\wedge^{0,1}(\Sigma)) $ when $k$ is large enough.
 This gives rise a bundle isomorphism ${\tilde{\cal L}}_{k,p}\simeq {\tilde{\cal T}}_{k,p}\otimes {\tilde \Omega}^{1}.$
 Here   the fibre of ${\tilde{\cal T}}_{k,p}$ at $f$ is $ L^p_{k}(\Sigma,  f^*(TM)), $  and ${\tilde \Omega}^{1}$ is the  trivial bundle whose fiber is $L_k^p(\wedge^{0,1}(\Sigma)) .$
 

We will show that  certain  smooth section of ${\tilde \Omega}^{1}$ on a local slice  is in  fact $C^{m_0}$-smooth viewed in other slices after composed with the transition functions between the slices.  
To this end,  we restrict  the $G$-bundle ${\tilde \Omega}^1$  to  the  orbit  $G\times W_1(f, H) $ of a local slice $W_1(f, H) $.


  Now let $\phi$ be a $C^{\infty}$ $(0,1)$-form on $\Sigma$ considered  as a element of the fiber of ${\tilde \Omega^1}.$ 
  Since the bundle is trivial, we get a constant section $\psi_1$ on the slice $W_1(f)$ defined by $\psi_1 (\xi)=\phi.$ It extends to a $G$-equivariant 
  section over $G\times W_1(f)$ by pull-backs. Denote the extended section by ${\tilde \psi_1}$.  Even assume that $\phi$ is  a $C^{\infty}$ $(0,1)$-form, it is still not immediately
  clear that this extended $G$-section is smooth on the open set $G\times W_1(f). $  On each $G$-slice, the extended section
  ${\tilde \psi_1}:G\times \{h_0\} \rightarrow L_k^p(\wedge^{0,1}(\Sigma))$ has the form ${\tilde \psi_1}(g, h_0)=g^*(\phi),$  which is smooth as a function  on $G.$ Presumably, this should imply that  the extended section is smooth  as usually the lack of differentiability comes from $G$-direction. 
  However,  as  mentioned in the introduction,  the best one can get is its $C^{m_0}$-smoothness.
  
  Instead of repeating  a similar argument, we directly prove the weaker statement that the restriction of the extended section
  to another  slice $W_2(f)$ is of class $C^{m_0}$.

  Recall that the coordinate transformation between the two slice $t_{1, 2}: W_2(f)\rightarrow W_1(f\circ g_0)$ is defined to be $t_{1, 2}(\xi)=\xi\circ T_{12}(\xi).$ Here $T_{12}:W_2(f)\rightarrow G$ is of class $C^{m_0}$ and $G$ is  the group of reparametriaztions of ${\Sigma}.$
  Then $\psi_2(\xi)=\psi_1\circ t_{1, 2} (\xi)={\tilde \psi_1}(\xi \circ T_{1, 2} (\xi))=(T_{1, 2} (\xi))^*({\tilde \psi_1})(\xi)=(T_{1, 2} (\xi))^*(\phi).$
  Here since ${\tilde \Omega^1} $   is a trivial bundle, the sections $\psi_i, i=1,2 $ are considered to be maps $W_i:\rightarrow  L_k^p(\wedge^{0,1}(\Sigma))$.

  In other words, $\psi_2=O_{\phi}\circ T_{1, 2}$ where the orbit map $O_{\phi}:G\rightarrow L_k^p(\wedge^{0,1}(\Sigma))$ is  defined by $O_{\phi}(g)=g^*(\phi).$ Note that since $\phi$ is $C^{\infty}$, $O_{\phi}$ is of class $C^{\infty}.$ Therefore, $\psi_2$ is $C^{m_0}$-smooth. Note that  there are two more general cases that the  argument above implies the same conclusion:  (i)  the constant map  $\psi_1:W_1(f) \rightarrow L_k^p(\wedge^{0,1}(\Sigma))$ can be replaced by any smooth map with image lying in $C^{\infty}(\wedge^{0,1}(\Sigma)); $   (ii)  $\phi$  is only of class $C^{2m_0}$. 
  


 What we did above gives a way to extend the part of the section of  $\phi\otimes \xi$ , $\phi$, to a smooth section $\psi_1$ over $W_1(f)$ which is of class  $C^{m_0}$ viewed in any other slices.

 As for the $\xi$, we do not extend it, but approximate it first under the assumption that  the covering is very fine so that we can approximate ${\xi}$ by the constant section  $\xi_0$ on $D_{\delta_i}(x_i)$ defined by $\xi_0(x)=\xi(x_i)$  for  $x\in D_{\delta_i}(x_i)$. Note that  this local  constant section can be approximated 
 by the pull-back $f|_{D(x_i)}^*({\tilde \xi}_0)$. Here ${\tilde \xi}_0$ is the local constant section of the bundle  $
 TM\rightarrow M$  in a neighbourhood of $f(x_i)$ given by transporting the value $\xi_0(x_i)$ in $T_{f(x_i)}M$ to the neighbourhood using $Dexp_{(f(x_i))}.$

  Clearly  ${\tilde \xi}_0$ is smooth  and we  may assume that it is  the restriction of a global smooth section, still denoted by ${\tilde \xi}_0$ of the bundle $TM\rightarrow M.$ Therefore, the section $\xi$  can be approximated by the pull-back of a global section ${\tilde \xi}_0$ of $TM$ by  $f|_{D(x_i)}$ .


 
  Now  $\xi$  can be  approximately extended  to  a  global section ${\tilde \xi}$ over ${\tilde {\cal B}}$ defined by ${\tilde \xi}(g)=g^*({\tilde \xi}_0). $   Since ${\tilde \xi}$ comes from the pull-backs of the smooth  "geometric" section of $TM\rightarrow M$,  it is  automatically  
  $G$-equivariant and  smooth. 
  
  More specifically, since the $G$-actions only act on the domain,   the $G$-equivariancy of ${\tilde \xi}$ is clear.  On the other hand, using the embedding $i:M\rightarrow {\bf R}^d,$  the proof of  its smoothness, which is supposed to be well-known,  can be reduced further to the case that $M$ is just ${\bf R}^d$.  In next section, we include the elementary proof of smoothness  for  this latter case for  completeness.

 At this point, we need to make sure that  after  putting  all these extended section together,  we get the extended section $GE(\eta)$ of $\eta$ such that $\|GE(\eta)(f)-\eta\|_{0, 2}<\delta$ for a prescribed $\delta>0.$
  To  this end, we write everything in local charts of $\Sigma$ given by the covering
  ${\Sigma}=\cup_{i\in I}D_{\delta_i}(x_i) $ and fixed  local trivializations of $TM$ near
  $f(x_i), i\in I.$

Then ${\eta}=\Sigma_{i\in I}\beta_i \eta.$  
Let $\phi_i$ be the $C^{\infty}$ local frame of the bundle $\Lambda^{0, 1}\rightarrow \Sigma$ over $D_{\delta_i}(x_i) $ with $\|\phi_i\|_{C^0}=1$. 
Denote ${\eta}|_{D_{\delta_i}(x_i)}$ by ${\eta}_i$.  Then ${\eta}_i=\phi_i\cdot \xi_i$ for some $\xi_i$ in $L_k^p(\Sigma, f^*(TM), h)$. In fact, we only need each $\xi_i$  to be defined over $D_{\delta_i}(x_i)$.
Using the local trivializations of $TM$ near $T_{f(x_i)}, $ we  define the corresponding local constant section $\xi_i^{x_i}$ over $D_{\delta_i}(x_i)$ given by $\xi_i^{x_i}(x)=\xi_i(x_i)$.

Then $$\eta= \Sigma_{i\in I}(\beta_i\eta)=\Sigma_{i\in I} (\beta_i\eta_i)$$
$$=\Sigma_{i\in I} (\beta_i\phi_i\cdot \xi_i). $$
We define a $K_k^p$-section $\eta'$ which is $C^0$-close to $\eta$ as follows,
$$\eta'=\Sigma_{i\in I} (\beta_i\phi_i\cdot \xi^{x_i}_i). $$
Note that $\eta'$ behaves like a $C^{\infty}$ section:$\beta_i\phi_i$ is smooth on 
$\sigma$ and $\xi^{x_i}_i$ is obtained by pulling back a smooth  vector field on $M$  by  $f|_{D_{\delta_i}(x_i)}$.

Then there constants $C$ and $C'$ such that when the covering is fine enough, for any prescribed $\epsilon'$, 
$$\|\eta-\eta'\|_{0, 2}\leq C \|\eta-\eta'\|_{C^0}
$$$$ \leq C Max_i\{\|\beta_i\phi_i\cdot (\xi_i-\xi_i^{x_i})\|_{C^0}\}$$
$$\leq C' Max_i\{ \|(\xi_i-\xi_i(x_i)|_{D_{\delta_i}(x_i)}\|_{C^0}\}\leq \epsilon'.$$  

Using  this $\eta'$ to replace $\eta$, the argument before gives the desired approximate extension $GE(\eta)$.  

In above argument, we only assume that $f$ and $\eta$ is of class $L_k^p$ and $L_{k-1}^p$ respectively.  If we assume that they are of class $C^{\infty}$ or of class $C^{2m_0},$ we get better result. We will only give the proof   for  $C^{\infty}$ case. 

 To this end, we change notation and  denote $\beta_i \eta$ by $ \eta_i$.
Then each $ \eta_i=\Sigma_{j\in J} \phi_{j}^i\cdot e^i_j(f_i)$. Here $f_i=f|_{D_{\delta_i}(x_i)} ,$ $\phi_{j}^i$ is a $C^{\infty}$-smooth $(0, 1)$-form on ${\Sigma}$ supported on ${D_{\delta_i}(x_i)}$ and  the collection of $e^i_j, j\in J$ is a local  $C^{\infty}$-frame of $TM$ near $T_{f(x_i)}M.$ From this local expression, one immediately see that each $e^i_j(f_i)=f_i^*(e^i_j), $ and hence it extends to a smooth 
section ${\tilde e^i_j}$ on  $W_f$ defined by ${\tilde e^i_j}(g)=g_i^*(e^i_j), $ where $g_i=g|_{D_{\delta_i}(x_i)} .$  We can extend each $\phi_{j}^i$ into ${\tilde \phi_{j}^i}$ as before.  Let 
$GE(\eta)=\Sigma_{i\in I} GE(\eta_i)=\Sigma_{i\in I, j\in J}{\tilde \phi_{j}^i}\cdot {\tilde  e^i_j}.$ Then $GE(\eta)$ so defined is  a true extension of $\eta$ rather than just   a  approximate extension.



In this way, we  extend    each element $\eta$ in the basis  of $K_f$ to a smooth section over  $W_1(f)$,  which is $C^{m_0}$-smooth viewed in any other
slices.
The results so far in this section are sufficient  to resolve the  difficulty of lacking of differentiability and to establish local transversality 
for the case that there is only
one stratum.

In the rest of this section, we will finish the construction of perturbed moduli space and prove that it is a smooth compact manifold  with the expected  dimension.

The argument below for regularizing a compact moduli space ${\cal M}$  works   quite generally  under the assumptions that  ${\cal M}$ has only one stratum, all isotropy groups for elements in ${\tilde {\cal B}}$ are trivial and the Fredholm section $s: {\cal B}\rightarrow {\bf L}$ is proper. 

\medskip
\vspace{2mm}
\noindent
 ${\bullet}$ Regularization of the moduli space  ${\cal M}$:  
 
 Let ${\cal M}=s^{-1}(0)$ be the moduli space of  unparametrized stable  ${J}$-holomorphic   spheres.  By our assumption ${\cal M}$ is compact.  For each $f$ in ${\cal M}$, let $W_f=W_{\epsilon}(f, H_f)$ be a local slice regarded as a local chart of ${\cal B}$ containing $f$. In the local chart $W_f $ and the trivialization of ${\bf L}$ over $W_f, $ we denote 
 the Fredholm section $s$ by $[s]=[s_f].$ The derivative of $s$ at point f  written  the local trivialization over $W_f, $ 
 $D[s_f]_f:
 L_k^p(\Sigma, f^*(TM), h)\rightarrow L_{k-1}^p(\Sigma, \Lambda^{0, 1}(f^*(TM))), $ is given by the well-known formula,
 $D[s_f]_f(\xi)=\nabla \xi+ J(f)\nabla \xi\circ i+ N({\partial} f, \xi). $ Here $s$ is the ${\bar {\partial }}_J$-section,   $\nabla$ is the unique $J$-invariant connection preserving
 the $g_J$-metric   whose torsion is equal to the torsion of $J$, $N=N_J$ and $g_J$ is  defined by $g_J=\omega (\cdot, J \cdot). $   There is a similar formula
 for ${\bar {\partial }}_{J, H}$-section. Although, without introducing a connection it only makes sense to take derivatives of a section  invariantly at its zeros, the above formula is still applicable  at  a general point $g$ in $W_f$ if we use  the local trivialization above. In particular, from these explicit formulas, we have that when $\epsilon$ is  small enough, for any $g$ in $W_{f, \epsilon}, $ $D[s_f]_g$ and  $D[s_f]_f$  are  close to each other with respect to the operator norm.
 
 ${\bullet}$ ${\bullet}$ Choice of  $\epsilon= \epsilon_f$ for $W_{f, \epsilon}:$
 
 We need to make choice of $\epsilon= \epsilon_f$ for $W_{f, \epsilon}$ such that the perturbed sections  of ${\bar {\partial }}_J$ by certain collection  of  sections related to the cokernels  achieve the "controlled"  transversality.

 Let  $K_f$ be the cokernel of $D_f=D[s_f]_f$,  $C_f$ be  its  kernel. We denote the $L^2$ orthogonal complement of $C_f$ in $L_k^p(\Sigma, f^*(TM), h)$ by  $N_f$,  and the $L^2$ orthogonal complement of $K_f$ in $L_{k-1}^p(\Sigma, \Lambda^{0,1} (f^*(TM)) )$ by ${\bar N}_f$. Then $D[s_f]_f:N_f\rightarrow {\bar N}_f$ is an isomorphism between the two Banach spaces.
 
 When $\epsilon=\epsilon_f$ and $\delta=\delta_f $  are  small enough, for any $g$ in $W_{f, \epsilon}$ and   bounded linear operator 
 $T:  L_k^p(\Sigma, f^*(TM), h)\rightarrow L_{k-1}^p(\Sigma, \Lambda^{0, 1}(f^*(TM))) $ with operator norm less than $\delta$, 
 $D[s_f]_g+T:N_f\rightarrow {\bar N}_f$ is still an isomorphism. This $T$ corresponds to the derivatives of perturbation sections along $W_f$ directions. Since  the perturbations will be made  small enough so that all relevant derivatives along $W_f$ directions are ignorable, we let $T$ to be zero first. 
 
 By definition, we have  $$D[s_f]_f\oplus I_{K_f}:
 L_k^p(\Sigma, f^*(TM), h)\oplus K_f\rightarrow $$ $$L_{k-1}^p(\Sigma, \Lambda^{0, 1}(f^*(TM)))$$ is surjective with the same kernel $C_f$.
 The right inverse of this map $G_{f, K_f}: L_{k-1}^p(\Sigma, \Lambda^{0, 1}(f^*(TM)))\rightarrow N_f\subset L_{k}^p(\Sigma, f^*(TM))$
 is a linear bounded operator with operator norm $|| G_{f, K_f}||$.
  Here $I_{K_f}:K_f\rightarrow L_{k-1}^p(\Sigma, \Lambda^{0, 1}(f^*(TM)))$ is the inclusion map.
 
  Fix a basis $(\eta_1, \cdots, \eta_{n_f})$ of $K_f$ as  finite dimensional subspace of $$L_{k-1}^p(\Sigma, \Lambda^{0, 1}(f^*(TM)))$$.
  Let $G(K_f)$ be linear space spanned by approximate sections $(G(\eta_1), \cdots, G(\eta_{n_f}))$ of  $L_{k-1}^p(\Sigma, \Lambda^{0, 1}(f^*(TM)))$. Assume that $||\eta_i-G(\eta_i)||_{k-1, p}\leq \delta|| $ for all $i$.  When $\delta$ and $\epsilon$ are small enough,
  for any $g$ in $W_{\epsilon, f}$ and $ ||  G(\eta_i)- \eta_i||_{k-1, p}< \delta $  $i=1, \cdots\ n_f, $  the linear map $D[s_f]_g\oplus I_{G(K_f)}:
 L_k^p(\Sigma, f^*(TM), h)\oplus G(K_f)\rightarrow L_{k-1}^p(\Sigma, \Lambda^{0, 1}(f^*(TM)))$ is still surjective. Moreover,  
   its right inverse, denoted by  $G_{g, G(K_f)}$, has  the operator norm which is almost the same as the fixed one, $|| G_{f, K_f}||$.
  
   Now we fix such  $\epsilon=\epsilon_f$   and $\delta=\delta_f$ temporarily, and 
  assume that for each point $f$  in ${\cal M}$ such a $W_f$ is already chosen.
 
 ${\bullet}$ ${\bullet}$ The space $GE(K_f)$ of geometric perturbations derived from $K_f$:
 
 First fix a cut-off function ${\gamma}_f$ supported in $W_f$, that is  equal to one on a smaller neighbourhood  $W'_f\subset W_f.$
 
 Recall that  in  the usual construction of abstract perturbation, one considers  each element of $\eta$ of $K_f$ as a section of ${\bf L}\rightarrow {\cal B}$ at the point $f$ and extend it over the local slice $W_f$ by parallel transport to get a constant section, denoted by $ E(\eta)$.
  Let $LE(\eta)=\gamma \cdot  E(\eta)$ be the corresponding localized section. Then in the local chart and trivialization over $W'_f$, $LE(\eta)(\xi)= LE(\eta)(0)=\eta$ so that  $LE(\eta)$ is still a constant section. Let $LE(K_f)$ to be the collection of all such $LE(\eta)$. Note that Since  the two operations   used in  extending $K_f $ to  $LE(K_f)$  are  linear,  $LE(K_f)$ is a finite dimensional vector space inside the space of smooth local sections of ${\bf L}$ over $W_f,  $  which  is isomorphic to $K_f$ and hence has the same dimension as that of $K_f.$  
  Because of linearity,  $K_f$ is obtained  from the  extensions of   $\eta_i, i=1, \cdots, n_f$ in a basis of $K_f$. 
   Note that since  on $W'_f,$ any section $LE(\eta)$ in  is a constant section   in the local trivialization, 
   the  collection of all evaluations of the sections  in  $LE(K_f)$  at any point $g$ in $W'(f)$ is just  the cokernel $K_f$, which is independent of $g$. 
  
   We have already proved that for the fixed basis $(\eta_1\cdots \eta_{n_f}),$
 each element $LE(\eta_i)$ of $LE(K_f)$ can be approximated by  a localized geometric section over $W_f$, denote by  $GE(\eta_i)$.
 Recall that the approximation is obtained by approximating the corresponding constant section first when multiplying the resulting section
 by the cut-off function. 
 
 Let  $GE(K_f)$ be  the linear space spanned by the geometric section  $GE(\eta_i), i=1, \cdots,n_f.$ 
 The main reason to switch to $GE(K_f)$ is that its elements  are   not only  smooth over $W_f$ but also $C^{m_0}$-smooth viewed in  other local slices. This may not be true for the elements in $LE(K_f)$. 

  Now we bring  the fixed positive constant $\delta_f$ above into the discussion.  Once $\delta=\delta_f$ is fixed,   for sufficiently small  ${\epsilon}_f, $ we may assume that  each  approximated section  $GE(\eta_i)$ satisfies that for any  $g$ in $W_f$, $ ||  GE(\eta_i)(g)- LE(\eta_i)(g)||_{k-1, p}< \delta_f $  $i=1, \cdot\ n_f .$  Here we have considered $GE(\eta_i)(g)$ and $    LE(\eta_i)(g)$  as  elements in $L_{k-1}^p(\Sigma, \Lambda^{0, 1}(f^*(TM)))$ by using the trivialization.  Note that since on $W'_f$, $LE(\eta_i)(g)=\eta_i$ in the local trivialization.
 the choice of $\delta$  implies that for any $g$ in $W_f', $ $(GE(\eta_i)(g), \cdots,(GE(\eta_{n_f})(g)) $  is a basis of the corresponding linear subspace in $L_{k-1}^p(\Sigma, \Lambda^{0, 1}(f^*(TM)))$,which is very  close to $K_f$.
Consequently, 
with the choices of $\epsilon$ and $GE(\eta_i)$,     for any point $g$ in $W'_f$,   the linear map 
$D[s_f]_g\oplus [ev]_g:L_k^p(\Sigma, f^*(TM), h)\oplus GE(K_f)\rightarrow L_{k-1}^p(\Sigma, \Lambda^{0, 1}(f^*(TM))) $  is surjective.  Moreover, its right inverse, denoted by $G_{g, GE(K_f)(g)}, $ has almost the same operator norm as the fixed one $||G_{f, G(K_f)}||.$
Here $[ev]_g: GE(K_f)\rightarrow L_{k-1}^p(\Sigma, \Lambda^{0, 1}(f^*(TM))) $ is the obvious evaluation map at $g$ written in the local trivialization of ${\bf L}$ over $W_f$.

 ${\bullet}$ ${\bullet}$ Size of  local perturbation space:
 
 Our next task is to  decide the "size" of each local perturbation space so that a quantitative version  of  implicit function theorem is applicable to the perturbed ${\bar{\partial }}_J$-operators.  To this end, we need (i) to find the derivatives of the perturbed ${\bar{\partial }}_J$-section over each local slice , and  (ii) to understand how these derivatives  are transformed between  different slices.


   Consider the local perturbed section ${\bf ps}=s\oplus ev:   W_f\times GE(K_f)\rightarrow {\bf L}|_{W_f},  $ defined by $s\oplus ev(g, {\lambda})=s(g)+\lambda(g).$
    Let $[s_f]\oplus [ev]$ be the corresponding map written in the  local trivialization.  Then  its   derivative  at $(g, {\lambda})$, $D([s_f]\oplus [ev])_{g, { \lambda}}=(D[s_f])_g+( D[ev])_{g, { \lambda}}.$  The partial derivative  of $[ev]$ at $(g, {\lambda})$ along 
 $GE(K_f)$-directions  is just the linear inclusion map  $[ev]_g:GE(K_f)(g)\rightarrow L_{k-1}^p(\Sigma, \Lambda^{0, 1}(f^*(TM))) $ introduced before, which is independent  of  ${\lambda}. $ We already know  that $(D[s_f])_g\oplus [ev]_g$ is surjective with right
  inverse whose norm is bounded by $||G_{f, G(K_f)}||.$

  
  Now denote the partial derivatives of $[ev]$ along $W_f$-directions by ${\partial}^W [ev]$. Then for any fixed $g$ in $W_f$, ${\lambda }$ and  $c{\lambda }$ in $GE(K_f)$ for some positive constant $c$,
 we have  ${\partial}^W [ev]_{(g, c{\lambda})}(\xi)=c{\partial}^W [ev]_{g, {\lambda}}(\xi)$  for any $\xi$ in the tangent  space of $W_f$ at $g.$
   Consequently the operator norm  $||{\partial}^W [ev]_{(g, c{\lambda})}||=c||{\partial}^W [ev]_{g, { \lambda}}||.$
 In other words, the operator norm of  ${\partial}^W [ev]_{(g, {\lambda})}$ get rescaled by $c$ if  the sized of $\lambda$ is rescaled by $c$.  This is the crucial fact that we need to get local transversality for the section $s\oplus ev$. However, in order to get desired global perturbation of $s$, we need to apply Picard method to get  more quantitative information on the  local extended moduli spaces.
 To this end, we need to compute higher derivatives. Clearly, higher partial derivatives of $[ev]$ along $W_f$-directions also rescale
  in the same way as the first derivative does. This is better than what we need.   To find second partial derivatives of $[ev]$ along the other
 directions, note that $[ev](g, \cdot)$ is already linear along $GE(K_f)$-directions, therefore we only need to find the mixed partial derivatives
  of $[ev]$.
  Since $$(\frac {{\partial}^{2}}{{{\partial}W}{{\partial} GE(K_f))}} [ev])_{(g, {\lambda})}(\xi, {\tilde \eta})=({\partial} ^W {\tilde \eta})|_g(\xi).$$

  We have  $$||(\frac {{\partial}^{2}}{{{\partial}W}{{\partial} GE(K_f))}} [ev])_{(g, {\lambda})}(\xi, {\tilde \eta})||_{k-1, p}$$$$\leq  (max_{g\in W_f}(\Sigma _{i\in I}|| ({\partial} ^W (GE (\eta_i))_g|| )\cdot ||{\tilde \eta}(g)||_{k-1, p}\cdot ||  \xi||_{k, p}.$$
  Here $|| ({\partial} ^W (GE (\eta_i))_g||$ is  the operator norm of the partial derivative at $g$ and each  $GE (\eta_i)$ is the approximated extension of the element $\eta_i, i\in I$ in the fixed basis of $K_f$.

  It follows from this that

  (I) when $\rho=\rho_f$ is small enough, for all ${ \lambda}$ in the small $\rho$-ball $Bl(K_f, \rho)$ of $GE(K_f) $ and $g$ in $W_f$, the operator norm of   $||{\partial}^W [ev]_{g, {\lambda}}||$  is less than $\delta_f$ specified before. Consequently,  in the local trivialization, the derivative the section  ${\bf ps}=s\oplus ev$, $(D[s_f]_g\oplus (D[ev])_{g, \lambda}:   W'_f\times Bl(K_f,\rho)\rightarrow {\bf L}|_{W'_f} $  is a surjective  map at any point. This  solves our problem to achieve the local transversality by using perturbation form  $Bl(K_f,\rho)$ only.
  In other words, the solution set of the local equation $[s_f]\oplus [ev](X)=0$  is a smooth submanifold in $W'_f\times Bl(K_f,\rho).$
  
  (II) The right inverse of the derivative $(D[s_f]_g\oplus (D[ev])_{g, \lambda},$ denote by  $G_{g, \lambda}, $ still has  almost the same operator norm as the fixed one $||G_{f, G(K_f)}||.$ 
  
  (III) Let $N$ be the non-linear term appeared in the Taylor expansion of $[s_f]_g\oplus [ev]$ at $(f, 0)$, then it satisfies the condition
  on $N$ required by the following Lemma.

\begin{lemma}{Picard method}

Assume that a smooth map $F:E\rightarrow L$ from Banach spaces $(E,\|\cdot\|)$
to $L$ has a Taylor expansion
$$F(\xi)=F(0)+DF(0)\xi +N(\xi)$$
such that $DF(0)$ has a finite dimensional kernel and a right inverse
$G$ satisfying
$$\|GN(\xi)-GN(\eta)\|\leq C(\|\xi\|+\|\eta\|)\|\xi-\eta\|$$
for some constant $C$. Let $\delta_1=\frac{1}{8C}.$ If $\|G\circ F(0)\|
\leq\frac{\delta_1}{2}$, then the zero set of $F$ in $Bl_{\delta_1}=
\{\xi,\,|\,\|\xi\|<\delta_1\}$ is a smooth manifold of dimension equal
to the dimension of $ker DF(0).$ In fact, if
$${\bf K}_{\delta_1}=\{\xi\,|\xi\in ker DF(0),\,\|\xi\|<\delta_1\}$$
and ${\bf K}^{\perp}=G(L),$ then there exists a smooth function
$$\phi:{\bf K}_{\delta_1}\rightarrow {\bf K}^{\perp}$$
such that $F(\xi+\phi(\xi))=0$ and all zeros of $F$ in $Bl_{\delta_1}$
are of the form $\xi+\phi(\xi).$
\end{lemma}

The proof of this Lemma is an elementary application of Banach's fixed point
theorem.

 Now applying this lemma to our case with the obvious interpretations of the notations,  we  conclude that the perturbed section $s\oplus ev$ defined on $W'_f\times Bl(K_f, \rho_f)$ is transversal to the zero section so that the zero locus, the local extended  moduli space $(s\oplus ev)^{-1}(0)$ in  $W'_f\times Bl(K_f, \rho_f)$, is a smooth manifold. Denote this extended moduli space by ${\bf E}{\cal M}^{f, \rho_f}$. The key point is that  it is  realized as a graph over  a disc of radius $\delta_{1, f}$ in the kernel of $D[s_f]_f$ in the local trivialization.
 Therefore by shrinking ${\bf E}{\cal M}^{f, \rho_f}$ a little bit corresponding to taking  $\delta_{1, f}$ to be a smaller  $\delta_{2, f}$, we get a corresponding  spaces  ${\bf E}{\cal M}_2^{f, \rho_f}$. Denote the original larger space as ${\bf E}{\cal M}_1^{f, \rho_f}.$
 Then closure of ${\bf E}{\cal M}_2^{f, \rho_f}$  is compact.  
 We may assume that ${\bf E}{\cal M}_2^{f, \rho_f}$ is the corresponding
 solution space of the equation $[s]\oplus [ev](X)=0$ in the smaller space $W''_f\times Bl(K_f, \rho_f)$ for some $W''_f\subset \subset W'_f.$

 Note that each element $\lambda_f$ of the perturbation space $GE(K_f)$ has  support
 inside  a local slice $W_f$,  hence can be considered  as a global section of ${\bf L}\rightarrow {\cal B}.$ In particular, $\lambda_f$ can be viewed as a local section over another local slice $W_{f'}$. We will cover ${\cal M}$ by finitely many such local slices $W_{f_i}, i\in I$ and consider the corresponding
 $Bl(K_{f_i}, \rho_{f_i})$ and the resulting global perturbation space
 $\oplus_{i\in I} Bl(K_{f_i}, \rho_{f_i}).$ Our goal is to show that for proper choice of the covering and small enough generic $\nu$ in $\oplus_{i\in I} Bl(K_{f_i}, \rho_{f_i}), $ the perturbed moduli space ${\cal M}^{\nu}$ is a compact  $C^{m_0}$-manifold.
 
 To this end, we need to  modify the above discussion to  incorporate the effect of the perturbations from $Bl(K_{f_j}, \rho_{f_j})$ with $j\not = i$ on the slice
 $W_f=W_{f_i}.$
 
  For that purpose, we   give a more general discussion first.

 Consider a ball of radius $r$ centered at origin  in a finite dimensional vector space, denoted by $Bl( r), $  which is linearly mapping into
 the space of "global" weakly smooth sections of the bundle ${\bf L}\rightarrow {\cal  B}. $ Here a "global" section of the bundle is  simply
  a collection of  compatible local sections with respect a covering data that we will specified in a moment.  This $Bl( r) $ will play the role of the space of global perturbations. Assume that the $W_f$ is 
  one of the local slice, and  $U_f\subset\subset W''_f$ is a small open neighbourhood of $f$.  Denote  corresponding sections of  $Bl( r) $ over $U_f$ by $Bl(f,  r)$.

  Consider the solution space of the equation about $(g, \eta) $ in $U_f\times Bl(f,  r)$ given by $[s](g)+\eta (g)=0. $
   Clearly the  above  solution space is covered by the solution space  about $(g, \xi, \eta) $ with $(g, \xi)$ in $U_f\times Bl(K_f, \rho_f)$ and $\eta$ in $ Bl(f,  r)$ given by $[s](g)+\xi(g) =-\eta (g). $   By the implicit function theorem above and the way we define $[s]\oplus [ev]$, there is a sufficient  small positive $r_f$ such that  when $r<r_f$, for any fixed $\eta $ in $ Bl(f,  r_f)$, all solutions of the above equation about $(g, \xi)$ in  $U_f\times Bl(K_f, \rho_f)$  is  homoemorphic to a finite dimensional
  disc, still denoted  by ${\bf K}_{\delta_2} ,$ the same notation used in the Picard method. Therefore, the closure of the solution space is contained in a compact set in $W''_f\times cl(Bl(K_f, \rho_f))\times cl(Bl(f,  r_f))$  which is homeomophic to  $cl({\bf K}_{\delta_2}) \times  cl(Bl(f,  r_f)).$ 
  
    To summary what we have done here,  consider the  solution space of the equation about $(g, \eta) $ in $U_f\times Bl(f,  r_f)$ given by $[s](g)+\eta (g)=0 $, and denote it by ${\bf E}{\cal M}(Bl(f,  r_f))$. 
    Then (1) $cl({\bf E}{\cal M}(Bl(f,  r_f)))\subset \pi_{W_f}(cl({\bf E}{\cal M}_2^{f, \rho_f}))\times cl(Bl(f,  r_f))$ which  is compact in $W''_f\times cl(Bl(f,  r_f))$.  Here $\pi_{W_f}:W_f\times Bl(K_f, \rho_f)\rightarrow W_f$ is the projection to $W_f$. (2) For any fixed $\nu_f$ in $Bl(f,  r_f), $  Let ${\cal M}^{\nu_f}$ be the solution space of the equation $[s](g)+\nu(g)=0 . $ Then $cl( {\cal M}^{\nu_f})$ is compact in $W''_f$. Since $W''_f$ is a local slice for ${\cal B}, $ we  conclude that the closure in ${\cal B}$,  $cl_{\cal B}( {\cal M}^{\nu_f})$  is compact.

  Since ${\cal M}$ is compact, we can select a finite covering ${\cal W}=\{W_i=W_{f_i}, i\in  N \}$ such that ${\cal M}$ is already covered by   the corresponding  open subsets $U_i\subset\subset W''_i, i\in  N.$ 
  
 

 To finish the construction of the perturbed moduli space, we need to make two assumptions:

  \medskip
\vspace{2mm}
\noindent
$A1:$ Let $U=\cup_{i\in I}U_i$ be the  open subset of ${\cal B}$. Then the boundary of $U$, denoted by $Bd(U)$, has no intersection with ${\cal M}$.
 Here $Bd(U)=cl_{\cal B}(U)\setminus U.$

  \medskip
\vspace{2mm}
\noindent
 $A2:$  There is a positive constant $C_0$, such that for any $g$ on $Bd(U), $ $||s(g)||>C_0.$
 
 We will prove that these two assumptions can be achieved latter in this section.
 
   Let  $B_i=BL(K_{f_i}, r_{f_i})$ be the  ball of radius $r_i=\rho_{f_i} $ in the linear space $GE(K_{f_i})$.
   
  Then each element of $GE(K_{f_i})$ as a section of ${\bf L}$ over $W_i$ is supported in a closed subset of  $W_i$ and hence can be regarded as a 
   global  section of the bundle ${\bf L}\rightarrow {\cal B}.$

   
Set $r=\Sigma _{i\in I} r_i, $ and    let $B(r)=\oplus_i  B_i $. 
   
   Consider  the  global  perturbation map $s\oplus ev:U\times B(r)\rightarrow {\bf L}$. Here we still use $ {\bf L}$ to denote the pull-back of the bundle over $U\times B(r).$  Clearly, by our construction,  $s\oplus ev$ is $C^{m_0}$-smooth when it is viewed in any  admissible local charts $U_i\times B(r), i\in I, $ and trivializations. 
   We will show that 
   
\medskip
\vspace{2mm}
\noindent
$\bullet$ $C1:$ when $r$ is small enough, $s\oplus ev:U\times B(r)\rightarrow {\bf L}$ 
  is  transversal to the zero section on each  $U_i\times B(r), i\in I $.
 
 Assume that $C1$ is true.  Let  ${\bf E}{\cal M}^{B(r)}$  be the  collection of the zero loci of the section inside $U_i\times B(r), i\in I $. Then it is a weakly $C^{m_0}$-smooth submanifold in ${\cal B}\times  B(r)$ in the sense defined before in this paper. Therefore, with the induced smooth structure, it is a $C^{m_0}$-manifold with expected dimension$=ind (Ds_f)+dim (B(r))$.   Let  $\pi: {\bf E}{\cal M}^{B(r)}\rightarrow B(r)$  be  the projection map between the  two finite dimensional manifolds, which is of class $C^{m_0}$. Then the "index" of $D\pi$ is equal to the $ind (s)$ which is fixed on each stratum.  Assume that $m_0> ind (s)$.  Then for generic choice of $\nu$ in $B(r)$,  perturbed moduli space ${\cal M}^{\nu} =(s^{\nu})^{-1}(0)$    is a
 $C^{m_0}$-submanifold of ${\bf E}{\cal M}^{B(r)}$ with dimension equal to $ind (s).$ Here 
 $\nu=\oplus_{i\in I} \nu_i $ in $B(r) , $  and 
   $s^{\nu}=s+\nu$  is the $\nu$-perturbed section of $s$  over   $U$. 
 
  To finish the construction, we need to show that ${\cal M}^{\nu}$ is compact for $r$  small enough.
 
To this end,   assume that $r<< C_0$.  Then  $s^{\nu}(X)=0 $ has no solution on $Bd(U)$. On the other hand,  since $s^{\nu}$ is continuous on ${\cal B}$,  any point in     $cl_{\cal B}({\cal M}^{\nu}) $  is still a solution of  $s^{\nu}(X)=0. $                          Since  $cl _{\cal B}(U) =U\cup Bd(U) $, this implies that $cl_{\cal B}({\cal M}^{\nu}) $ is inside $U$. Therefore, $cl_{\cal B}({\cal M}^{\nu})= {\cal M}^{\nu}.$

 Now ${\cal M}^{\nu}=\cup_{i\in I}{\cal M}^{\nu}_i$ and $cl_{\cal B}({\cal M}^{\nu})=\cup_{i\in I}cl_{\cal B}({\cal M}^{\nu}_i). $ We  already  proved that each $cl_{\cal B}({\cal M}^{\nu}_i)$  is compact in a more general setting  as a
  application of Picard method.  Therefore, ${\cal M}^{\nu}=cl_{\cal B}({\cal M}^{\nu})$ is compact for any $\nu$ in $B(r).$

  This finishes the construction of perturbed moduli space and proves the following theorem under the assumption that ${\cal B}$ has only one stratum and all isotropy groups are trivial.
 
 \begin{theorem}\label{second}
  When  $r$ is sufficiently small, for  any $\nu$ in $B(r)$, the perturbed moduli
  space ${\cal M}^{\nu} $ is  compact. For a generic choice of $\nu$,  ${\cal M}^{\nu} $ is a finite dimensional compact topological  manifold.  Moreover,  in the latter case, as a topological  submanifold  of ${\cal B}$ with induced weakly smooth structure with respect to the equivalent class of the covering data $[{\cal U}]$, 
${\cal M}^{\nu} $ is in fact an honest smooth manifold of class $C^{m_0}$.

\end{theorem}

 \medskip
\vspace{2mm}
\noindent $\bullet$ Proof of $C1:$
Fix a local slice $U_i\subset \subset W''_i\subset \subset W_i=W_{f_i}.$
Let $[s_i]\oplus [ev]:U_i\times B(r)\rightarrow L_{k-1}^{p}(\Sigma, \Lambda^{0,1}(f_i^*TM))$  be the  section $s\oplus ev:U\times B(r)\rightarrow {\bf L}$ written in the local trivialization over $U_i\times B(r).$ In the local trivialization, each
$\nu$ in $B(r)$ still takes a form $\nu=(\nu_1, \cdots, \nu_l)$, where $l$ is the cardinality of the index set $I.$
 But each $\nu_j$ with $j\not= i$ is obtained from the corresponding one composed with transition functions between the two slices and trivializations.  Since  there are only finitely many such $\nu_j$, for a fixed $i$, there exists a   constant $C_i>0$, such that   for  each such $\nu_j$,  $\|\nu_j(g)\|_{k-1, p}$ is  bounded by $C_i\cdot r$ for all $g\in U_i$.
 
 To prove $C1$, we only need to show the stronger statement that when $r$ is small enough,  $D([s_i]\oplus [ev])_{g, \lambda}:L_{k}^{p}(\Sigma, f_i^*TM, h_i)\times T_{\lambda}B(r)\rightarrow L_{k-1}^{p}(\Sigma, \Lambda^{0,1}(f_i^*TM))$  is surjective at any point $(g, \lambda)$ in  $U\times B(r).$  
 
We already calculated the this derivative and gave the related estimate
for quantitative transversality before. We recall the computation.

$$D([s_i]\oplus [ev])_{g, \lambda}(\xi, \eta)=D[s_i]_g(\xi)+[ev]_g(\eta)+({\partial}^{W_i}\lambda)_g(\xi). $$

Now it was prove that $$D[s_i]_g+[ev]_g:L_{k}^{p}(\Sigma, f_i^*TM, h_i)\times \oplus_{i\in I} GE(K_{f_i})\rightarrow L_{k-1}^{p}(\Sigma, \Lambda^{0,1}(f_i^*TM))$$  is surjective.  In fact, even $$D[s_i]_g+[ev_i]_g:L_{k}^{p}(\Sigma, f_i^*TM, h_i)\times  GE(K_{f_i})\rightarrow L_{k-1}^{p}(\Sigma, \Lambda^{0,1}(f_i^*TM))$$  is already surjective with a right inverse whose operator norm is bounded above for all $g$ in $U_i$. Therefore, there is a positive constant $\delta_i$ such that for any bounded linear operator $T_g:T_gU_f=L_{k}^{p}(\Sigma, f_i^*TM, h_i)\rightarrow L_{k-1}^{p}(\Sigma, \Lambda^{0,1}(f_i^*TM))$, if the operator norm $\|T_g\|<\delta_i, $ $$D[s_i]_g+T_g +[ev_i]_g:L_{k}^{p}(\Sigma, f_i^*TM, h_i)\times  GE(K_{f_i})\rightarrow L_{k-1}^{p}(\Sigma, \Lambda^{0,1}(f_i^*TM))$$  is  still  surjective 

  We already proved before  that when $r$ is small enough  for $\|\lambda\|<r,$ the operator norm of the partial derivative  $({\partial}^{W_i}\lambda)_g:L_{k}^{p}(\Sigma, f_i^*TM, h_i)\times  GE(K_{f_i})\rightarrow L_{k-1}^{p}(\Sigma, \Lambda^{0,1}(f_i^*TM))$   is less that the prescribed value $\delta_i$ for any $g$ in $U_i$.
 
 Put this together, we have proved that for sufficiently small $r$, $D([s_i]\oplus [ev])_{g, \lambda}$ is surjective for any ${g, \lambda}$ in $U_i\times B(r).$

 It remains to prove that the two assumptions can be arranged.

\medskip
\vspace{2mm}
\noindent
 ${\bullet}$ Proof of $A1$: 
 
 We assume that ${\cal M}$ is already covered by ${\cal U'}=\cup_{i\in I} U'_i, $ were $U'_i\subset\subset U_i$ is a smaller ball inside the ball $U_i$ with the same center $f_i.$  Here we have already assume that the radius of each $U_i$ is much smaller than the injective radius of $M$ so that each  $U_i$ is really  identified with  a  ball in  the model space, and $U'_i$ is a strictly smaller ball.

 By our assumption, ${\cal M}$ is contained in $\cup_{i\in I}  {\cal M'}_i, $ where $  {\cal M'}_i={\cal M}\cap U_i'.$ 
 To prove $A1$, we only need to show that the intersection of each ${\cal M'}_i$ with $Bd (U)$ is empty.
 
 Recall that by definition $Bd ( U)=cl_{\cal B}(U)\setminus U$. Then  $cl_{\cal B}(U)=\cup_{i\in I} cl_{\cal B}(U_i),$ and  each $cl_{\cal B}(U_i)$ is just the corresponding closed  ball.   Clearly,  $Bd ( U)=\cup_{i\in I} Bd_i $ where $ Bd_i=Bd (U_i)\setminus (\cup_{j\in I} U_j).$
 Now  consider ${\cal M'}_i$ with  $i$ being fixed.  Since  the intersection of ${\cal M'}_i$ with $Bd (U_i)$ is empty, so is the intersection of ${\cal M'}_i$ with $Bd_i. $ On the other hand, for any $j\not = i, $  $Bd_j$  has no intersection with $U_i$. Since
 ${\cal M'}_i$ is inside  $U'_i$. We still have that the intersection of ${\cal M'}_i$ with $Bd_j$ is empty. Therefore,  
 the intersection of each ${\cal M'}_i$ with $Bd (U)$ is empty.

 

 \medskip
\vspace{2mm}
\noindent 
 ${\bullet}$ Proof of $A2$:
  
A equivalent form of $A2$ is the statement  that if $\{g_k\}_{k=1}^{\infty}$  is a sequence in  $W'$ such that $\lim _{k\mapsto\infty}||s(g_k)||_{k-1, p} =0$, then there is a subsequence, still denoted by  $\{g_k\}_{k=1}^{\infty}$, such that it is $L_k^p$-convergent to an element in ${\cal M}.$
  
  Assume this is true, and  that $A2$ is not true. Then there is a sequence $\{g_k\}_{k=1}^{\infty}$ in $Bd(U)$  such that $\lim _{k\mapsto\infty}||s(g_k)||_{k-1, p} =0$. Then we have that $L_k^p$-limit $\lim _{k\mapsto\infty}g_k =g_{\infty}$ exits and is in ${\cal M}$
  
Since $Bd(U)$ is  closed in ${\cal B}={\cal B}_{k, p},$ $g_{\infty}$ is also in  $Bd(U)$.  Therefore, ${\cal M}$ and $Bd(U)$ has non-empty
intersection, which contradicts to $A1$.

  The  equivalent form of $A2$ for $s={\bar{\partial}}_J$  on a fixed stratum  was proved by Floer in [F1] based on local elliptic estimate for the non-linear ${\bar{\partial }}_J$-operator  near each point on the domain ${\Sigma}.$ The method there is the "standard" one. 
 
 Here we give a more "global"  proof based on the Taylor expansion for ${\bar{\partial }}_J$-operator used in Picard method.
  This proof works  for  general  Fredholm sections as long as those reasonable conditions specified in the proof are satisfied.

\medskip
\vspace{2mm}
\noindent 
 ${\bullet}$ Proof of the  equivalent statement of $A2$: 
 
  Since there only finitely many $W_i$'s, we only need to look at those $W'_i$ which contains infinitely many $g_k$'s. So we may assume
  that the sequence $\{g_k\}_{k=1}^{\infty}$ is contained in $W'_i\subset\subset W_i=W_{f_i, \epsilon_{f_i}}$ for a fixed $i.$ However, we need
   to assume that $\epsilon_i$  is  small enough so  that some conditions to be specified can be satisfied. 
   
   We remark that once ${\cal M}$ is given, for each $f$ in ${\cal M}$,  we need to choose these $\epsilon _f$'s satisfying  the  conditions  below  together with the conditions before.  We then define the corresponding $W_f$, $U_f$, etc. and select the finite covering of ${\cal M}$ by
   $W_i$ and  $U_i$ with $i\in I$ as before.

   Now we work with  the local slice, $W'_f\subset\subset W_f$ with $f$ in ${\cal M}$. Identify them with two open  balls of radius $\delta_2< \delta_1$ in $L_k^p(\Sigma, f^*(TM), h)$ with radius $\delta_2< \delta_1.$ Then   $f$ has local coordinate $f=0.$
 In the local trivialization of ${\bf L}$ over  $W_f$, we have the Taylor expansion at $f$ of $[s]=[s_f]$, the ${\bar{\partial}}_J$-section written in the trivialization as follows.
 For any $\xi$ in  $W_f$,  $$[s](\xi)=[s](0)+D[s]_f\xi +N(\xi)=D[s]_f\xi +N(\xi),$$ since $[s](0)=0.$
Moreover,  there is a constant $C$ depending on $f$, $\delta_2$ and geometric data of $M$ only such that 
$N$ satisfying
$$\|N(\xi)-N(\eta)\|_{k-1, p}\leq C(\|\xi\|_{k, p}+\|\eta\|_{k, p})\|\xi-\eta\|_{k, p}.$$

For our purpose, we introduce an equivalent metric on $E=L_k^p(\Sigma, f^*(TM), h).$  Let
${\bf K}=ker Ds_f$ and ${\bf K}^{\perp}$ be its $L^2$ orthogonal complement. 
 Since  ${\bf K}$ is finite dimensional,we have the decomposition $E={\bf K}\oplus{\bf K}^{\perp}$.

For any $\xi=(\gamma, \eta)$ in $E$, we define  a weaker but equivalent norm, $\|\xi\|^w_{k,p}=\|\gamma\|_{0,2}+\|\eta\|_{k,p}$.
By elliptic estimate for the linear elliptic operator $Ds_f$, we know that there exists a constant of similar nature, still denoted by $C$, such that $\|\gamma\|_{k,p}\leq C\|\gamma\|_{0,2}.$ Therefore, we get the same estimate for $N$ as above by replacing all $L_k^p$-norms on the right hand side by the new weaker norms and changing the constant accordingly.

 The upper semi-continuity of the dimension of the $ker Ds_f$ with respect to the  variable $f$ moving in ${\cal M}$ together the fact that
 ${\cal M}$ is compact, implies that the constants above  are bonded above, hence  independent of $f$ any more if we work on $W=\cup_{i\in I}W_i$.
 
 Go back to $W_f$. In local charts, the sequence in $W'_f$ has the form $g_k=(\gamma_k, \eta_k)$ with respect to the decomposition above.

 Note that  $\gamma_k$ is in ${\bf K}_{\delta_2}\subset\subset {\bf K}_{\delta_1}$. Here ${\bf K}_{\delta_2}$ and $ {\bf K}_{\delta_1}$
are the two discs in ${\bf K}$ of radii  ${\delta_2}$ and ${\delta_1}$ respectively measured in $L_k^p$-norm. Therefore, after taking a subsequence of  $\{g_k\}_{k=1}^{\infty}$, we may assume that $\{\gamma_k\}_{k=1}^{\infty}$ is already convergent in $L_k^p$ and hence $L^2$-norm.
In particular, $\|\gamma_k-\gamma_l\|_{0, 2}$ goes to zero when $k$ and $l$ go to infinity.

Now by elliptic estimate for  $Ds_f$, there is a (uniform) constant $C_1$ such that
$$\|\eta_k-\eta_l\|_{k, p}\leq C_1 \|D[s]_f(\eta_k-\eta_l)\|_{k-1,p}
=C_1 \|D[s]_f(\gamma_k,\eta_k)-D[s]_f(\gamma_l, \eta_l)\|_{k-1,p}$$ $$=C_1 \|D[s]_f(g_k)-D[s]_f(g_l)\|_{k-1,p}$$
$$ \leq C_1 \|N(g_k)-N(g_l)\|_{k-1, p}\leq C_1\cdot C (\|g_k\|_{k, p}+\|g_l\|_{k, p})\|g_k-g_l\|^w_{k, p}$$$$= C_1\cdot C (\|g_k\|_{k, p}+\|g_l\|_{k, p})(\|\eta_k-\eta_l\|_{k, p}+\|\gamma_k-\gamma_l\|_{0, 2}).$$

We now assume that  $\delta_2$ and $\delta_1$ are chosen in such a way that $2C_1\cdot C \delta_1\leq \frac {1}{2}.$ Since $g_k$ is in $W'_f$ which is a ball of radius equal to $\delta_2$, we have  $C_1\cdot C (\|g_k\|_{k, p}+\|g_l\|_{k, p})\leq \frac {1}{2}.$
Therefore, we have $$\|\eta_k-\eta_l\|_{k, p}\leq \frac {1}{2}\|\eta_k-\eta_l\|_{k, p}+\frac {1}{2}\|\gamma_k-\gamma_l\|_{0, 2}. $$
 
This implies that  $$\|\eta_k-\eta_l\|_{k, p}\leq \|\gamma_k-\gamma_l\|_{0, 2}.$$ Consequently,  after taking a subsequence, 
$\{g_k\}_{k=1}^{\infty}$ is a Cauchy sequence with respect to $L_k^p$ or its equivalent norm. Therefore  $\lim_{ k\mapsto\infty}g_k=g_{\infty}$ exits and  satisfies $s(g_{\infty})=0$ by continuity.




%





 \section{Smoothness of Banach Norm}

In this section we  collect some results related to the smoothness of $L_k^p$-norms.  Almost all the results are well-known and mentioned in various books on geometric analysis. Despite of the elementary nature of these results,   it seems that  not all of the  proofs are  widely known. 
 We  outline the key ideas of the proofs of these results here.

\begin{lemma}\label{first}
Let $\pi: V\rightarrow \Sigma$ be a $C^{\infty}$-smooth metric vector bundle on  a $C^{\infty}$-smooth Riemannian manifold $\Sigma,$  and $G$ be a Lie sub-group of the $C^{\infty}$- differomorphisms of $\Sigma.$ The action of G  on $\Sigma$ induces a action  on $ L_k^p(\Sigma, V)$.
Let  $\Psi: G\times L_k^p(\Sigma, V)\rightarrow L_k^p(\Sigma, V)$ be the induced action map. Denote  composition of $\Psi$  with the $p$-th power of $L_k^p$-norm by   $F:G\times L_k^p(\Sigma, V)
\rightarrow {\bf R}^1.$  Then $F$ is smooth if $p$ is an even non negative integer. Here we have assumed that $\Sigma$ is oriented and $G$-action preserves the orientation.

\end{lemma}

 \proof


The main issue here is the smoothness along $G$-direction.  In view of the lack of smoothness of the action map $\Psi$ along this direction,
It seems that this part of the lemma can not  be true.

To argue that this  is plausible, observe that  in the case that $\Sigma$ is ${\bf R}^n$ or ${\bf T}^n $ and  $G$ is the group of translations,  for a fixed section $\xi\in L_k^p(\Sigma, V), $ $F_{\xi}(g)=F(g, \xi):G\rightarrow {\bf R}^1$ is a constant map, hence trivially smooth.

This simple observation immediately suggests that the general case follows from the  changing variable formula  for integrations in calculus.
We only give the formula for the case that $k=1$. The case of  general k  can be proved similarly with more complicated notations.

We will write $F=F_0+F_1,$ and deal with each term separately.

Then, 
$$F_0(g, \xi)=\int_{\Sigma} ||g\cdot \xi||^p dx=\int_{\Sigma} ||\xi(g(x))||^p dx$$$$=\int_{g(\Sigma)} ||\xi(y)||^p det(g^{-1})dy=\int_{\Sigma} ||\xi(y)||^p det(J_{g^{-1}}(y))dy.$$ Here  $x$ and $y$ are two coordinate systems on $\Sigma$ such that 
$y=g(x)$ and $J_{g^{-1}}(y)$ is the Jacobian matrix of the differomorphism 
$g^{-1}:\Sigma\rightarrow  \Sigma.$

Clearly, the last identity show that $F_0$ is smooth in $g.$

$$F_1(g, \xi)=\int_{\Sigma} \Sigma_i||{\partial }_{x^i}(\xi(g(x)))||^p dx $$$$=\int_{\Sigma} \Sigma_i||\Sigma_j\frac {{\partial }y^j} {{\partial }x^i}(g^{-1}(y))
{\partial }_{y^j}(\xi(y))||^p det (J_{g^{-1}}(y)) dy$$ $$= \int_{\Sigma} \Sigma_i||\Sigma_j(J_g)^j_i(g^{-1}(y))
{\partial }_{y^j}(\xi(y))||^p det (J_{g^{-1}}(y)) dy$$

Here ${\partial }_{x^i}$ is the covariant derivatives on $V$ along the $x^i$-direction and $(J_g)^j_i$ is the $(j, i)$-entry of the Jacobian of the orientation preserving $C^{\infty}$ differomorphism $g:\Sigma \rightarrow \Sigma.$

Again, the last identity shows the smoothness of $F_1$ in $g$ when $p$ is an even non-negative integer.

Note that the role of above changing  variable formulas is to switch the $G$-action on $\xi$, which is only $L_1^p$ hence may "lose" derivative under $G$-action,  to the terms like $(J_g)^j_i(g^{-1}(y))$.  Since
$g$ is a $C^{\infty}$-differomorphism so that terms like $J_g$ have  $C^{\infty}$ smooth entries, it is well known that the lack of derivatives does not happen in this case.

The  $C^{\infty}$-smoothness along the $\xi$-direction  can be deduced  from the following well-known fact:

\begin{lemma}\label{second}
 The $p$-th power of the $L_k^p$-norm  on $L_k^p(\Sigma, V)$ is smooth for  even p. When  for   p is positive but not even,
it is at least $[p-1]$-smooth.

\end{lemma}


 \proof

Note that  the $L_k^p$-norm that we have used here   is $||\xi||_{k,p}=(\Sigma_{i\leq k}\int_{\Sigma} 
|\nabla^i \xi|^pdx)^{1/p},$  not the one obtained as the summation of the k semi-norms.
Clearly it is sufficient to  prove this only for $k=0.$ Assume that $p=2l$.
Then  $F(\xi)=\int_{\Sigma}< \xi, \xi>^ldx $   as a function from  $L_k^p(\Sigma, V)$ to ${\bf R}.$ If the derivative of $F$ exists at $\xi$, it has the form
$DF_{\xi}(\eta)=2l\int_{\Sigma}< \xi, \xi>^{l-1}<\xi, \eta>dx .$  One needs  to show that this formal derivative is the real one.
To proof this, one of key ingredients  is to show that it  is a bounded linearly map from  $L_k^p(\Sigma, V)$ to ${\bf R}.$
This follows form Horder inequality as following:

$$|DF_{\xi}(\eta)|\leq 2l\int_{\Sigma}| \xi|^{2(l-1)+1}|\eta|dx $$$$\leq  2l(\int_{\Sigma}| \xi|^{(2l-1)q} dx)^{1/q} \cdot (\int_{\Sigma}
|\eta|^pdx)^{1/p} $$$$=2l(\int_{\Sigma}| \xi|^{p} dx)^{1/q} \cdot (\int_{\Sigma}
|\eta|^pdx)^{1/p} $$$$ =2l|| \xi||_{p}^{p/q} \cdot ||\eta||_{p}. $$

Here we have used the identity $(p-1)q= p$ which follows from  $1/ p+1/q=1.$

\QED

There are several well-known corollaries of the above lemmas that we will use.

\begin{cor}\label{first}
Let $E =L_k^p(\Sigma, V)$ with p being  even and positive. There exists smooth cut-off function $\phi:E\rightarrow [0,\,\,1]$ supported in the unit ball of $E$ such that
$\phi(x)=1$ on the ball  of radius half.

\end{cor}

\proof

Let $\beta:{\bf R}\rightarrow [0,\,\,1]$ be a  "standard" smooth bump function supported on $[-1, \, \, 1]$ and equal to 1 on $[-1/2,\,\, 1/2].$ 
Then $\phi=\beta\circ F$ does the job, where $F$ is the p-th power of the $L_k^p$-norm.

\QED

\begin{cor}\label{first}
Let ${\cal B}$ be a paracompact Banach manifold of class $C^m$ modelled on  $E =L_k^p(\Sigma, V)$ with p  even and positive. Then ${\cal B}$
admits $C^m$ partition of unit.

\end{cor}

\proof

The proof in Lang's book for Hilbert manifold works equally well for this case. The key fact used in Lang's proof is that in Hilbertian case
the square of Hilbert norm is smooth.

\QED






\QED
\begin{lemma}\label{first}

 Let ${\cal B} =L_k^p(\Sigma, V)$ be the Banach space of all $_k^P$-maps from $\Sigma $ to $V={\bf R}^d$, $f:V\rightarrow {\bf R}^d$ is a $C^{\infty}$  function, regarded as a smooth vector field on $V$. Assume that the Sobolev  weight  $m_0 = k-\frac {n}{p} $ is large enough  so 
 that $ L_k^p(\Sigma, {\bf R^1})$ is a Banach algebra. Then the pull back of $f$, $\Psi_f:{\cal B} =L_k^p(\Sigma, V)\rightarrow L_k^p(\Sigma, {\bf R^d})$ is also smooth. Here $  \Psi_f(h)=f\circ h.$

\end{lemma}

The proof is similar to the one for the well-known fact that the induced map on the space of Sobolev maps under the coordinate transformation of the target space is   smooth. The detail of the proof will be given somewhere else.


\begin{thebibliography}{W}







\bibitem[C]{c}
 X. Chen,
Note on Partition of Unit,
{\it Electronic file, dated Feb. 8, 2013.}

%

%


%

%
%
\bibitem[F]{f}
A. Floer, Symplectic fixed points and holomorphic spheres,  
{\it  Commun. Math. Phys. } {\bf 120,} 575-611  (1989).

%


%

%
%
\bibitem[FH]{fh}
 H. Hofer,
A General Fredholm Theory and Applications,
{\it Preprint} (2005), arXiv: math.0509366[math.SG].
%


%
\bibitem[Ll]{Ll}
 S. Lang,
Differential Manifolds,
{\it Springer-Verlag} 1972.
%
%
\bibitem[L0]{l0}
G. Liu,      Associativity of Quantum  Multiplications,
{\it Comm. Math. Phys. 191 (1998),  265-282. }.

%
%
\bibitem[L1]{l1}
G. Liu,      Abstract Smoothness on Orbit Spaces,
{\it Electronic file,  dated Feb. 23, 2013}.

%

%
\bibitem[L2]{l2}
G. Liu,    Weakly Smoothness in GW Theory and  Smoothing of the Extended Moduli Spaces,
{\it In Preparation}.

%

%
\bibitem[L3]{l3}
G. Liu,   Flat  Structures in GW Theory and Smoothness of the Moduli Spaces,
{\it In Preparation}.

%
%
\bibitem[L4]{l4}
G. Liu,   Gluing  and Local Smoothness of ${\bar{\partial}}_J$-Section,
{\it In Preparation}.
%


%
%

%
%


%
%










%
%
\bibitem[LT]{lt}
G. Liu and G. Tian,    Floer homology and Arnold cojecture,
{\it J. Diff. Geom.} {\bf 49} No.2  (1998).

%

%
\bibitem[McW]{mcw4}
D. McDuff and  K. Wehrheim, Smooth Kuranishi     Structures with trivial isotropy, 
{\it  Preprint} 2012.
%

%

\end{thebibliography}
\end{document}